\documentclass[10pt,a4paper]{article}
\usepackage[utf8]{inputenc}
\usepackage{ucs}
\usepackage[english]{babel}
\usepackage{xcolor}
\usepackage[square,sort,comma,numbers]{natbib}			%
\usepackage{mathtools}
\usepackage{bm}
\usepackage{amsthm}
\usepackage{mathrsfs}
\usepackage{makecell, multirow, booktabs}
\usepackage{graphicx}
\usepackage{float}

\newtheorem{Remark}{Remark}
\newtheorem*{Remark*}{Remark} %
\newtheorem{Example}{Example}

\usepackage[createShortEnv]{proof-at-the-end}

\newcommand{\del}[1]{}

\newcommand{\myoption}{end}

\usepackage[left=2cm,right=2cm,top=2cm,bottom=2cm]{geometry}
\title{Optimal Linear Interpolation under Differential Information: application to the prediction of perfect flows}
\author{
\begin{tabular}{ccc} Soumyodeep~Mukhopadhyay$^{1,2,3}$ & Didier~Rulli\`ere$^{1,2}$ & Rodolphe~Le~Riche$^{2,1}$ \\
		Xavier~Bay$^{1,2}$ & Laurent~Genest$^{3}$ & David~Gaudrie$^{3}$
 \end{tabular}%
 }

\date{
\begin{tabular}{ccc}
    $^1$ \emph{Mines Saint-Étienne}%
    & & $^2$ \emph{CNRS, UMR 6158 LIMOS, FR} \\
    & $^3$ \emph{Stellantis, FR} &
\end{tabular}
}

\input{notation.sty}

\begin{document}
\maketitle
\noindent\rule{\linewidth}{0.5pt}
\abstract{%
\noindent Approximation of functions satisfying partial differential equations (PDEs) is paramount for %
\del{is fundamental to}%
simulation of physical fluid flows and other problems in physics. 
Recently, physics-informed machine learning \del{\cite{Karniadakis2021}} approaches have proven useful as a data-driven %
\del{alternative}%
complement to numerical models for partial differential equations, bringing faster responses and allowing us to capitalize on past observations%
\del{in mechanics}. %
\del{The data-driven surrogates leverage partial functional observations or simulation based design of experiments to approximate the full solution in a fraction of the typical computational time incurred by numerical models.}%
However, their efficiency and convergence depend on the availability of vast training datasets. %
\del{An alternate perspective is to look at this problem as a specialized interpolation in the presence of auxiliary differential information.}%
For sparse observations, Gaussian process regression or Kriging has emerged as a powerful interpolation model, offering principled estimates and uncertainty quantification. Several attempts have been made to condition Gaussian processes on linear PDEs via artificial or collocation observations and kernel design.%
These methods suffer from scalability issues in higher dimensions and %
limited generalizability. %
The aim of this study is to explore the extension of the Kriging predictor %
in the presence of %
linear PDE information at a finite number of collocation points. 
\del{Using a well-defined covariance relationship among a random process and its derivatives,}%
Two approaches are proposed%
: 1) A collocated co-Kriging \cite{Wackernagel2003} with primary observations of the physical field and auxiliary differential observations; 
2) A constrained Kriging optimization problem strongly satisfying linear PDE constraints at the points of prediction through a Lagrangian formulation. Numerical experiments are given for ordinary differential equations, 2D harmonic PDEs and an application to perfect flows around a cylinder.
This work highlights a trade-off between the computational efficiency of the Lagrange multipliers approach and the strict interpolation of observations.%
\\[2em]
\noindent{{\bf Keywords}: Physics-informed Machine Learning, %
\del{Constrained Kriging}%
Co-Kriging, Constrained Best Linear Unbiased Predictor}%
\newline
}
\noindent\rule{\linewidth}{0.5pt}

\pagebreak
\tableofcontents

\section{Introduction}

\newcommand{\citesth}{\textcolor{red}{[CITE REF(S)]~}}

Predicting fields governed by partial differential equations (PDEs) is of relevance in a wide range of scientific and industrial processes. Due to the lack of analytical resolutions, numerical methods are instrumental in practice. The existing literature on numerical approaches is extensive and focuses on stitching together solutions on spatial-temporal discretizations. Finite elements, finite volume and spectral methods are some of the most widely used numerical solvers. \cite{ames1977numerical} presents a detailed review on this subject. %
A general differential operator $\mathcal{A}$ \cite[see][]{Evans2010} 
 is typically defined as, 
\begin{eqnarray}
\mathcal{A}\left[ f\right](\vecx) &\coloneqq & g\left( D^\setM f, \vecx \right)~,   
\end{eqnarray} 
where $f: \setX \to \R$ is a smooth functional and $D^\setM$ represents mixed derivative orders, $\vecm$, that belong to a set of multi-indices $\setM$, i.e., $D^\setM \coloneqq \Big \lbrace \frac{ \partial^{\vert \vecm \vert} }{\partial \vecx_1^{m_1} \dots \partial \vecx_d^{m_d}}
~:~ 
 \vecm \in \setM \Big \rbrace$. %
The PDEs are linear if  $\mathcal{A}$ is linear. For the purposes of this paper\footnote{We will ignore the temporal aspect in our methodology but note that it is easily extendable using an anisotropic kernel \cite[$\S 4.2$]{Rasmussen2005}, where the time $t \in \R^+$ is an extra dimension in the input space and comes with additional cross-covariance requirements}, we limit ourselves to the boundary value problem formulation of PDEs,
\begin{eqnarray*}
\accolade{ \mathcal{A} \left[ f \right](\vecx) &=& v(\vecx) \quad \vecx \in \Omega~, \\
\mathcal{B} \left[ f \right](\vecx) &=& h(\vecx) \quad \vecx \in \partial\Omega ~,} 
\end{eqnarray*}
where $\mathcal{A}$ and $\mathcal{B}$ are differential operators applied on a compact domain $\Omega$ and its boundary $\partial \Omega$. %
\noindent Since they are based on a discretization, high-fidelity computational simulations are often expensive. Repeated computations over varying $\mathcal A$ or $\mathcal B$ quickly become impossible. Yet, it is a common task: for example, the design of an aerospace vehicle involves repeated calculations of the flow around it because the effect of some boundary conditions parameters is investigated \cite{airfoilbook,CADparams}. %
Data-driven approaches serve as surrogates to traditional numerical models that aim at faster computations. The forward problem of estimating $f$ given governing PDEs alongside true observations of $f$ or its derivatives, have garnered considerable interest recently. We present a brief review of such methods classified on the basis of the learning models used. 

\paragraph{Physics-informed neural networks (PINNs).}
The article \cite{Raissi2019} demonstrated that the data-driven approach using neural networks is viable for a variety of non-linear PDEs including the Burgers's equation and the 2D Navier-Stokes. 
The physical laws are injected as regularization terms for the neural network loss function. These additional terms depend on the residuals $\mathcal{A} \left[ f\right] - v$ and $\mathcal{B} \left[ f\right] - h$ computed at specific design points commonly referred to as collocation points. 
The authors first coined the term ``Physics-informed Neural Network'' (PINN) and, since then, data-driven approaches addressing PDEs and physical problems are referred to as physics-informed machine learning. 
There has been rapid advancement since the conception of PINNs in terms of neural network architectures, convolutional neural networks \cite{Gao2021}, generative adversarial networks \cite{Yang2020}, Bayesian neural networks \cite{Yang2021} among others. Other developments include neural operators such as \cite{kovachki2021neural}, a Fourier neural operator based on the fast Fourier transform \cite{Li2023} and DeepONet \cite{lu2021learning}, 
all of which focus on learning an operator mapping the forcing term and PDE parameters to the the PDE solutions. 
Some detailed reviews and benchmarks on the PINNs include \cite{Cuomo2022, PIMLreview, PINNtoPIKANsreview2024, hao2024pinnacle, PDEBENCH_An_Ex_Takamo_2022}. %
Despite its efficacy, the reliance on large training datasets is a challenge for PINNs. 
A direction taken in this regard is physics constrained neural networks as described in \cite{Sun2020, Liu2021} which use no or little labeled training data. The novelty is in the handling of initial and boundary conditions as hard constraints via special, output-layer encodings instead of the usual soft enforcement as regularization terms. 
Another direction, specifically aimed towards better uncertainty quantification, is to consider a Bayesian prior over unknown or noisy terms in a PDE presented as a PINN-type learning problem via Bayesian neural networks \cite{Yang2021}. To study convergence guarantees of PINNs, \cite{When_and_why_PI_Wang_2022} neural tangent kernel representation of neural networks attribute inaccuracy to spectral biases. 
\paragraph{Physics-informed conditional Gaussian processes.}
In the same spirit, there is a body of literature that uses Gaussian process (GP) priors \cite{Rasmussen2005} conditioned on PDE information at specific collocation points. GPs are well-established interpolation models and extensively used as surrogate models for black-box functions and they offer two major advantages over PINNs: 1) they perform better when observations are scarce; 2) by construction, they provide uncertainty estimates for their predictions. This motivated our choice of learning model. The application of GPs to approximating PDE solutions is made possible by the fact that GPs are closed under linear transformations. Hence, one can find closed form expressions for GPs under the action of linear differential operators. This inspired some of the early work in physics informed machine learning such as \cite{Skilling1992} for ordinary differential equations (ODEs). 
\cite{Solak2002} used this in the context of dynamical systems. 
\cite{Graepel2003234} presented possible applications to both ODEs and PDEs. Extensions such as \cite{Srkk2011} and  latent force models \cite{Alvarez2011} and \cite{Raissi2019} demonstrate the versatility of GPs.
There is also extensive research on constrained GPs \cite{DaVeiga2012, LpezLopera2018} which have considered constraints involving higher-order derivatives. These works involve conditioning the Gaussian process with linear observations or bounds of derivatives at collocation points. We adopt similar ideas as derivative observations or artificial observations or virtual observations as termed by \cite{Jidling2017}. Naturally, this drives up the computational cost given the $\mathcal{O}(n^3)$ cost of dense matrix inversion, where $n$ is the number of observations. However, in \cite{High_Dimensiona_de_Roo_2021} the Gram matrix structure of the covariance matrix is exploited to drive down this cost significantly. Despite ongoing efforts, computational cost remains a substantial challenge and merits further investigation. 
\paragraph{Physics-informed kernel design.}
Another line of work is about designing matrix-valued covariance kernels that can satisfy PDEs. 
For example, the divergence free kernel proposed in \cite{Narkowich1994} was applied to fluid dynamics in \cite{Divergence_Free_Wendla_2009}, and a curl-free kernel is described in \cite{scheuerer2012covariance}. 
A review article on multi-output kernels provided a unified view on such kernels \cite[$\S 5$]{Alvarez2012}. 
More recently, \cite{Jidling2017} proposed linearly constrained GPs and reduced it to computing the inverse of a linear differential operator. 
\cite{Physics_derived_Cross_2021} developed a physics-informed covariance kernel for applications in structural dynamics specifically addressing the forced linear oscillator equation. 
As shown in \cite{henderson2023}, if a kernel satisfies a specific linear PDE, it induces a second order random field, i.e. a Gaussian process, whose trajectories satisfy this PDE almost surely.
This result was later applied to the 3D-wave equation \cite{henderson2023application}.
Although it is of course preferable to build a model from adapted components, constructing kernels that satisfy PDEs is problem specific. Furthermore, such kernels are applicable to linear PDEs only. 
A possible path towards more general model constructions was proposed in \cite{Wang2022PIDKL} as learning kernel functions via deep learning. 
\newline
\newline
Building on the idea of physics-informed conditional Gaussian processes, we make use of the optimization perspective offered by best linear unbiased prediction (BLUP), commonly referred to as Kriging \cite{Wackernagel2003}. The general contribution of this work is to extend best linear predictors to make surrogates of PDE solutions. 
We first demonstrate the extension of collocated co-Kriging \cite{Wackernagel2003} to pass differential information as secondary observations. \cite[$\S 5.5$]{Chils2012} briefly hints at the possibility of such a model. Second, we formulate a constrained Kriging predictor considering differential constraints at points of prediction. 
The utility and limitations of these two BLUP approaches are compared through the prediction quality and the computational cost incurred. 
\newline \newline
\noindent This paper begins, in Section~\ref{sec: krigingintro}, with an overview of Kriging as BLUP and a modified form of ordinary Kriging relevant to our application. It moves on to a brief description of linear differential operators and associated notations in Sections~\ref{sec: linops} and \ref{sec: extds}.
In particular, we introduce the extended design space that has derivative indices as an extra dimension of the original spatial domain. 
We describe the two proposed methodologies in Sections~\ref{sec: derivativeobs} and \ref{sec: ConstrainedKriging}. 
In Section~\ref{sec: comparison}, we examine the theoretical and empirical differences between the two methods. Further numerical results on synthetic datasets and an application to the prediction of perfect flows in 2D are given in Section~\ref{sec: results}. A summary of the work along with perspectives, notably for nonlinear PDEs, can be found in Section~\ref{sec: conclusion}.

\section{Background}

Approximating the best predictor of a quantity of interest given some data about other related quantities, lies at the heart of modern machine learning techniques. The most widely used loss function or metric to quantify the prediction quality in this context is the mean squared error (MSE). Specifically, this loss function is the expected squared error incurred by a certain predictor function of the given data. The frequentist definition of MSE is simply:
\begin{equation}
MSE(Y, \hat{Y}) = \Esp{\left(Y - \hat{Y}\right)^2}
\end{equation}
where $Y$ is a real-valued $L^2$ random variable %
and $\hat{Y}$ is an estimator of $Y$. Given partial information in the form of a vector of real-valued $L^2$ random variables $\mathbf{X} = (X_1, \ldots, X_n)$ %
we intend to minimize,
\begin{equation}
\underset{\varphi(\mathbf{X}) \in L^2} {\min} \Esp{\left(Y - \varphi(\mathbf{X})\right)^2}
\end{equation}
It is well known \citep{Blitzstein2014} that the conditional expectation $\tau(\mathbf{X}) = \Esp{Y\vert\mathbf{X}}$ is the unique minimizer of the mean squared error given $\mathbf{X}$. This follows from the orthogonal projection with respect to the inner porduct 
$\langle \varphi_1 , \varphi_2 \rangle \coloneqq \Esp{\varphi_1 \varphi_2}$. In fact, the residual $Y-\Esp{Y\vert\mathbf{X}}$ is orthogonal to any function $\varphi(\mathbf{X}) \in L^2$ with respect to this inner product. Furthermore, the minimal mean squared error attained, with the conditional expectation as estimator, is the expected conditional variance $\Esp{\Var{Y\vert\mathbf{X}}}$.
Most machine learning tasks can be interpreted as approximating the conditional expectation, for example, Neural networks parameterize the problem with $\tau_{\theta}(\mathbf{X})$ where $\theta$ are the network parameters and optimize the empirical MSE or loss functions that scale monotonically with respect to the MSE. In this Section, we discuss standard results when restricted to the class of linear predictors based on $\mathbf{X}$. We will also discuss some mathematical background to tackle differential equations in this regard. 

\subsection{Kriging as the Best Linear Unbiased Predictor}
\label{sec: krigingintro}
Originally developed in the field of geostatistics, Kriging \cite{Wackernagel2003, Chils2012} is an interpolation model for a real-valued, %
second order random field $(Y(\vecx))_{\vecx\in \setX}$, where $\setX \subset \R^d$ denotes the spatial domain %
or design space. It is commonly introduced as a generalization of linear regression, providing the best linear unbiased predictor (BLUP) \cite{Kolmogorov1941, weiner1949} assuming a known mean and by incorporating spatial correlation through a specified covariance structure. %
The linear predictor for an unknown location $\vecx\etoile$, given $n$ locations %
$\vecX = (\vecx_i)_{i=1}^n$ with observations, $\vecY = Y(\vecX)$,
\begin{eqnarray}
\label{eqn: linearestimator}
\widehat{Y}(\vecx\etoile) &\coloneqq & \transpose{\vecalpha(\vecx\etoile)} \vecY    
\end{eqnarray}
where $\vecalpha(\vecx\etoile) \in \mathbb{R}^n$. An unbiased predictor, such that $\EspSymbol\left[\widehat{Y}(\vecx\etoile)\right] = \EspSymbol\left[ Y(\vecx\etoile)\right]$, should satisfy,
\begin{eqnarray}
\label{eqn: unbiasedness}
\transpose{\vecalpha(\vecx\etoile)} \vecmu &=& \mu\etoile
\end{eqnarray}   
where $\vecmu = \EspSymbol \left[\vecY\right]$ and $\mu\etoile = \EspSymbol \left[ Y(\vecx\etoile)\right]$. The linear coefficients, $\vecalpha(\vecx\etoile)$, are estimated as the minimizer of the error variance defined below,
\begin{eqnarray}
\label{eqn: minerrorvar}
\vecDelta(\vecx\etoile) &\coloneqq & \EspSymbol\left[ \left( \widehat{Y}(\vecx\etoile) - Y(\vecx\etoile)\right)^2 \right]~, 
\end{eqnarray}
allowing us to set up the final constrained optimization problem, 
\begin{eqnarray}
\label{eqn: krigingopt}
\accolade{&\underset{\vecalpha(\vecx\etoile) \in \mathbb{R}^n}{\argmin}& \vecDelta(\vecx\etoile)  \\
&\text{subject to: }& \transpose{\vecalpha(\vecx\etoile)} \vecmu  = \mu\etoile }
\end{eqnarray}
Extending to the case of predicting at multiple locations $\vecX\etoile = (\vecx\etoile_j)_{j=1}^q$, would imply a multi-objective optimization problem, equation \eqref{eqn: krigingopt} for each point. Fortunately, we can simplify the problem by considering the sum of the expected squared errors, $\vecDelta(\vecx\etoile_1) + \dots + \vecDelta(\vecx\etoile_q)$, as variance is non-negative. The optimization problem, in this case, writes out as, 
\begin{eqnarray}
\label{eqn: krigingoptmultiple}
\accolade{ &\underset{\matAlpha \in \mathbb{R}^{n \times q}}{\argmin}& \left(\sum_{j=1}^q \vecDelta(\vecx\etoile_j)\right) \\
&\text{subject to: \quad}& \transpose{\matAlpha} \vecmu = \vecmu\etoile~,} 
\end{eqnarray}
where $\matAlpha = \transpose{\left(\vecalpha(\vecx\etoile_1), \dots, \vecalpha(\vecx\etoile_q)\right)} \in \mathbb{R}^{n \times q}$ and $\vecmu\etoile = \EspSymbol\left[ Y(\vecX\etoile)\right]$. To solve this optimization problem, we assume the covariance structure of $Y$ in terms of a positive semi-definite (PSD) kernel $k: \setX \times \setX \to \mathbb{R}$, 
\begin{eqnarray}
\CovSymbol\left[ Y(\vecx), Y(\vecx')\right] &=& k(\vecx, \vecx') 
\end{eqnarray}
We can rewrite the objective function in equation \eqref{eqn: krigingoptmultiple} (see proof of proposition \ref{prop: OrdinaryKriging} for details) as, 
\begin{eqnarray}
\label{eqn: krigingopttraceform}
\accolade{&\underset{\matAlpha\in \mathbb{R}^{n\times q}}{\argmin}& \Trace{\transpose{\matAlpha}\matK \matAlpha} - 2 \Trace{\transpose{\matAlpha} \matH} + \Trace{\matK\etoile} \\
&\text{subject to: \quad}& \transpose{\matAlpha} \vecmu = \vecmu\etoile }
\end{eqnarray}
We solve equation \eqref{eqn: krigingopttraceform} by introducing a vector of Lagrange multipliers $\veclambda$ of size $1 \times q$. This is similar to the computation of ordinary Kriging expressions. The final expressions are presented in the following proposition. 
\begin{thmE}[Ordinary Kriging without constraint][\myoption]
\label{prop: OrdinaryKriging}

Let $\matAlpha=\cbind{\vecalpha(\vecx\etoile_1), ..., \vecalpha(\vecx\etoile_q)}$ the $n \times q$ matrix of the linear predictor's optimal weights minimizing \eqref{eqn: krigingopttraceform} at all prediction locations $\vecx\etoile_1, \ldots, \vecx\etoile_q$, under the unbiasedness constraint. Then
\begin{equation}\label{eqn: ordinarykriging}
\matAlpha = \matK^{-1} \left( \matH +   \vecmu \transpose{\veclambda}
\right) \,,
\end{equation}
with $q \times 1$ Lagrange multiplier $\veclambda$
\begin{equation*}
\transpose{\veclambda} = \left( \transpose{\vecmu} \matK^{-1} \vecmu \right)^{-1}
\left( \transpose{\vecmu\etoile} - \transpose{\vecmu} \matK^{-1} \matH \right) \, ,
\end{equation*}
Covariance matrices are $\matK \in \R^{n \times n}$, with 
 $\matK_{ij} = \Cov{Y(\vecx_i), Y(\vecx_j)}$, $i, j \in \set{1, \ldots, n}$,  and  \newline $\matH=\cbind{\vech(\vecx\etoile_1),\ldots, \vech(\vecx\etoile_q)} \in \mathbb{R}^{n \times q}$, with $\vech_i(\vecx\etoile) = \Cov{Y(\vecx_i), Y(\vecx\etoile)}$.
\end{thmE}
\begin{proofE}
~Consider a prediction location $\vecx\etoile$. The weights $\vecalpha(\vecx\etoile)$ are depending on $\vecx\etoile$, but for the sake of simplicity we write $\vecalpha$ when this single prediction location is considered. Using Equation~\eqref{eqn: krigingopttraceform}, one easily shows that 
$$\Delta(\vecx\etoile) =  \transpose{\vecalpha} \matK \vecalpha - 2 \transpose{\vecalpha}\vech(\vecx\etoile) + \sigma^2(\vecx\etoile) $$
where $\matK_{ij} = \Cov{Y(\vecx_i), Y(\vecx_j)}$, $\vech_i(\vecx\etoile) = \Cov{Y(\vecx_i), Y(\vecx\etoile)}$ , $i, j \in \set{1, \ldots, n}$, and $\sigma^2(\vecx\etoile)=\Var{Y(\vecx\etoile)}$.  

The minimization under the constraints %
leads to minimize
$$\Delta_\lambda(\vecx\etoile) =  \transpose{\vecalpha} \matK \vecalpha - 2 \vech(\vecx\etoile) -2 \lambda \left( \transpose{\vecalpha} \vecmu - \mu\etoile\right) \, .$$
Taking the gradient with respect to $\vecalpha$, one gets

$$ \matK \vecalpha -  \vech(\vecx\etoile) - \lambda  \vecmu =0 \, .$$
So that finally

$$ \vecalpha = \matK^{-1} \left( \vech(\vecx\etoile) + \lambda  \vecmu \right) \, ,$$
 
and as $\transpose{\vecalpha}\vecmu = \mu\etoile$, one gets

$$
\lambda = \left( \transpose{\vecmu} \matK^{-1} \vecmu \right)^{-1}
 \left( \mu\etoile - \transpose{\vecmu} \matK^{-1} \vech(\vecx\etoile) \right)  \, .$$

Hence the optimal weights $\vecalpha$ (at the specific location $\vecx\etoile$) for the prediction $\widehat{Y}(\vecx\etoile)$ (Kriging mean), and its quadratic error $\Delta(\vecx\etoile)$ (Kriging variance).\\

Remark that both $\vecalpha$ and $\lambda$ are depending on $\vecx\etoile$ in the previous equations, and should have been written $\vecalpha(\vecx\etoile)$ and $\lambda(\vecx\etoile)$.
One can gather predictions at $q$ several points in a vector. This allows to do all the predictions at once by a single matrix expression for a $n \times q$ matrix of weights $\matAlpha$:

\begin{equation}\label{eqn: ordinarykrigingproof}
\matAlpha = \matK^{-1} \left( \matH +   \vecmu \transpose{\veclambda}
\right) \,,
\end{equation}

where $\matAlpha=\cbind{\vecalpha(\vecx\etoile_1), ..., \vecalpha(\vecx\etoile_q)}$, $\matH=\cbind{\vech(\vecx\etoile_1), ..., \vech(\vecx\etoile_q)}$, and where $\veclambda =({\lambda(\vecx\etoile_1), ..., \lambda(\vecx\etoile_q)})$.
\end{proofE}
In the case when the underlying process is Gaussian, the Kriging minimizer and minimum associated error correspond to the Gaussian conditional expectation and conditional variance expressions. It is one case where the best linear predictor is also the best predictor.
\subsection{Linear differential operator equations}
\label{sec: linops}
  
A linear partial differential operator $\mathcal{L}$ is generally expressed as,
\begin{eqnarray}
\label{eqn: linop form}
\LinOp{f}(\vecx) &=& \underset{\vecm \in \setM }{\sum} u_{\vecm}(\vecx) D^{\vecm} f (\vecx)~,   
\end{eqnarray}     
where $u_{\vecm}(\vecx) \in \R$ are real coefficients for the mixed derivatives of $f$ and $\setM$ is a set of relevant multi-indices. A linear differential operator equation enforced on a discrete set of $p$ locations $\vecX = \big \lbrace \vecx_i \in \setX \big \rbrace_{i=1}^p$, 
\begin{eqnarray}
\label{eqn: linopeqn}
\underset{\vecm \in \setM }{\sum} u_{\vecm}(\vecx_i) D^{\vecm} f~(\vecx_i) &=& v(\vecx_i)\quad \text{for } i = 1,\dots,p
\end{eqnarray} 
where $v: \setX \to \R$ is some other functional, also referred to as the forcing term in the PDE literature, can be expressed as vector equations,
\begin{eqnarray}
\transpose{\vecu_i} D^\setM f~(\vecx_i) &=& v(\vecx_i)  \quad \text{for } i = 1,\dots,p
\end{eqnarray}
where $\vecu_i \coloneqq \lbrace u_{\vecm}(\vecx_i)\rbrace_{\vecm \in \setM}$. Arranging the coefficient vectors $\vecu_i$ as columns of a matrix, 
\begin{eqnarray}
\transpose{\matU} \mathrm{vec} \left(D^\setM f~(\vecX)\right) ~=~ v(\vecX)~;~ \matU~=~\left[ \vecu_1, \dots, \vecu_p \right]~,  
\end{eqnarray} 
where $\mathrm{vec}~(\cdot)$ is the vectorization operation \cite{Loan2000}, stacking the columns of the matrix to create a single vector, and $D^\setM f~(\vecX) = \left \{  D^{\vecm}f~(\vecx); ~\vecm \in \setM, \vecx \in \vecX \right \}$ is a $\m \times p;~\m = \card{D^\setM} f, p = \card{\vecX}$, sized matrix. Each row represents a specific mixed derivative $\vecm$, applied over the set of points $\vecX$ and each column represents all mixed derivatives applied to a specific point. %
The matrix $\matU$ can be made to represent a variety of information such as the average estimator, multiple PDEs at several points and average of PDEs over several points. We present several examples to elucidate the construction of the $\matU$ matrix. 

\begin{Example}[One PDE satisfied over two points]
A simple example to illustrate the coding of a system of differential equations defined pointwise. Consider the linear differential equation: $f(\vecx) + f''(\vecx) = 0$. The matrix encoding accounting for this information at two spatial points $\vecx_1, \vecx_2 \in \setX$,  
\begin{eqnarray}
\label{eqn: example1}
\begin{pmatrix}
1 & 0 & 1 & 0 & 0 & 0 \\
0 & 0 & 0 & 1 & 0 & 1 
\end{pmatrix} \begin{pmatrix}
f(\vecx_1)\\
f'(\vecx_1)\\
f''(\vecx_1)\\
f(\vecx_2)\\
f'(\vecx_2)\\
f''(\vecx_2)
\end{pmatrix} &=& \bm{0}~,
\end{eqnarray}
where $\transpose{\matU}$ is the $2 \times 6$ left matrix, $vec\left(D^\setM f~(\vecX)\right)$ is the $6 \times 1$ right matrix, and $\bm{0}$ is the corresponding $v(\vecX)$.  
\end{Example}
\begin{Example}[Two PDEs satisfied at one point]
Alternatively, suppose that a single spatial point $\vecx \in \setX$ verifies another differential equation, $2f'(x) - f''(x) = 0$, in addition to the one formerly described. The matrix representation for the system of differential equations at $\vecx$ would be,  
\begin{eqnarray}
\begin{pmatrix}
1 & 0 & 1 \\
0 & 2 & -1  
\end{pmatrix} \begin{pmatrix}
f(\vecx)\\
f'(\vecx)\\
f''(\vecx)\\
\end{pmatrix} &=& \bm{0}
\end{eqnarray}
\end{Example}

\begin{Example}[One PDE satisfied in average]
An interesting aspect is to consider the sample average \cite{Rulliere2025}. Suppose it is desirable to impose conditions on the sample mean of our predictions, possibly over several locations and consisting of mixed orders of derivatives. This is feasible with only minor modifications to the $\matU$ matrix. For example, a PDE satisfied in average over a set of locations $\{\vecx_i: i = 1,\ldots,q\}$ implying $\frac{1}{q}\sum_{i=1}^q f(\vecx_i) + f''(\vecx_i) = \bm{0}$ can be coded using an extension of the $\transpose{\matU}$ in eq. \eqref{eqn: example1}, replacing $1$'s with $1/q$ and extending the matrix to $q$ rows. Similarly, extending the vector of the function derivatives at spatial locations to $3q$ rows using the same arrangement as seen in eq. \eqref{eqn: example1}.
\end{Example}

\noindent The description above, applies to a sufficiently differentiable random field in the mean-squared differentiability sense \cite{Christakos1992}. In the context of a Kriging model $Y$, we work towards developing the expression for the BLUP predictor given auxiliary information such as, 
\begin{eqnarray}
\label{eqn: constraints1}
\transpose{\matU} \mathrm{vec} \left(D^\setM Y~(\vecX)\right) &=& v(\vecX)
\end{eqnarray} 
This is well-defined if $\Esp{Y} = \vecmu$ is differentiable and the variance exists for all considered derivative orders of $Y$.
 
\subsection{Extended design space}
\label{sec: extds}

This Section explores the handling of spatial locations and orders of derivatives of the random function $Y$ with a covariance structure made explicit by the kernel function $k_\setX: \setX \times \setX \to \R$. To formalize the spatial coordinates and the order of derivative as an integrated mathematical object, we introduce an extended design space defined as the Cartesian product $\setS \coloneqq \setX \times \setM$, where $\setX$ denotes the original spatial domain and $\setM$ represents the set of relevant multi-indices. This offers a unique viewpoint, derivatives of $Y$ are points in the extended design space. %
Consider the extended random function $Z: \setS \to \R$ which maps, 
\begin{eqnarray}
Z&:& (\vecx, \vecm) \mapsto Y^{(\vecm)}(\vecx)~, 
\end{eqnarray}
where the vector $(\vecx, \vecm) \in \setS$ will be referred to as $\vecs$ henceforth. 
Equation \eqref{eqn: constraints1} re-writes as a compact matrix equation, 
\begin{eqnarray}
\label{eqn: constraints2}
\transpose{\matU} \vecZ &=& \vecv
\end{eqnarray} 

To fully model $Z$ as random function, we need to quantify the covariance. Assuming $\Esp{Y}$ is sufficiently differentiable, we seek a modified kernel of the form $k_\setS: \setS \times \setS \to \R$ to quantify the covariance between $Z(\vecs) \coloneqq Y^{(\vecm)}(\vecx)$ and $Z(\vecs') \coloneqq Y^{(\vecm')}(\vecx')$ where $\vecs = (\vecx, \vecm)$ and $\vecs' = (\vecx', \vecm')$. The following result on covariance functions \cite{Abrahamsen1997, Christakos1992, Rasmussen2005}, facilitates exactly that.
\begin{eqnarray}
\label{eqn: CovAbrahamsen}
\Cov{Z(\vecs), Z(\vecs')} &=& \frac{\partial^{|\vecm|}}{\partial x_1^{\compm_1} \cdots \partial x_d^{\compm_d}} \frac{\partial^{|\vecm'|}}{\partial x'^{\compm'_1}_1 \cdots \partial x'^{\compm'_d}_d}
 k_\setX~(\vecx,\vecx')~,
\end{eqnarray}
where we require that $\displaystyle \partial^{|\vecm| + |\vecm'|}k_\setX~(\vecx, \vecx')$ exist for all relevant $\vecs, \vecs' \in \setS$ to ensure that all expressions that follow are well-defined. For sufficiently regular kernels, such as the squared exponential kernel (smooth), this holds true. We provide the analytical derivative expressions in the apppendix \ref{sec: kernelderiavtives}. This is a special covariance relationship since it involves kernel derivatives of arbitrary order but this is well-explored in the GP literature since the GPs are closed under the action of linear differential operators \cite{Graepel2003234, Srkk2011, DaVeiga2012, Raissi2017, Jidling2017, Solak2002}. 
\newline \newline
\noindent Notational convenience aside, the extended design space allows for the direct application of the usual Kriging formulas, such as Proposition \ref{prop: OrdinaryKriging}, or any other covariance based prediction model, without any adaptation to derivatives. We highlight that the extended design space characterization and the unique perspective it offers, to the best of our knowledge, is an original contribution. An example of how this can prove to be a powerful tool when dealing with derivatives in presented below.
\begin{Example}[Kriging in the extended design space]
In a one-dimensional spatial domain, $\setX = \R$, consider the following problem: predict the derivative $Y'(x\etoile)$ from observations $Y(x_1)$ and the second-order derivative $Y''(x_2)$. The extended design space considerably simplifies the problem since we can directly apply the classical Kriging formulas (Proposition \ref{prop: OrdinaryKriging}) given observations $\vecZ = \transpose{(Z(\vecs_1), Z(\vecs_2))}$ where $\vecs_1 = (x_1, 0)$ and $\vecs_2 = (x_2, 2)$ to predict $Z(\vecs\etoile)$ at $\vecs\etoile = (x\etoile, 1)$. The necessary covariance matrices, $\Cov{Z(\vecs_i), Z(\vecs_j)}, i,j \in \{1,2\}$ and $\Cov{Z(\vecs_i), Z(\vecs\etoile)}, i \in \{1, 2\}$ can be constructed using \eqref{eqn: CovAbrahamsen}.  
\end{Example}

\noindent  To summarize, the  extended design space allows us to work with a process $Z$, defined on a modified design space that intrinsically handles the involved derivative operations and therefore, encode differential equations as simple linear combinations of $Z$.

\section{Methodology}

This Section addresses the two proposed treatments of pointwise differential information. The first method, an adaptation of the collocated co-Kriging methodology and the second method presenting a Lagrangian based constraint resolution at the points of prediction.

\subsection{Collocated co-Kriging}
\label{sec: derivativeobs}
Several attempts have been made at using derivative information as observations in the context of the usual Gaussian process regression methodology. The earliest related work traces back to applications in Ordinary Differential Equations \cite{Skilling1992} and dynamical systems \cite{Solak2002}. In the past, multiple studies have explored the use of GPs in problems involving linear differential operators, specifically \cite{Graepel2003234, Srkk2011} which account for derivative information as virtual observations. This has inspired a large variety of work in the recent past for applications in solving linear PDEs such as, \cite{Raissi2017} and \cite[\S 3.1]{Jidling2017}. The latter refers to this approach as artificial observations and highlights the computational challenge in higher dimensions. Constrained GPs involving derivative, inequality constraints enforced at finite points, \cite{DaVeiga2012} is a prominent example. For a detailed review on this subject, the interested reader may refer to \cite[\S 6.1]{Swiler2020}.   
\newline \newline  
\noindent We utilize a collocated co-Kriging \cite{Wackernagel2003} methodology to obtain the best linear pointwise predictor given primary observations and derivative information as secondary observations. A very similar approach has been discussed in \cite[\S 5.5]{Chils2012}. A notable difference is that, \cite{Chils2012} mostly addressed the use of direct derivative observations as secondary information whereas we construct a generic co-Kriging predictor composed of any linear transformation of mixed derivatives. In our framework, co-Kriging with direct derivative observations reduces to simple or ordinary Kriging (Proposition \ref{prop: OrdinaryKriging}) with observations in the extended design space. However, \cite{Chils2012} also proposed the use of finite differences approximation to account for linear Neumann boundary conditions via co-Kriging, the implications of which we did not explore yet. We dedicate the rest of this Section to explain the mathematical details of co-Kriging with pointwise linear PDE information as secondary variables.    
\newline\newline
\noindent Suppose we are given $n$ primary observations of $Z$ at locations $\vecs_1,\ldots,\vecs_n$ and differential equations imposed at $p$ collocation points in $\setX$ (similar to physics-informed neural networks \cite{Raissi2019}) consisting of $c$ locations as part of the PDEs
$\vecs^+_1,\ldots,\vecs^+_{c}$ where the superscript `$+$' denotes collocation nature. For example, consider the differential equation $\partial^2_{xx} Y(\vecx^+) + \partial^2_{yy} Y(\vecx^+) = 0$ along with an observation of $Y$ at point $\vecx$. Here, $n=1$ (observation $\vecx$), $p=1$ (differential equation at one collocation point $\vecx^+$) and $c=2$, considering points $\vecs^+_1 = (\vecx^+, \partial^2_{xx})$ and $\vecs^+_2 = (\vecx^+, \partial^2_{yy})$ in the extended design space $\setS$. Even though $c$ does not show up in the final expressions, it necessary to define the co-Kriging predictor in eq-\eqref{eqn:KrigingPredictorDerivativeObs}.
The Kriging predictor at a point $\vecs\etoile$ involving both the primary observations $\vecZ = \transpose{\left( Z(\vecs_1),\ldots,Z(\vecs_n)\right)}$ and the secondary observations $\transpose{\matU} \vecZ^+ = \vecv$ where $\vecZ^+ = \transpose{\left(Z(\vecs^+_1),\ldots,Z(\vecs^+_c)\right)}$, is given as,
\begin{eqnarray}
\label{eqn:KrigingPredictorDerivativeObs}
\Zck(\vecs\etoile) &=& \sum_{i=1}^n \alpha_i(\vecs\etoile) ~ Z(\vecs_i) + \sum_{j=1}^p \alpha^+_j(\vecs\etoile)~\transpose{u}_j~\vecZ^+ ~. 
\end{eqnarray}
where $\alpha^+$ is to distinguish coefficients related to the differential equations and $u_j$ is the $j^{th}$ column of $\matU$. This predictor is the same as the multicollocated ordinary co-Kriging predictor in described in \cite{Chils2012} except that we use  an arbitrary set of collocation points, leading to an optimization problem analogous to the one described in \eqref{eqn: krigingopt}. Considering both the initial vector $\vecZ$ and the auxiliary vector $\transpose{\matU} \vecZ^+$ as observations, for $q$ points of prediction, $\vecs_1\etoile,\ldots,\vecs_q\etoile$ we solve for all linear coefficients, $\matAlpha^+ = \left[ \vecalpha(\vecs_1\etoile),\ldots,\vecalpha(\vecs_q\etoile)\right]$ as we did in proposition \ref{prop: OrdinaryKriging}. We can obtain zero error variance at some observation location $\vecs_i$ by setting linear coefficients $\vecalpha = e_i = (0,\ldots,1,\ldots,0)$, and hence satisfying the interpolation property $\Zck(\vecs_i) = Z_i$. We summarize the final expressions in the following proposition.

\begin{thmE}[Collocated ordinary co-Kriging][\myoption]
\label{prop: collocatedcokriging}
The coefficient matrix associated with the co-Kriging linear predictor, minimizing the predicted error variance under the unbiasedness constraint~\eqref{eqn: krigingopttraceform} is given by
\begin{equation}
\label{eqn: OrdinaryCollocatedcoKriging}
\matAlpha^+ =  \left(\matK^+\right)^{-1} \left(\matH^+ + \vecmu \transpose{\veclambda} \right)~, 
\end{equation}	
with $q \times 1$ Lagrange multiplier
\begin{equation*}
\transpose{\veclambda} = \left( \transpose{\vecmu} (\matK^+)^{-1} \vecmu \right)^{-1} \left( \transpose{\vecmu\etoile} - \transpose{\vecmu} (\matK^+)^{-1} \matH^+ \right)
\end{equation*}
with $\vecmu = \EspSymbol \transpose{\begin{bmatrix} \vecZ, & \transpose{\matU}\vecZ^+ \end{bmatrix}}$, $\vecmu\etoile = \EspSymbol \vecZ\etoile$ and covariance matrices, 
\begin{eqnarray*}
\matK^+ \coloneqq \begin{bmatrix}
\matK_{11} & \matK_{12} \matU   \\
\transpose{\matU} \matK_{21} & \transpose{\matU} \matK_{22} \matU
\end{bmatrix} &;&
\matH^+ \coloneqq \begin{bmatrix}
\Cov{\vecZ, \transpose{\vecZ\etoile}} \\
\transpose{\matU} \Cov{\vecZ^+, \transpose{\vecZ\etoile}}
\end{bmatrix}~,
\end{eqnarray*}
where $\matK_{11} = \Cov{\vecZ, \transpose{\vecZ}}, \matK_{12} = \Cov{\vecZ, \transpose{\vecZ^+}}, \matK_{21} = \Cov{\vecZ^+, \transpose{\vecZ}}$ and $\matK_{22} = \Cov{\vecZ^+, \transpose{\vecZ^+}}$.
\end{thmE}
\begin{proofE}
For a single prediction point, the co-Kriging predictor in \eqref{eqn:KrigingPredictorDerivativeObs} can be re-written as 
\begin{equation}
\transpose{(\vecalpha^+)} \begin{bmatrix} \vecZ \\ \transpose{\matU} \vecZ^+ \end{bmatrix}
\end{equation}
We proceed similarly to the proof of Proposition \ref{prop: OrdinaryKriging} where we wrote the Kriging error variance in \eqref{eqn: krigingopttraceform} as, 
$$\Delta(\vecs\etoile) =  \transpose{\vecalpha^+} \matK^+ \vecalpha^+ - 2 \transpose{\vecalpha^+}\vech^+(\vecs\etoile) + \sigma^2(\vecs\etoile) $$
note that the covariance matrices are modified based on the block structure of the observation vector $\transpose{[\vecZ, \transpose{\matU} \vecZ^+]}$. More explicitly, the block covariance matrices write as, 
$$
\matK^+ = \begin{bmatrix}
\Cov{\vecZ, \vecZ} & \Cov{\vecZ, \vecZ^+} \matU \\
\transpose{\matU}\Cov{\vecZ^+, \vecZ} & \transpose{\matU}\Cov{\vecZ^+, \vecZ^+} \matU
\end{bmatrix}
$$
and, 
$$
\vech^+(\vecs\etoile) = \begin{bmatrix}
\Cov{\vecZ, Z(\vecs\etoile} \\
\transpose{\matU} \Cov{\vecZ^+, Z(\vecs\etoile)}
\end{bmatrix}
$$ 
The minimization under unbiasedness constraints leads to the same constrained minimization problem as described in the proof of Proposition \ref{prop: OrdinaryKriging}, 
\begin{equation*}
\Delta_\lambda(\vecs\etoile) =  \transpose{\vecalpha^+} \matK^+ \vecalpha^+ - 2 \vech^+(\vecs\etoile) -2 \lambda \left( \transpose{\vecalpha^+} \vecmu - \vecmu\etoile\right) \, .
\end{equation*}
Since we retrieve the expressions of the same form as in the proof of Proposition \ref{prop: OrdinaryKriging}, we can follow the exact same procedure and set the gradient equal to zero and solve for the resulting system of equations to find the Lagrange multiplier, 
$$
\lambda(\vecs\etoile) = \left( \transpose{\vecmu} (\matK^+)^{-1} \vecmu \right)^{-1}
 \left( \mu\etoile - \transpose{\vecmu} (\matK^+)^{-1} \vech^+(\vecs\etoile) \right)  \, .$$
and optimal matrix of weights (for multi-point at $q$ prediction locations $\vecs\etoile_1,\ldots,\vecs\etoile_q$) as, 
\begin{equation}
\matAlpha^+ =  \left(\matK^+\right)^{-1} \left(\matH^+ + \vecmu \transpose{\veclambda} \right)~, 
\end{equation}
where $\matAlpha^+ = \cbind{\vecalpha^+(\vecs\etoile_1), ..., \vecalpha^+(\vecs\etoile_q)}$, $\matH^+ = \cbind{\vech^+(\vecs\etoile_1), ..., \vech^+(\vecs\etoile_q)}$ with the vector of Lagrange multipliers $\veclambda = \left(\lambda(\vecs\etoile_1),\ldots,\lambda(\vecs\etoile_q)\right)$. 

\end{proofE}
\noindent The final prediction is simply given as $\Zck^{\etoile} = \transpose{(\matAlpha^+)} [\vecZ,~\transpose{\matU} \vecZ^+]$. Notice that the final expression for the optimal coefficients is exactly the same as Proposition \ref{prop: OrdinaryKriging} with the matrices $\matK$ and $\matH$ replaced by $\matK^+$ and $\matH^+$. Another observation is that assuming $\vecZ$ as a centered process, the prediction corresponds to that of the conditional mean of a GP with observations from the transformed GPs after the action of linear differential operators. An important consideration is that the matrix inversion of $\matK^+$ directly, incurs a computational cost of $\mathcal{O}\left( (n + p)^3 \right)$ which can be significantly more compared to simple Kriging since it is desirable to have $p \gg n$. This is because intuitively, a sufficiently dense set of collocation ensures constraint satisfaction with higher probability.

\begin{Example}[Working example]
\label{example: CKWorkingExample}
As a gentle introductory example, consider a simplistic, one dimensional problem setup where the true function, $f(x) = sin~x$. We proceed to assess the impact of the following differential equation as auxiliary information:
\begin{equation}
f + f'' = 0
\end{equation}
To limit the dimension of the problem, suppose there exist three points: a true observation at $\vecs^1 = (x_1, 0)$; a point of prediction $\vecs\etoile = (x\etoile, 0)$; and another point $\vecs^c = (x\etoile, 2)$, part of the differential equation at the target location $x\etoile$.
\begin{equation}
\label{eqn: 1Ddiffconst}
Z(\vecs\etoile) + Z(\vecs^c) = 0
\end{equation}
Furthermore, we operate under the simple Kriging assumption, which is to say that the extended process is mean zero, $E Z = 0$ everywhere in the extended design space.
The Kriging predictor is built as a linear combination of the true observation and the differential equation at $x\etoile$,
\begin{eqnarray}
\label{eqn: ObservationbasedWorkingExample}
\Zck (\vecs\etoile) &=& \alpha_1(\vecs\etoile)~ Z(\vecs^1) + \alpha_2(\vecs\etoile) ~\left( Z(\vecs\etoile) + Z(\vecs^c) \right)
\end{eqnarray}
$\alpha_1(\vecs\etoile), \alpha_2(\vecs\etoile)$ can be resolved as the minimizers of predictive mean square error (MSE) at $\vecs\etoile$,
\begin{eqnarray}
	\underset{{\alpha_1, ~\alpha_2~\in \mathbb{R}^2}}{\argmin}  
		\EspSymbol\left[ \left( \Zck(\vecs\etoile) - Z(\vecs\etoile)\right)^2 \right]~,
		\label{eqn: SimpleKrig}
\end{eqnarray}
Setting $\vecmu, ~ \vecmu\etoile = 0$ in \eqref{eqn: OrdinaryCollocatedcoKriging}, the final prediction can be expressed as, 
\begin{eqnarray}
	\label{eqn: ToyExamplecollocatedcoKriging}
	\widehat{Z}(\vecs\etoile) &=& \transpose{(\matH^+)} (\matK^+)^{-1} \transpose{\left[ Z(\vecs^1), 0 \right]}~,
\end{eqnarray}
  
where,
\begin{eqnarray*}
\matK^+ &\coloneqq& 
		\begin{bmatrix}
			\VarSymbol Z(\vecs^1) & \CovSymbol\left[ Z(\vecs^1),  Z(\vecs\etoile) + Z(\vecs^c)\right] \\
			. & \VarSymbol \left[Z(\vecs\etoile) + Z(\vecs^c)\right]
		\end{bmatrix} \\
		\matH^+ &\coloneqq& 
		 	\begin{bmatrix}
				\CovSymbol\left[ Z(\vecs^1), Z(\vecs\etoile)\right]  \\ \CovSymbol\left[ Z(\vecs\etoile) + Z(\vecs^c),  Z(\vecs\etoile)\right]
			\end{bmatrix}
\end{eqnarray*}

For the covariance computations, we have used the isotropic squared exponential kernel as described in eq-\eqref{eqn: sqexp}. The relevant derivative expressions are given in Section \ref{sec: kernelderiavtives}. Figure \ref{fig: ToyExample1}, illustrates this example. Given that we consider only one derivative observation, there is only a minor difference between the collocated co-Kriging and the usual Kriging prediction.

\begin{figure}[hbtp]
\centering
\includegraphics[width = 230px]{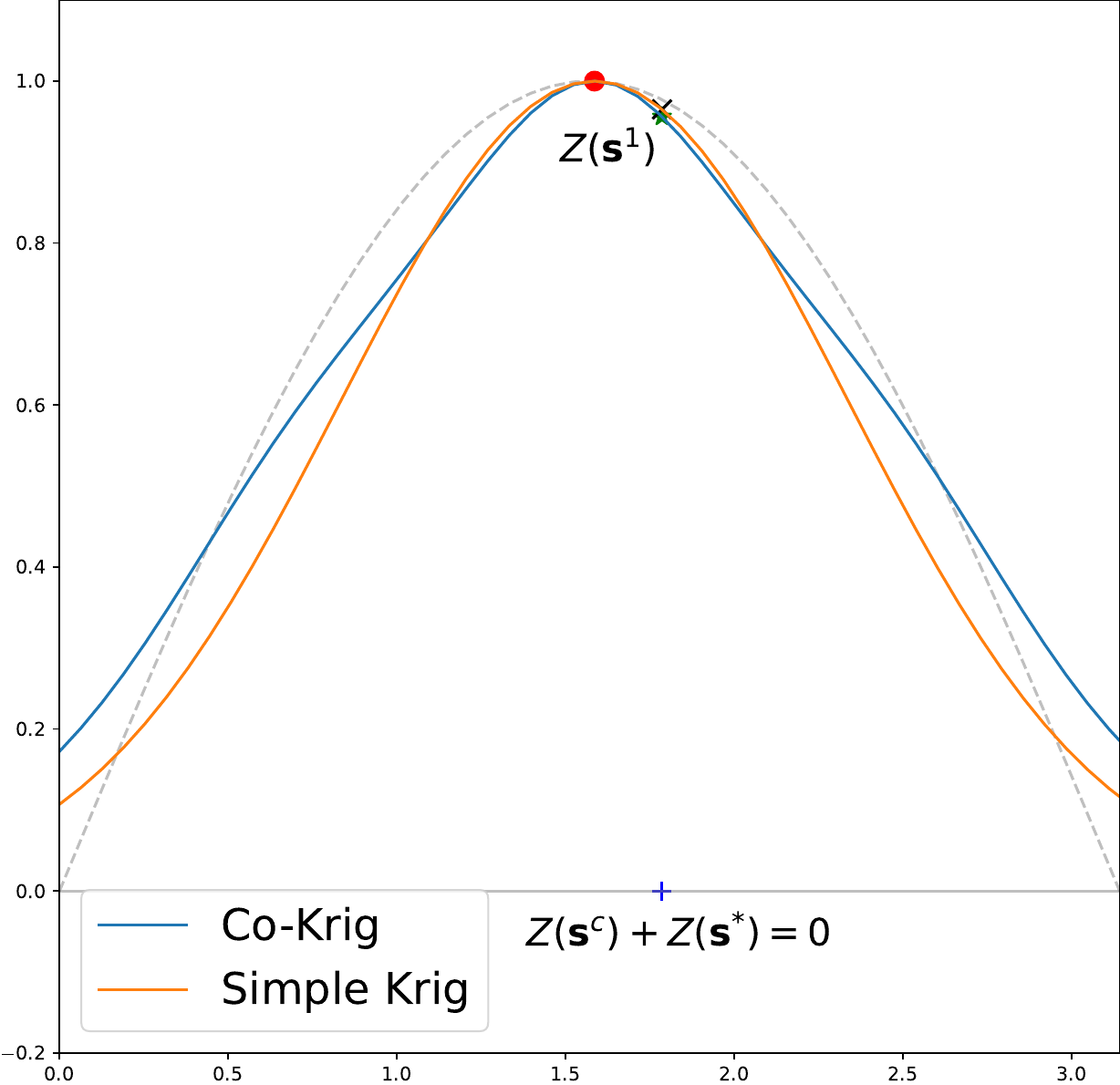}
\caption{Demonstrating the difference between simple Kriging (in orange) and the prediction (in green and the blue curve) with a single collocation point (in blue).}
\label{fig: ToyExample1}
\end{figure}
\end{Example}
\begin{Remark}
It is not necessary to separate the pure and differential equations related observations. We can incorporate the entire vector of observations using a unified $\matU^+$ matrix,
\begin{eqnarray}
(\matU^+)^T \left(\vecZ\right)_{i=1}^{n+c} &=& \left(\vecv\right)_{i=1}^{n+p} \\
\matU^+ &=& \begin{bmatrix}
	\bm{I}_{n\times n} & \bm{0} \\
	\bm{0} & \matU
	\end{bmatrix}_{(n+c) \times (n+p)}	  
\end{eqnarray}
The separation is important for compatibility of notation for the two methods we explain here. If we use the extended $\matU^+$, it is worth noting that the final expressions will be relatively simpler since we won't need the covariances between pure observations $\matK$ and the cross-covariances between pure $\vecZ$ and PDE-related observations $\vecZ^+$. 
\end{Remark}
\subsection{Linear Lagrangian Kriging}
\label{sec: ConstrainedKriging}
It is common to encounter constraints in the applications of Kriging. For example, the ordinary Kriging estimator \eqref{eqn: ordinarykriging} involves the unbiasedness constraint as part of its formulation \cite{Wackernagel2003}. Similarly, unbiasedness constraints for non-linear functionals of the Kriging estimator, was first addressed in \cite{Cressie1993}. The constrained Kriging predictor, preserves both the expectation and the variance which is sufficient to retain the same distribution as the original Gaussian field, under a non-linear transformation. This idea was furthered in \cite{Aldworth2003} where they developed an estimator that preserves the all covariances. This class of estimators are aptly named as covariance matching constrained Kriging estimators. These estimators have been coded as part of an R package \cite{Hofer2011} and have found applications for the problem of coal quality estimation \cite{Ertunc2013} among others. Lagrange multipliers have been used in constrained stochastic optimization \cite{Bhatnagar2011} for constraints on long-run averages of the cost function and inequality relations. They use a derivative-based optimization algorithms on the Lagrangian. 
\newline \newline
\noindent Despite a vast volume of articles on this subject, the nature of constraints is almost always limited to simple equality and inequality constraints. To the best of our knowledge, there has been no investigation on operator equations or differential operator equations as constraints. We adopt a simple Lagrange multiplier based resolution of differential constraints, similar to \cite{Cressie1993, Aldworth2003, Bhatnagar2011}. We consider Kriging predictors such that they satisfy pointwise linear PDE constraints specifically at the points of prediction. We exploit the matrix formulation as described in \eqref{eqn: constraints2} to obtain equality constraints. 
\newline \newline
A crucial difference from Section \ref{sec: derivativeobs} is that the distinction between prediction points and collocation points is not clear anymore. The linear Lagrangian predictor at an unknown point $\vecs\etoile$ is $\Zlk(\vecs\etoile) = \transpose{\alpha(\vecs\etoile)} \vecZ$. It would be preferable to impose constraints at all prediction points. Hence we have $q$ constraint equations and $\lagrange$ locations $\vecZ^{\star} = \left( Z(\vecs_{1}\etoile), \ldots, Z(\vecs_{\lagrange}\etoile)\right)^T$ involved in the constraint equations with $\lagrange \geq q$ since each constraint equation is composed of atleast one location.     
\begin{eqnarray}
\transpose{\matU} \vecZ\etoile &=& \vecv\etoile~, 
\end{eqnarray}  
where $\matU\in \R^{\lagrange \times q} $ $\vecv\etoile \in \R^q$. We will predict all $\lagrange$ points as part of our model output. The mean squared error (MSE) optimization problem in this case can be stated as, 
\begin{align}
\label{eqn:LagrangianKrigOpt}
\underset{\matAlpha \in \R^{n \times \lagrange}}{\argmin} \quad &\Trace{\transpose{\matAlpha}\matK \matAlpha} - 2 ~\Trace{\transpose{\matAlpha} \matH} + \Trace{\matK\etoile} \\
\text{subject to: \quad} & \transpose{\matAlpha} \vecmu = \vecmu\etoile  
 \tag*{(C1)}~,  \\
  &\transpose{\matU} \transpose{\matAlpha} \vecZ = \vecv\etoile \tag*{(C2)}~,
\end{align}
where C1 is the unbiasedness constraint and C2 represents the differential constraints and  $\matH \in \R^{n \times \lagrange}$. 

\newcommand{\gammamu}{\gamma_1}
\newcommand{\gammaZ}{\gamma_2}
\newcommand{\gammamuZ}{\gamma_{3}}

\begin{thmE}[Ordinary Kriging with differential constraints][\myoption]
\label{prop: CosntrainedOK}
The matrix of the linear predictor's optimal weights minimizing the optimization program~\eqref{eqn:LagrangianKrigOpt} under both unbiasedness and differential constraints~(C1, C2) is given by
	\begin{equation}
\matAlpha = \matK^{-1} \left(\matH +  \vecmu \transpose{\veclambda} +  \vecZ \transpose{\veclambda'}\transpose{\matU}\right)
	\end{equation}
with Lagrange multipliers, in the case where $\transpose{\matU}\matU$ is invertible, $\gamma_1 \neq 0$ and $\gammaZ- \frac{\gammamuZ^2}{\gammamu} \neq 0$,
\begin{equation*}
\accolade{
		{\veclambda} & =& {\gammamu}^{-1} \left({\vecmu\etoile} - \transpose{\matH} \matK^{-1} \vecmu - \gammamuZ {\matU}\veclambda'\right)\\
{\veclambda'}  &=& 	 \left(\gammaZ- \frac{\gammamuZ^2}{\gammamu} \right)^{-1} \left( \transpose{\matU}\matU\right)^{-1} \left({\vecv} - 
\transpose{\matU}\transpose{\matH} \matK^{-1} \vecZ
 - \frac{\gammamuZ}{\gammamu} \transpose{\matU}\left( {\vecmu\etoile} - \transpose{\matH} \matK^{-1} \vecmu \right)\right)
}
\end{equation*}
With scalar values $\gammamu= \transpose{\vecmu} \matK^{-1} \vecmu$,
$\gammaZ= \transpose{\vecZ} \matK^{-1} \vecZ$ and $\gammamuZ= \transpose{\vecmu} \matK^{-1} \vecZ$. 
\end{thmE}
\begin{proofE}
The resolution is quite similar to the one in~\cite{Rulliere2025} where constraint of order $\vecm=(0, \ldots, 0)$ are considered (prescribed mean of predicted values, e.g. for adverse modelling).
Here, using Lagrange multipliers $2\veclambda$ and $2\veclambda'$, with respective sizes $q\times 1$ and $p \times 1$, we have to solve 
\begin{equation*}
\min_{\matAlpha} \left[ \Delta\etoile(\matAlpha) - 2 \left( \transpose{\vecmu} \matAlpha - \transpose{\vecmu\etoile}\right)\veclambda - 2 \left( 	\transpose{\vecZ } \matAlpha \matU - \transpose{\vecv}\right) \veclambda' \right]
\end{equation*}
Using classical matrix differentiation results, using the symmetry of $\matK$, we can write:
\begin{equation*}
\frac{ \partial
}{\partial A} \Delta\etoile(\matAlpha) = 	2\matK\matAlpha -2 \matH
\end{equation*}
And finally a set of necessary conditions for optimality, to be solved in $\matAlpha$, $\transpose{\veclambda}$ and $\transpose{\veclambda'}$
\begin{equation*}
\accolade{
 	\matK\matAlpha - \matH -  \vecmu \transpose{\veclambda} -  \vecZ \transpose{\veclambda'} \transpose{\matU} &=& \veczero\\
	\transpose{\vecmu} \matAlpha - \transpose{\vecmu\etoile} &=& \veczero\\
\transpose{\vecZ} {\matAlpha} \matU - \transpose{\vecv} &=& \veczero 
 }
\end{equation*}
From the first Equation, we get
$$ \matAlpha = \matK^{-1} \left(\matH +  \vecmu \transpose{\veclambda} +  \vecZ \transpose{\veclambda'}\transpose{\matU}\right)$$
Then injecting this value into both constraints
\begin{equation*}
\accolade{
	\transpose{\vecmu} \matK^{-1} \left(\matH +  \vecmu \transpose{\veclambda} +  \vecZ \transpose{\veclambda'}\transpose{\matU}\right)  &=& \transpose{\vecmu\etoile}\\
	\transpose{\vecZ} \matK^{-1} \left(\matH +  \vecmu \transpose{\veclambda} +  \vecZ \transpose{\veclambda'}\transpose{\matU}\right) \matU  &=& \transpose{\vecv} 
}
\end{equation*}

Denoting the scalar values $\gammamu= \transpose{\vecmu} \matK^{-1} \vecmu$ and 
$\gammaZ= \transpose{\vecZ} \matK^{-1} \vecZ$ and $\gammamuZ= \transpose{\vecmu} \matK^{-1} \vecZ$, we get
\begin{equation*}
\accolade{
	\transpose{\vecmu} \matK^{-1} \matH + \gammamu \transpose{\veclambda} + \gammamuZ \transpose{\veclambda'}\transpose{\matU} &=& \transpose{\vecmu\etoile}\\
	\transpose{\vecZ} \matK^{-1} \matH \matU +  \gammamuZ \transpose{\veclambda} \matU +  \gammaZ \transpose{\veclambda'} \transpose{\matU}\matU  &=& \transpose{\vecv} 
}
\end{equation*}
Hence
\begin{equation*}
\accolade{
&&	\gammamu \transpose{\veclambda}  = \transpose{\vecmu\etoile} - \transpose{\vecmu} \matK^{-1} \matH - \gammamuZ \transpose{\veclambda'}\transpose{\matU}\\
&&	\transpose{\vecZ} \matK^{-1} \matH \matU+  \frac{\gammamuZ}{\gammamu} \left( \transpose{\vecmu\etoile} - \transpose{\vecmu} \matK^{-1} \matH - \gammamuZ \transpose{\veclambda'}\transpose{\matU}\right) \matU+  \gammaZ \transpose{\veclambda'} \transpose{\matU}\matU  = \transpose{\vecv} 
}
\end{equation*}
and
\begin{equation*}
\accolade{
		\transpose{\veclambda} & =& \frac{1}{\gammamu} \left(\transpose{\vecmu\etoile} - \transpose{\vecmu} \matK^{-1} \matH - \gammamuZ \transpose{\veclambda'}\transpose{\matU} \right)\\
		 \left(\gammaZ- \frac{\gammamuZ^2}{\gammamu} \right)\transpose{\veclambda'} \transpose{\matU}\matU &=& \transpose{\vecv} - \transpose{\vecZ} \matK^{-1} \matH \matU - \frac{\gammamuZ}{\gammamu} \left( \transpose{\vecmu\etoile} - \transpose{\vecmu} \matK^{-1} \matH \right)\matU
}
\end{equation*}
and in the case where $\transpose{\matU}\matU$ is invertible
\begin{equation*}
\accolade{
	\transpose{\veclambda} & =& \frac{1}{\gammamu} \left(\transpose{\vecmu\etoile} - \transpose{\vecmu} \matK^{-1} \matH - \gammamuZ \transpose{\veclambda'}\transpose{\matU} \right)\\
	 \left(\gammaZ- \frac{\gammamuZ^2}{\gammamu} \right)\transpose{\veclambda'}  &=& \left(\transpose{\vecv} - \transpose{\vecZ} \matK^{-1} \matH \matU - \frac{\gammamuZ}{\gammamu} \left( \transpose{\vecmu\etoile} - \transpose{\vecmu} \matK^{-1} \matH \right)\matU\right)\left( \transpose{\matU}\matU\right)^{-1}
	}
\end{equation*}
Hence the result.\\
Note that it corresponds to the result in \cite{Rulliere2025}, Proposition 3 ``\textit{Joint Kriging weights under a predicted values constraint}'' with notational changes: (this paper) $\leftrightarrow$ (cited paper), $\vecZ \leftrightarrow \transpose{\mathbb{Y}}$, $\matU \leftrightarrow \boldsymbol{\pi}$, $\vecv \leftrightarrow \transpose{\boldsymbol{m}}$, $\veclambda' \leftrightarrow \transpose{\veclambda'}$, $\vecmu \leftrightarrow \vecun_n$, $\vecmu\etoile \leftrightarrow \vecun_q$, $\gamma_1 \leftrightarrow \delta$, $\gamma_3 \leftrightarrow \vecu$.
\end{proofE}
Like usual Kriging, the final predictions can be written as $\Zlk^{\etoile} = \transpose{\matAlpha} \vecZ$. In terms of complexity, there are two matrix inversions, $\matK^{-1}$ and the sparse matrix $\left(\transpose{\matU}\matU \right)^{-1}$ with computational costs of $\mathcal{O}(n^3)$ and $\mathcal{O}(q^3)$ respectively. Typically, $q \gg n$. This is already an improvement over Section \ref{sec: derivativeobs} (complexity: $\mathcal{O}((n + p)^3)$). With Lagrangian Kriging, owing to the analytical resolution via Lagrange multipliers, the predictions will always be consistent with the differential constraints. However, this does not imply anything about the prediction quality. For example, going back to the problem in Example $\ref{example: CKWorkingExample}$, we might have unrealistic predictions for both $f$ and $f''$ while satisfying $\hat{f} + \hat{f}'' = 0$  at the points of prediction. For co-Kriging, the constraint satisfaction is related to the proximity from collocation points.
\begin{Example}[Working example (continued)] 
\label{example: LKWorkingExample}
We borrow the same example as described in example \ref{example: CKWorkingExample}, and resolve the constrained optimization problem instead. When treating the differential information as a constraint, the predictor is composed of the true observation only. Since we have no information about $\vecs^c$, it is to be predicted as well, 

\begin{eqnarray}
	\label{eqn: ConstrainedKrigPredictor}
	\Zlk(\vecs\etoile) = \alpha(\vecs\etoile) Z(\vecs^1)~;~\Zlk(\vecs^c) = \alpha(\vecs^c) Z(\vecs^1) 
\end{eqnarray} 

Additionally, we have a constraint to satisfy at the prediction location which amounts to a constrained optimization problem, 

\begin{eqnarray}
	\label{eqn: ConstrainedKrigOpt}
	&\underset{\alpha_1,~\alpha_2~\in \mathbb{R}^2}{\argmin}& \EspSymbol \left[ \left( \Zlk(\vecs\etoile) - Z(\vecs\etoile)\right)^2\right]  +  \EspSymbol \left[ \left( \Zlk(\vecs^c) - Z(\vecs^c)\right)^2\right]  \\ 
	&\text{subject to: \quad}& 0 = \Zlk(\vecs\etoile) + \Zlk(\vecs^c)  \nonumber \\
	&& = \left( \alpha_1  + \alpha_2 \right) Z(\vecs^1) \nonumber
\end{eqnarray}

We can solve this using a Lagrangian multiplier $\lambda$ and setting the gradient with respect to $\alpha_1, \alpha_2$ and $\lambda$ to zero (mathematical working in the appendix \ref{sec: WorkingExampleLKmath}) to obtain, 
\begin{eqnarray}
	\label{eqn: ConstrainedKrigPrediction}
	\Zlk(\vecs\etoile) &=& \left(\frac{\CovSymbol \left[\vecs^1,\vecs\etoile\right] - \CovSymbol \left[\vecs^1,\vecs^c\right]}{2 ~ \VarSymbol \vecs^1 } \right)  Z(\vecs^1)
\end{eqnarray} 

Once again, we present an accompanying visual, Figure \ref{fig: ToyExampleCK}, illustrating the difference between the Lagrangian Kriging predictor and usual Kriging predictor for the working example. There is a noticeable difference between the two predicted curves.

\begin{figure}[ht!]
\centering
\includegraphics[width=230px]{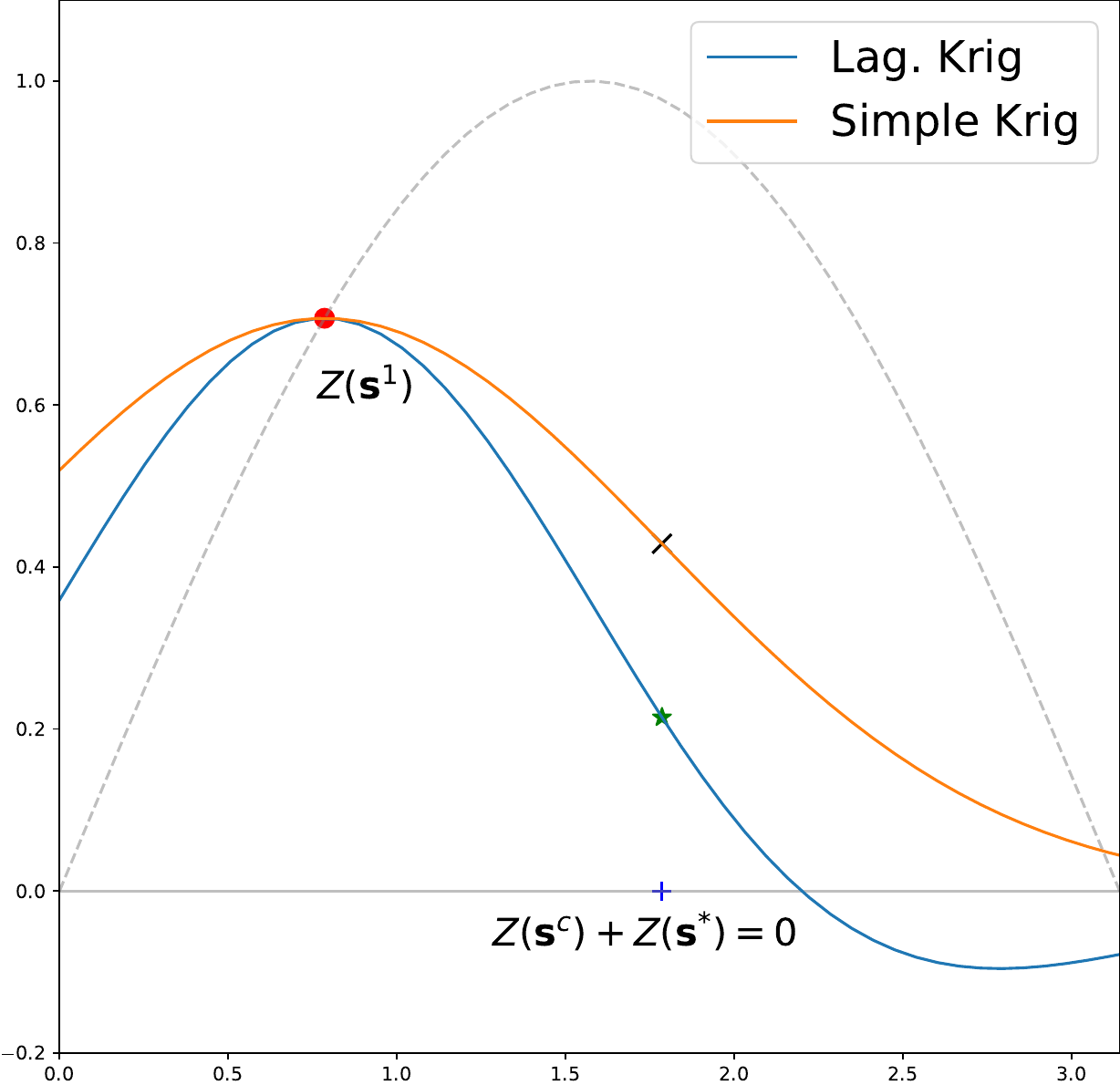} 
\caption{Comparing the simple Kriging (in orange) and Lagrangian Kriging prediction (in green and the blue curve) with a single point constraint.}
\label{fig: ToyExampleCK}
\end{figure}
\end{Example}

\begin{Remark}[Constraints without the primary function]
\label{remark: LK}
The Lagrangian Kriging predictor can only account for differential information which consists of the function of interest undifferentiated, in addition to its derivatives. If the differential information does not involve the primary function, the primary associated Kriging objective is essentially unconstrained and the predictions of the primary function would correspond to the usual Kriging predictions. This stems from the fact that the linear coefficients involved in the prediction of the primary function are independent of the linear coefficients involved in the prediction of its higher order derivatives. This being a recurrent concept, we detail this issue when predicting $f: \R^n \to \R$ (primary variable) under constraints involving only the Laplacian of $f$, i.e., $\nabla^2 f$. From \eqref{eqn:LagrangianKrigOpt}, recall that the Kriging objective is simply $\Trace{\transpose{\matAlpha}\matK \matAlpha} - 2 ~\Trace{\transpose{\matAlpha} \matH}$ with $\matAlpha \in \R^{n \times \lagrange}$. Suppose we want to predict $f$ at $q$ locations, in that case we can partition the matrix $\matAlpha$ of coefficients into two blocks,
\begin{equation*}
\matAlpha = \left[\matAlpha_f, \matAlpha_{\nabla^2f}\right]~,
\end{equation*}
with $\matAlpha_f \in \R^{n \times q}$ involved in the prediction of the primary variable and $\matAlpha_{\nabla^2f} \in \R^{n \times \lagrange - q}$ involved in the prediction of the Laplacian components $\{\partial^2_{x_i} f: i = 1,\ldots,n\}$. Therefore, by separating the objective, the optimization problem can be written as two sub-problems, an unconstrained optimization problem with respect to $\matAlpha_f$
\begin{eqnarray*}
\underset{\matAlpha_f \in \R^{n \times q}}{\argmin} \Trace{\transpose{\matAlpha_f}\matK \matAlpha_f} - 2 ~\Trace{\transpose{\matAlpha_f} \matH} 
\end{eqnarray*}
and a constrained optimization problem, of exactly the same form as presented in \eqref{eqn:LagrangianKrigOpt}, with respect to $\matAlpha_{\nabla^2f}$, 
\begin{align*} 
 \underset{\matAlpha_{\nabla^2f} \in \R^{n \times \lagrange - q}}{\argmin} & \Trace{\transpose{\matAlpha}_{\nabla^2f}~\matK~\matAlpha_{\nabla^2f}} - 2 ~\Trace{\transpose{\matAlpha}_{\nabla^2f}~\matH} \\
 \text{subject to: \quad} & \transpose{\matU}\transpose{\matAlpha}_{\nabla^2f}~\vecZ = \vecv
\end{align*}
and consequently, the final predictions of $f$ relying on $\matAlpha_f$ would correspond to that of classical Kriging.              
\end{Remark}
\subsection{Comparing the two methods}
\label{sec: comparison}

We delve deeper into the question of how co-Kriging and Lagrangian Kriging are different. Since, they are based on the classical Kriging approach, we expected similar results in the beginning. However, there are some key fundamental differences which we will address in detail in this Section. This is particularly relevant for a practitioner willing to make an informed choice among the two models. 

\paragraph{Analytical form.}
To explore the key mathematical difference in the two methods, we commence by reducing the expression for the collocated simple co-Kriging approach by exploiting the block structure of matrices $\matK^+$ and $\matH^+$. We assume that $\vecZ$ and $\vecZ\etoile$ are centered\footnote{The simple Kriging framework allows for simplified expressions since we can ignore the unbiasedness constraint}, i.e., $\vecmu, ~\vecmu\etoile = 0$. Given an observation vector of $n$ points $\vecZ$ and, for a fair comparison, differential information at all $q$ prediction points encoded by the matrix equation $\transpose{\matU} \vecZ\etoile = \vecv\etoile$ where $\vecv\etoile$ is of size $q$ and $\vecZ\etoile$ of size $\lagrange$, the $2 \times 2$ block inverse of $\matK^+$, assuming $\Cov{\vecZ, \transpose{\vecZ}}$ and the Schur's complement are invertible, 
\begin{eqnarray}
\begin{bmatrix}
A & B \\
C & D
\end{bmatrix}^{-1} &=& 
\begin{bmatrix}
A^{-1} + A^{-1} B (D - CA^{-1}B)^{-1} C A^{-1} & -A^{-1} B (D - CA^{-1}B)^{-1} \\ 
-(D - CA^{-1}B)^{-1} C A^{-1} & (D - CA^{-1}B)^{-1}
\end{bmatrix} 
\end{eqnarray} 
where, $A = \matK_{11},~B = \matK_{12} \matU,~C = \transpose{\matU} \matK_{21}$ and $D = \transpose{\matU} \matK_{22} \matU$ with the same definitions of $\matK_{ij}$ as in proposition \ref{prop: collocatedcokriging} except that $\vecZ^+$ replaced by $\vecZ\etoile$ since we consider all our prediction points as collocation points. %
It is also useful to recognize that $\matH^+ = \transpose{\left[  \matK_{12} , \transpose{\matU} \matK_{22} \right]}$ and on %
replacing $\matK_{11}$ with $\matK$ and $\matK_{12},~ \matK_{21}$ with $\matH,~ \transpose{\matH}$, where $\matK,~\matH$ are the usual Kriging covariance matrices (see proposition \ref{prop: OrdinaryKriging}) to compute the final predictor as $\transpose{(\matAlpha^+)} \transpose{\left[ \vecZ, \vecv\etoile\right]}$,
\begin{eqnarray}
\vecZck\etoile 
&=& \matK_{2\vert 1} \matU \left( \transpose{\matU} \matK_{2\vert 1} \matU \right)^{-1}  \left(\vecv\etoile - \transpose{\matU} \transpose{\matH} \matK^{-1} \vecZ\right) ~+~ \transpose{\matH} \matK^{-1} \vecZ
\end{eqnarray} 
where $\matK_{2\vert 1} = \matK_{22} - \transpose{\matH} \matK^{-1} \matH$ which is also the Schur's complement $\matK^+ / \matK$ and most interestingly $\matK_{2\vert 1} = \Cov{\vecZ\etoile \vert \vecZ}$. When comparing with the predictor in Proposition \ref{prop: CosntrainedOK}, under the simple Kriging assumption a lot of terms disappear, notably all terms involving $\veclambda, \gamma_1, \gamma_3$ in addition to $\vecmu,~\vecmu\etoile$ as mentioned before, rendering a simplified final predictor,

\begin{eqnarray}
\label{eqn: LKfinal}
\vecZlk\etoile &=& \matU \left(\transpose{\matU}\matU\right)^{-1} \left(\vecv\etoile - \transpose{\matU}\transpose{\matH} \matK^{-1} \vecZ \right) ~+~ \transpose{\matH} \matK^{-1} \vecZ 
\end{eqnarray}
 
where the scalar $\transpose{\vecZ} \matK^{-1} \vecZ$ has been assumed to be non-zero in the calculations. A quick inspection reveals that both predictors are affine in $\vecZ$ and it is easy to obtain $\vecZlk\etoile$ from $\vecZck\etoile$ by replacing $\matK_{2 \vert 1}$ by the identity matrix. The $\matK_{2 \vert 1}$ is responsible for augmenting the co-Kriging model complexity. It is also immediately obvious that Lagrangian Kriging does not account for cross-covariances among the constraint equations whereas co-Kriging does. 
\begin{Remark}[Constraints without the primary function (contd.)]
\label{remark: LK continued}  
The fact that the Lagrangian Kriging problem described in Remark \ref{remark: LK}, is consistent with \eqref{eqn: LKfinal} is briefly explained here. To reiterate what was explained in Remark \ref{remark: LK}, consider the task of predicting $f: \R^2 \to \R$ while imposing constraints on $\nabla^2 f$. The matrix $\matH$ would consist of three kinds of predictions $\hat{f}, \hat{f}_{xx}$ and $\hat{f}_{yy}$. Out these three, there are constraints imposed on $\hat{f}_{xx}$ and $\hat{f}_{yy}$ which will reflect in the $\matU$ matrix. However, since Lagrangian Kriging does not account for any cross-covariance between $f$ and its second order derivatives, we can safely separate out the task of predicting $\hat{f}$ as unconstrained Kriging given $n$ observations of $f$. This is a major fallacy since the predictions independent (even though $f, f_{xx}$ and $f_{yy}$ mathematically related) of each other unless they are explicitly involved in a constraint equation. The problem becomes evident in cases where the constraints does not explicitly involve the function of interest.  
\end{Remark}

\paragraph{Computational complexity.}
To discuss the computational aspects, we recall that $n$ is the number of observations in the extended design space, $p$ is the number of collocation points, $q$ is the number of prediction points.%

\begin{table}[ht]
\caption{Complexity (in one dimension) and the number of total covariance derivative computations.} \vspace{1em}
\label{tab: Complexity}
\centering
\begin{tabular}{c c c}
\toprule \small
&  \small{Complexity} & \makecell{\small Covariance \\ computations} \\
\midrule%
\multirow{2}{*}{\small Co-Krig.} & \multirow{2}{*}{\small{$\mathcal{O}((n+p)^2(n+p+q))$}} & \multirow{2}{*}{\small $\mathcal{O}((n+p+q)(n+p))$} \\ 
&&\\
\cline{1-3}
\multirow{2}{*}{\small Lag. Krig.} & \multirow{2}{*}{\small $\mathcal{O}((n + q)^2)$} & \multirow{2}{*}{\small$\mathcal{O}(n(n + q))$}
\end{tabular}
\end{table}

\noindent Table \ref{tab: Complexity}, summarizes the computational aspects of Co-Kriging and Lagrangian Kriging. Considering all the matrix operations, co-Kriging complexity is $\mathcal{O}((n+p)^2(n+p+q))$ (with or without block inverse). On the other hand, Lagrangian Kriging complexity is $\mathcal{O}((n+q)^2)$.
An interesting scenario is when $p, q \gg n$ wherein, co-Kriging would be $\mathcal{O}(p^3)$ versus Lagrangian Kriging would only be $\mathcal{O}(q^2)$ as it does not rely on collocation points. This is to emphasize the significant computational edge we obtain with Lagrangian Kriging. Apart from the size of the involved matrix inversions and multiplications, there are also fewer covariance computations involved in Lagrangian Kriging since the entire $\matK_{22}$ matrix only appears in $\vecZck\etoile$ expression and not in $\vecZlk\etoile$. This represents a significant portion of the cost since higher order kernel derivatives are computationally expensive. Lagrangian Kriging enjoys a relatively faster run time, noticeable in practical implementations as well (Table \ref{tab: CPUtimes}). Both these methods are prohibitively expensive for many observations and many prediction points. Scalability to higher dimensions $d$ comes at a cubic cost, a $d^3$ factor, times the one-dimensional costs since each dimension of a regular observation is counted as a separate observation resulting in $n\times d$ observations.

\subsection{Calibrating covariance parameters}
The following Section discusses the crucial hyperparameter tuning aspect for kernels used in Kriging. Typically, the Gaussian process prior assumption simplifies this task as it provides an analytical expression for the likelihood of observed data. Seeking covariance parameters that maximize the multivariate normal likelihood is called the maximum likelihood estimate (MLE) and has been used in physics-informed GPs \cite{Solak2002, Jidling2017, Raissi2017}. The unconditional MLE was used in \cite{DaVeiga2012} for a GP under linear constraints. A constrained MLE criteria was proposed by \cite{Bachoc2019}. However, we have avoided the GP prior assumption in favour of generality. Without the analytical likelihood expression, a more suitable choice is the Leave-one-out cross validation (LOOCV) based criteria \cite{Bachoc2013} where they demonstrate its efficacy under model misspecification which happens to be an added benefit. LOOCV is well established in the broader GP and Kriging literature but to the best of our knowledge, has not found use in physics-informed or derivative based GP approaches. In the case of simple Kriging, the usual cross validation loss can be found using a single evaluation of the inverse covariance matrix $\matK^{-1}$ and observations $\vecZ$, thanks to the virtual cross validation loss formulas first introduced by \cite{Dubrule1983}. For our experiments, we have chosen the infinitely differentiable, isotropic, squared exponential kernel (or RBF kernel),
\begin{eqnarray}
\label{eqn: sqexp}
k(\vecx,\vecx'\vert \sigma^2, \theta) &=& \sigma^2 \exp{\left(\frac{-\Vert \vecx - \vecx'\Vert_2^2}{2 \theta^2}\right)}   ~;~ \vecx, \vecx' \in \mathcal{X}~;~ \sigma^2, \theta \in \mathbb{R}^+ 
\end{eqnarray}      
where $\sigma^2$ is the noise or variance parameter and $\theta$ is the correlation or lengthscale parameter. A noteworthy property is that the predictions are independent of the noise parameter and only depend on $\theta$. The role of $\sigma^2$ is importance for the uncertainty quantification. Hence, the prescribed workflow is to find the optimal $\theta$ minimizing the LOOCV prediction error and then computing $\sigma^2$ that satisfies a predictive variance criteria. There is a detailed description and comparison of different predictive variance criteria in \cite{Marrel2024}. For our experiments, we compute the $\hat{\theta}$ that minimizes the predictive LOOCV MSE and then $\hat{\sigma}^2$ is chosen such that the predictive LOOCV variance criterion is equal to one, as suggested by \cite{Cressi1993}.
\begin{equation}
\label{eqn: LOOCVMSE}
\mathrm{LOOMSE}~(\theta) = \frac{1}{n}\sum_{i=1}^n (Z_i - \hat{Z}_{-i})^2
\end{equation}
\begin{equation}
\label{eqn: LOOCVVAR}
\mathrm{LOOVAR}~(\theta, \sigma^2) = \frac{1}{n} \sum_{i=1}^n \frac{(Z_i - \hat{Z}_{-i})^2}{\sigma^2 ~\Var{\hat{Z}_{-i}}\big\vert_{k_{\sigma^2 = 1, \theta}(.,.)}}
\end{equation}
where $\hat{Z}_{-i}$ is the prediction at target location $\vecs_i$ having observed $\vecZ  \backslash Z_i$ and $\Var{\hat{Z}_{-i}} = \sigma^2 ~\Var{\hat{Z}_{-i}}\big\vert_{k_{\sigma^2 = 1, \theta}(.,.)}$ is the  full conditional variance of the model factored into two independent components. The first of which is the noise variance $\sigma^2$ and the second component is the conditional variance for a chosen lengthscale $\theta$ assuming unit noise variance ($\sigma^2=1$) since it has already been accounted for in the first component, $\Var{\hat{Z}_{-i}}\big\vert_{k_{\sigma^2 = 1, \theta}(.,.)}$. Invoking the virtual cross validation loss formulas for simple Kriging from \cite[prop. 3.1]{Bachoc2013} we can write the LOOCV optimal parameters as,  
\begin{eqnarray}
\label{eqn: virtualLOOCVMSE}
\hat{\theta} &=& \underset{\theta \in \mathbb{R}^+}{\argmin} ~\frac{1}{n} \transpose{\left(\matK^{-1} \vecZ\right)} \mathrm{diag}(\matK^{-1})^{-2} \left(\matK^{-1}\vecZ \right)\big\vert_{k_{\sigma^2 = 1, \theta}(.,.)} \\
\label{eqn: virtualLOOCVVAR}
\hat{\sigma}^2 &=& \frac{1}{n} \transpose{\left(\matK^{-1} \vecZ\right)} \mathrm{diag}(\matK^{-1})^{-1} \left(\matK^{-1}\vecZ \right)\big\vert_{k_{\sigma^2 = 1, \hat{\theta}}(.,.)}
\end{eqnarray}
where $\mathrm{diag}(.)$ represents the diagonal matrix for a given matrix. Adapting to our model specifics is relatively straightforward. For Lagrangian Kriging, we can use these formulas directly if there are no constraints to consider when predicting at observation locations. The model prediction coincides with simple Kriging exactly. If this is not the case, we cannot use the virtual LOOCV formulas and need to compute the LOOCV loss as in \eqref{eqn: LOOCVMSE}-\eqref{eqn: LOOCVVAR}. In the case of simple co-Kriging, the primary observations are part of the total observations $[\vecZ, \vecZ^+]$. Hence, we use a filter matrix to retain the first $n$ elements only, the modified expressions being, 
\begin{eqnarray}
\label{eqn: CKvirtualLOOCVMSE}
\hat{\theta} &=& \underset{\theta \in \mathbb{R}^+}{\argmin} ~\frac{1}{n} \transpose{\left((\matK^+)^{-1} \begin{bmatrix}\vecZ \\ \vecZ^+\end{bmatrix}\right)} \mathrm{diag}((\matK^+)^{-1})^{-2} \begin{bmatrix}
\bm{I}_n & 0 \\
0 & \bm{0}
\end{bmatrix} \left((\matK^+)^{-1} \begin{bmatrix}\vecZ \\ \vecZ^+\end{bmatrix} \right)\Big\vert_{k_{\sigma^2 = 1, \theta}(.,.)} \\
\label{eqn: CKvirtualLOOCVVAR}
\hat{\sigma}^2 &=& \frac{1}{n} \transpose{\left((\matK^+)^{-1} \begin{bmatrix}\vecZ \\ \vecZ^+\end{bmatrix}\right)} \mathrm{diag}((\matK^+)^{-1})^{-1} \begin{bmatrix}
\bm{I}_n & 0 \\
0 & \bm{0}
\end{bmatrix} \left((\matK^+)^{-1} \begin{bmatrix}\vecZ \\ \vecZ^+\end{bmatrix} \right)\Big\vert_{k_{\sigma^2 = 1, \hat{\theta}}(.,.)} 
\end{eqnarray} 
where $\bm{I}_n$ is the identity matrix of size $n \times n$ with $\bm{0}'s$ of appropriate sizes. We use the \textit{Adam} optimizer to minimize the LOOCV MSE criterion in our implementations.   

\subsection{Uncertainty Quantification}
\label{sec: UQ}

Recall that the minimal attainable MSE is the expected conditional variance and for future reference we omit the term `expected' and assume this implicitly. This is the case when using the conditional expectation as predictor. Assuming $Z$ to be a Gaussian process is a very specific choice since the best linear predictor coincides with the conditional expectation and consequently the minimal MSE is the conditional variance. It is important to note that this is not always the case. The best linear predictor is, generally speaking, a sub-optimal predictor. Since, we refrain from making the Gaussian process prior assumption throughout the paper in favor of generalizability, the minimal MSE obtained with the best linear predictor is at best, an approximation to the conditional variance. To quantify the distance from conditional variance, denote the complete given data $\setD = (\vecZ, \transpose{\matU}\vecZ^+)$ (where `$+$' $= `\star$' in the case of Lagrangian Kriging), the sub-optimal  best linear predictor as $\tau = \tau(\setD)$ and consider the definition of $\Var{Z\etoile \vert\setD}$,  
\begin{eqnarray}
\Var{Z\etoile \vert\setD} &=& \Esp{(Z\etoile - \Esp{Z\etoile\vert\setD})^2 ~ \vert ~\setD} \\
&=& \Esp{(Z\etoile - \tau + \tau - \Esp{Z\etoile\vert\setD})^2 ~ \vert ~\setD} \\
&=& \Esp{(Z\etoile - \tau)^2 + 2(Z\etoile - \tau)(\tau - \Esp{Z\etoile\vert\setD}) + (\tau - \Esp{Z\etoile\vert\setD})^2  ~ \vert ~\setD} 
\end{eqnarray}
where $(Z\etoile - \tau)(\tau - \Esp{Z\etoile\vert\setD})$ can be expanded as $(Z\etoile - \Esp{Z\etoile\vert\setD} + \Esp{Z\etoile\vert\setD} - \tau)(\tau - \Esp{Z\etoile\vert\setD}) = (Z\etoile - \Esp{Z\etoile\vert\setD})(\tau - \Esp{Z\etoile\vert\setD}) - (\tau - \Esp{Z\etoile\vert\setD})^2$ and since $(\tau - \Esp{Z\etoile\vert\setD})$ is a function of $\setD$, the orthogonality of the residual implies that the first term would be zero. This simplifies the expression to,
\begin{eqnarray}
\Var{Z\etoile \vert\setD} 
&=& \Esp{(Z\etoile - \tau)^2 ~ \vert ~\setD} - (\tau - \Esp{Z\etoile\vert\setD})^2 
\end{eqnarray} 
where we use $\Esp{(\tau - \Esp{Z\etoile\vert\setD})^2 \vert\setD} = (\tau - \Esp{Z\etoile\vert\setD})^2$ since it is a function of $\setD$. On the right hand side, the first term can be recognized as the minimal MSE (MMSE). The second term represents the distance of the best linear predictor from the conditional expectation, the best predictor. The conventional approach is to use MMSE as a close approximation to the conditional variance itself allowing the construction of confidence intervals based on the Gaussian assumption or more generally, using the Vysochanskij–Petunin inequality \cite{Chils2012}. So far, however, there has been little discussion addressing this correction term. A possible resolution can be to use Cross-validation based approaches to quantify the MSE and use it as a measure of the correction term. We can also try a polynomial basis expansion of the conditional expectation and compute the contribution of quadratic and higher terms as a proxy for $\tau - \Esp{Z\etoile\vert\setD}$, assuming the linear term coincides with the best linear prediction. %
The main takeaway here is that the MMSE overestimates the conditional variance. Using MMSE as conditional variance is a cautious choice. For our results, we assume the MMSE as the conditional variance and construct $\pm2\sigma$ intervals. 
\\
\noindent In the case of co-Kriging, the MMSE is easy to compute since, it is the same expressions as simple or ordinary Kriging but with the extended matrices $\matK^+$ and $\matH^+$. For simple co-Kriging, 
\begin{equation}
\mathrm{\widehat{Var}}_{CK}\left[\vecZ\etoile\vert \setD\right] = \matK\etoile - \transpose{(\matH^+)} (\matK^+)^{-1} \matH^+
\end{equation}
where $\matK\etoile$ is the matrix $\Cov{\vecZ\etoile, \vecZ\etoile}$. Similarly, in the case of Lagrangian Kriging under the simple Kriging assumption, we can solve for the MMSE,
\begin{equation}
\mathrm{\widehat{Var}}_{LK}\left[\vecZ\etoile\vert \setD\right] = \matK\etoile - \transpose{(\matH + \vecZ \transpose{\veclambda} \transpose{\matU})} \matK^{-1} (\matH - \vecZ \transpose{\veclambda} \transpose{\matU})
\end{equation}
where $\veclambda = (\transpose{\vecZ} \matK^{-1} \vecZ)^{-1} (\transpose{\matU} \matU)^{-1} (\vecv - \transpose{\matU}\transpose{\matH}\matK^{-1}\vecZ)$ and the definitions of matrices $\matK, \matH, \matU$ and vector $\vecv$ are derived from Section \ref{sec: ConstrainedKriging}. In the case of 2D flow prediction, we can compute the uncertainty in the $x-y$ components of velocity separately. To quantify uncertainty of the squared magnitude however, is not that straightforward. Fortunately, the potential function formulation explicitly quantifies cross-covariances between velocity components. Under the normal distribution assumption for both velocity components, the squared sum is characterized by a generalized $\chi^2$ distribution, for which we compute all parameters and from the explicit covariance relations.

\section{Results}
\label{sec: results}

We present some numerical experiments to illustrate the performance on a synthetic ODE problem, and a synthetic PDE problem before addressing an application in fluid dynamics. All experiments have been implemented on a local machine with \texttt{Intel(R) Core(TM) i7-5600U CPU @ 2.60GHz}. 
To compute the derivatives of the covariance function we have tried two ways: automatic differentiation via \texttt{pytorch} \cite{pytorch} and analytic expressions (upto order $4$). Analytic expressions are the quickest to implement but tedious to code. 
For the results presented below, we use the \texttt{pytorch} implementation in 1D and for quicker computation, we use the analytical expressions in 2D (see appendix \ref{sec: kernelderiavtives}).  

\subsection{Synthetic ODE}
\label{sec: ODE}
Consider the ODE $f(x) + f''(x) = 0$. We have also covered this as part of the working Example \ref{example: CKWorkingExample} when describing the two methods. The problem is to estimate $f(x\etoile)$ given,
\begin{eqnarray*}
f(x_i) &=& y_i \quad \text{for $i = 1,\ldots ,n$} \\
f(x_j) + f''(x_j) &=& 0 \quad \text{for $j = n+1,\ldots, n+p$}  
\end{eqnarray*}
where $x_{[1:n+p]} \in [0, 2\pi]$. It is straightforward to see that all functions of the form $f(x) = a \cos (x) + b \sin (x)$ with $a,b \in \R$ will satisfy $f(x) + f''(x) = 0$. Imposing boundary condition, $f(0) = 0$ narrows it down to all functions of the form $a\sin(x)$. One extra observation in the domain would imply a unique solution. We start with the same number of collocation points $p=10$ and prediction points $q=10$ to allow for a fair comparison of simple Kriging, co-Kriging and Lagrangian Kriging. Additionally, we provide $n=4$ random observations of the function $f(x) = \sin(x)$.
\begin{figure}[H]
\begin{center}
\begin{tabular}{cc}
\centering
\includegraphics[trim = {0 0 0 0.5cm}, clip, width=0.3\textwidth]{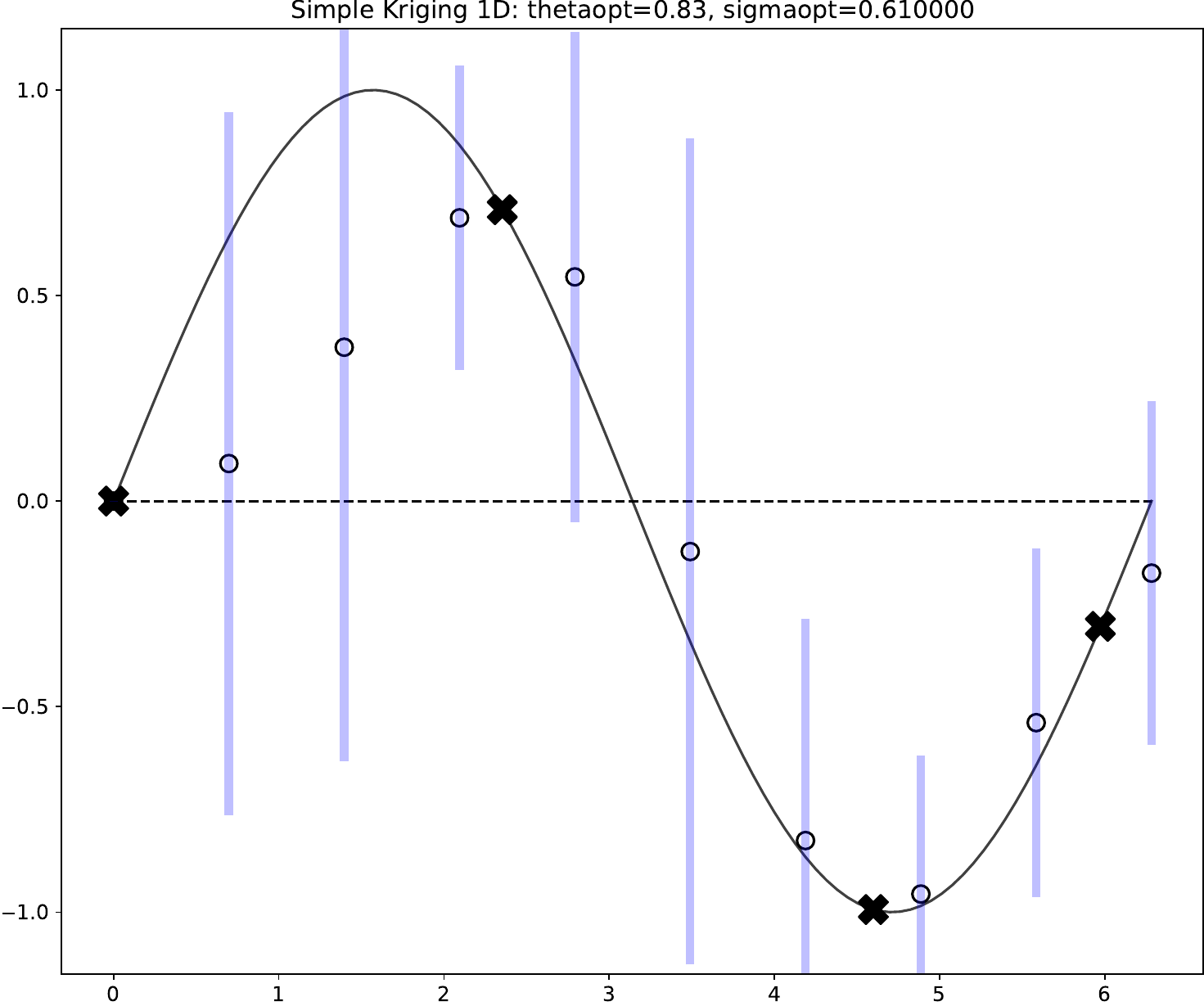} & \includegraphics[trim = {0 0 0 0.5cm}, clip,width=0.3\textwidth]{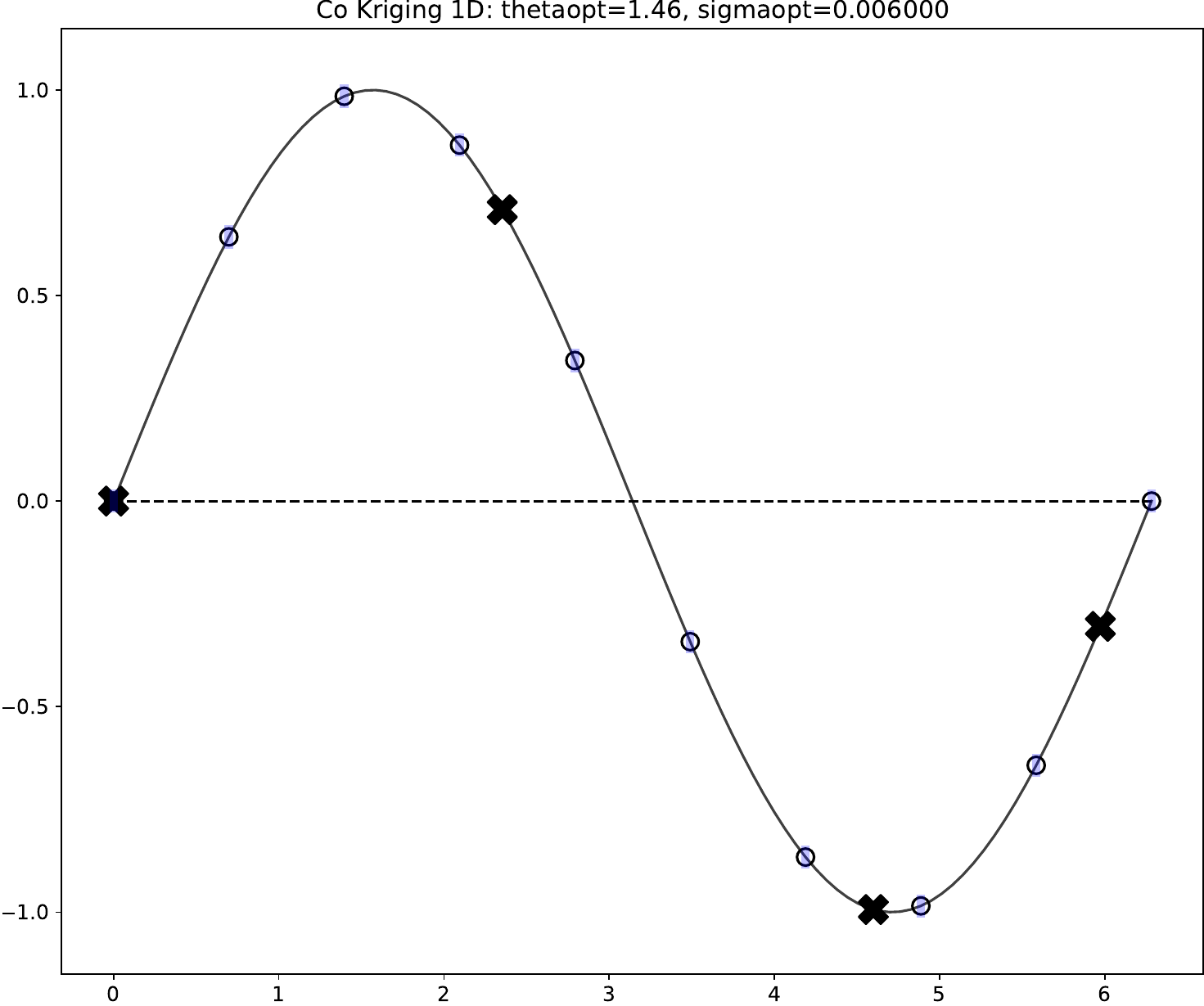} \\
\textbf{(a)} Simple Kriging & \textbf{(b)} Collocated co-Kriging \\[1em]
\includegraphics[trim = {0 0 0 0.5cm}, clip,width=0.3\textwidth]{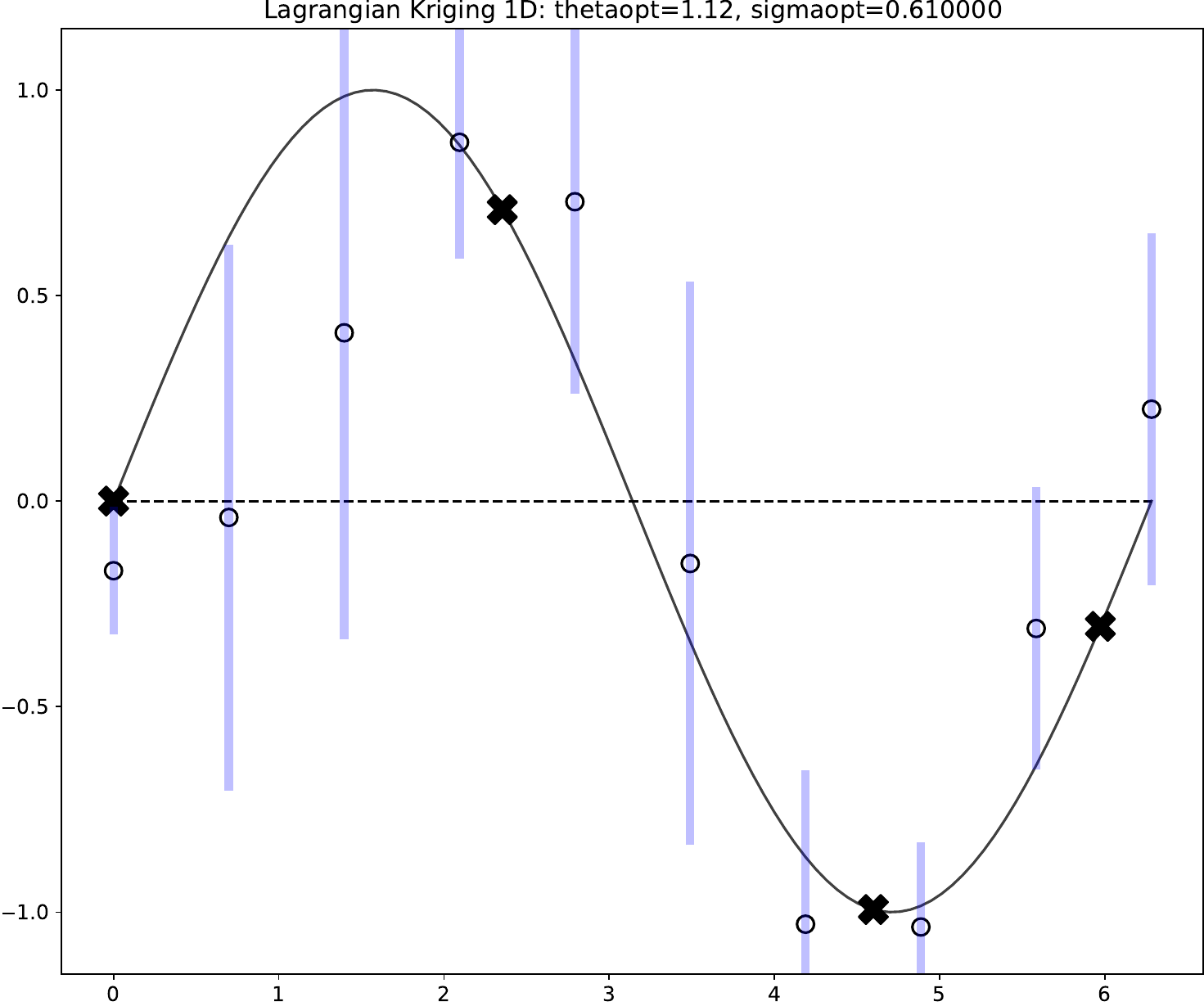} & \includegraphics[trim = {0 0 0 0.5cm}, clip,width=0.3\textwidth]{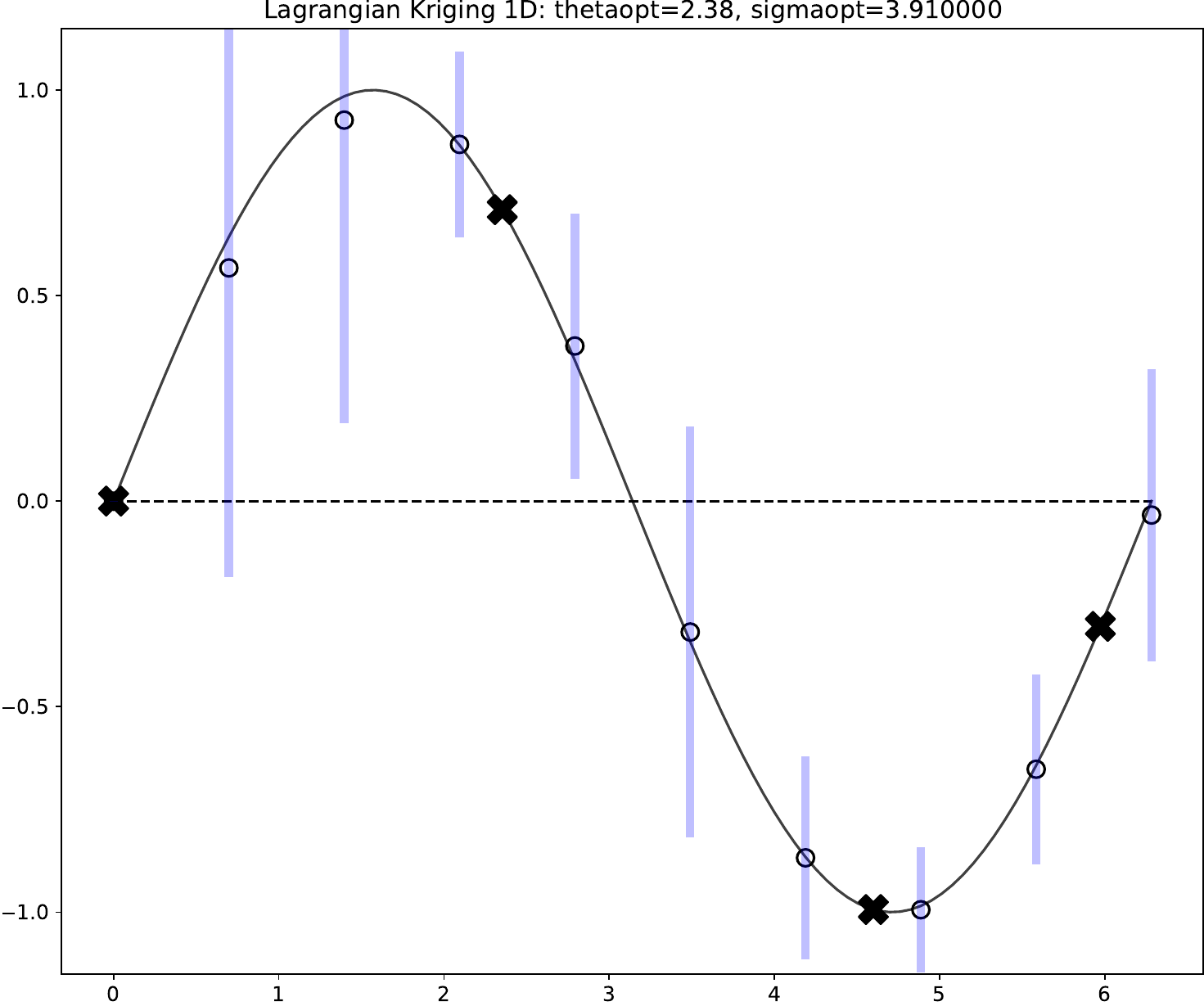} \\
 \textbf{(c)} Lagrangian Kriging (LOOCV) & \textbf{(d)} Lagrangian Kriging (improved interpolation)
\end{tabular}
\end{center}
\caption{All figures show the exact function in black, four observations as black crosses, collocation points / the predicted $f$ as empty circles with $\pm 2\sigma$ confidence intervals: \textbf{(a)} Simple Kriging ($\hat{\theta} = 0.83,~ \hat{\sigma} = 0.61$). \textbf{(b)} Co-Kriging ($\hat{\theta} = 1.46, ~\hat{\sigma} = 0.006$) \textbf{(c)} Lagrangian Kriging (LOOCV optimal parameters: $\hat{\theta} = 1.12, ~\hat{\sigma} = 0.61$.)  \textbf{(d)} Lagrangian Kriging (minimizer of the interpolation error: $\hat{\theta} = 2.38, ~\hat{\sigma} = 3.91$). Note that, the co-Kriging uncertainty intervals are too tiny to notice owing to the model accuracy.}
\label{fig: ODEcomparison}
\end{figure}

\noindent From the Figure \ref{fig: ODEcomparison}, derivative observations based co-Kriging is the clear winner in terms of accuracy. The Lagrangian Kriging approach has two variants. The first one, Figure \ref{fig: ODEcomparison}c, uses the LOOCV optimal hyperparameters  and notably loses the interpolation property which is a distinctive property of classical Kriging. To resolve this, in the second variant, Figure \ref{fig: ODEcomparison}d, we tune hyperparameters using the same criterion formulas as LOOCV (MSE and variance criterion) except that we use all observations while predicting (not leaving any out). In this scenario, the MSE in particular, corresponds to the interpolation error caused by the prediction constraints and provides a better estimate (interpolating) as is evident from Figure \ref{fig: ODEcomparison}d. Both methods improve model uncertainty when compared to simple Kriging. In Table \ref{tab: Performance1D}, we substantiate again that co-Kriging is undoubtedly the best for this particular task. Simple Kriging and Lagrangian Kriging with LOOCV parameters render a similar performance even though simple Kriging performs marginally better in this case. It is also noteworthy that Lagrangian Kriging with minimizer of the interpolation error is significantly better than the LOOCV variant. This is to stress on the fact that, nudging the parameters slightly or opting for modifying the hyperparameter tuning criterion can result in a major improvement in the Lagrangian Kriging performance with no computational overhead.  
\newline\newline
\begin{table}[H]
\small
\caption{Comparing the models based on the usual MSE mertic for the 1D prediction task. For Lagrangian Kriging, we report results for both variants, the LOOCV and the improved interpolation parameters.} \vspace{1em}
\label{tab: Performance1D}
\centering
\begin{tabular}{c c c c}
\toprule 
Simple Krig. & Co-Krig. &  Lag. Krig. $(\theta = 1.12)$ & Lag. Krig. $(\theta = 2.38$) \\
\midrule%
$0.08$ & $6.75 \times 10^{-11}$ & $0.12$ & $1.22 \times 10^{-3}$\\
\cline{1-4}
\end{tabular}
\end{table}

\noindent We also present the run times, in 1D, for varying $p$ in Table \ref{tab: CPUtimes} to further motivate the utility of Lagrangian Kriging. The computational time is divided into two parts, the time of covariance computations to construct all covariance matrices (construction) and the time to invert and multiply the covariance matrices (inversion) to obtain the final predictions. The maximum $p$ considered is $1000$ and matrix inversions are tractable. In all cases, the major bottleneck comes from the number of covariance calculations required to fill the covariance matrices. This is most likely a consequence of the dimension. The construction time for co-Kriging is significantly higher owing to the expanded covariance matrices, $\matK^{+}$ and $\matH^{+}$. In comparison, Lagrangian Kriging offers a speed-up of 100x - 1000x in several instances.   

\begin{table}[H]
\small
\caption{Computational times for varying $p$ (without accounting for kernel hyperparameter calibration). For Lagrangian Kriging, we consider $q = p$.} \vspace{1em}
\label{tab: CPUtimes}
\centering
\begin{tabular}{c c c c}
\toprule 
$n = 4$ &  &  Co-Krig. $(q=100)$ & Lag. Krig. \\
\midrule%
\multirow{2}{*}{\small $p = 100$} & \footnotesize{Construction} & 32.04s & 0.26s \\ 
& \footnotesize{Inversion} & 0.07s & 0.09s \\
\cline{1-4}
\multirow{2}{*}{\small $p = 250$} & \footnotesize{Construction} & 179.08s & 0.59s \\ 
& \footnotesize{Inversion} & 0.13s & 0.16s \\
\cline{1-4}
\multirow{2}{*}{\small $p = 500$} & \footnotesize{Construction} & 686.12s & 1.17s \\ 
& \footnotesize{Inversion} & 0.26s & 0.20s \\
\cline{1-4}
\multirow{2}{*}{\small $p = 1000$} & \footnotesize{Construction} & 2680.72s & 2.34s \\ 
& \footnotesize{Inversion} & 0.39s & 0.56s \\
\cline{1-4}
\end{tabular}
\end{table}
\begin{figure}[H]
\centering
\begin{tabular}{c c}
\includegraphics[trim = {0 0 0 0.5cm}, clip,width=0.45\textwidth]{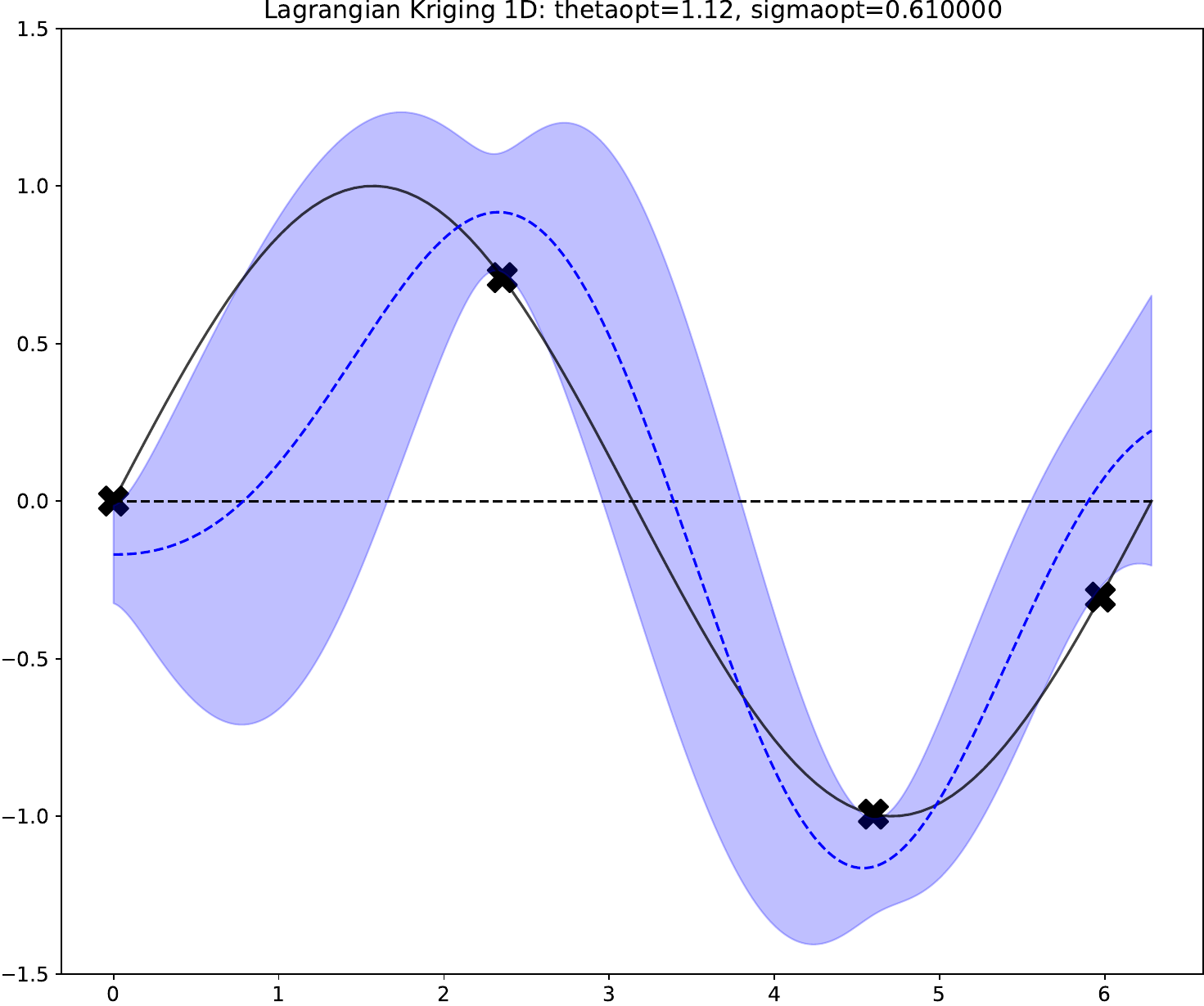} &
\includegraphics[trim = {0 0 0 0.5cm}, clip,width=0.45\textwidth]{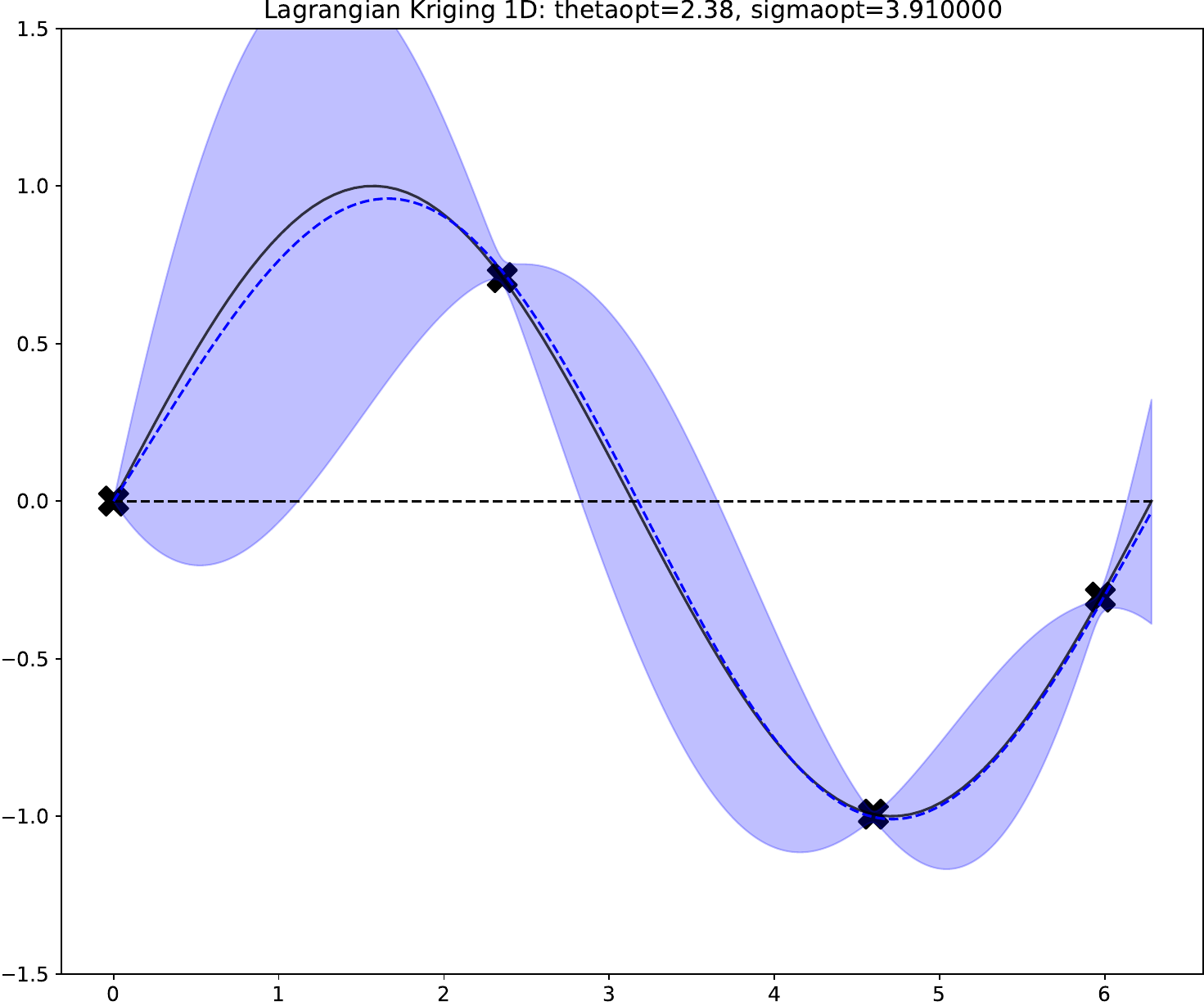} \\
(a) LOOCV & (b) Improved interpolation \\
\end{tabular}
\caption{Lagrangian Kriging with the same observations but predicting at $q=1000$ locations. (a) Using the optimal LOOCV parameters. (b) Using parameters minimizing the interpolation error.}
\label{fig: CKn4p1000}
\end{figure}

\begin{figure}[H]
   \centering
   \begin{tabular}{ccc}
   \includegraphics[width=0.3\textwidth]{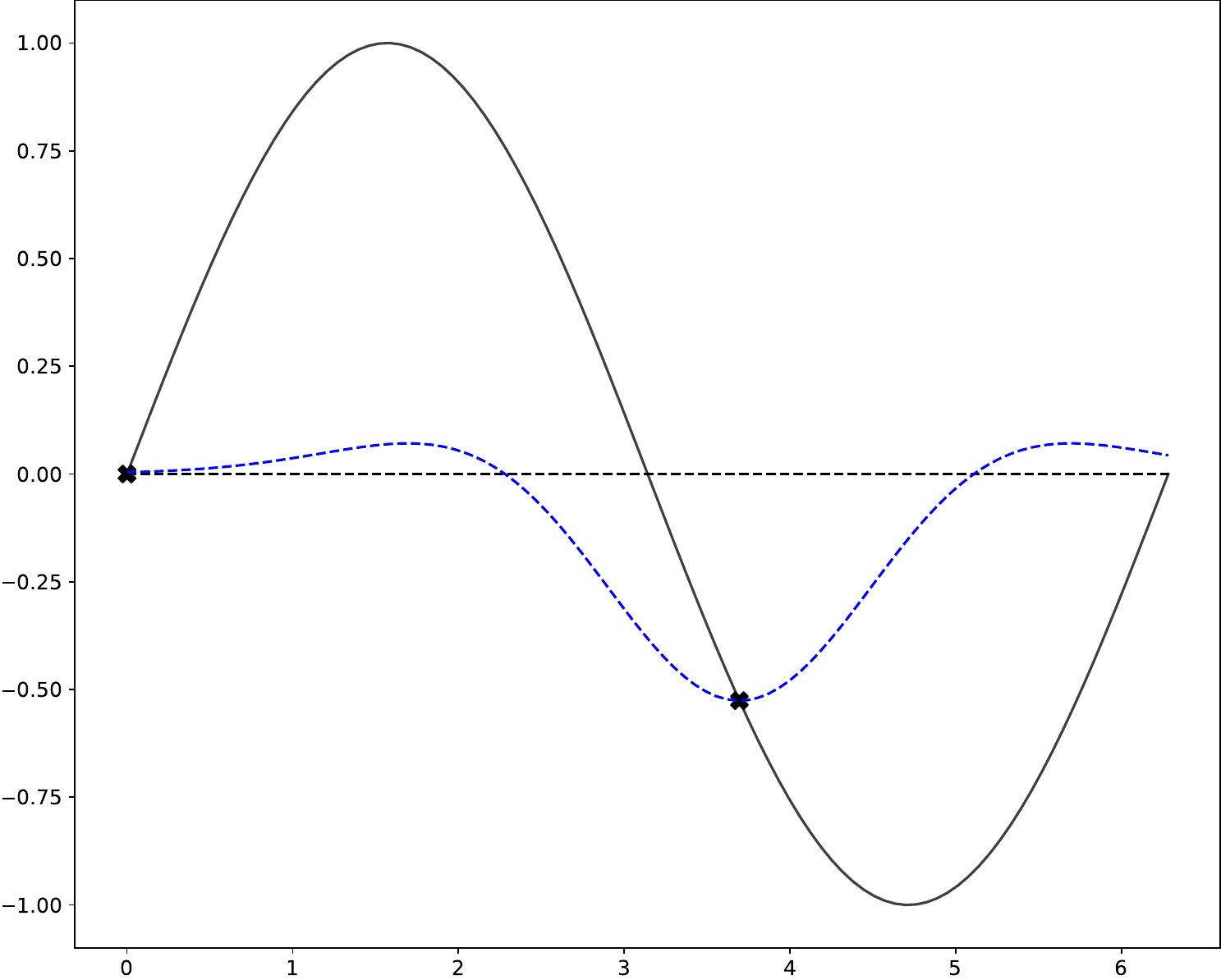} & \includegraphics[width=0.3\textwidth]{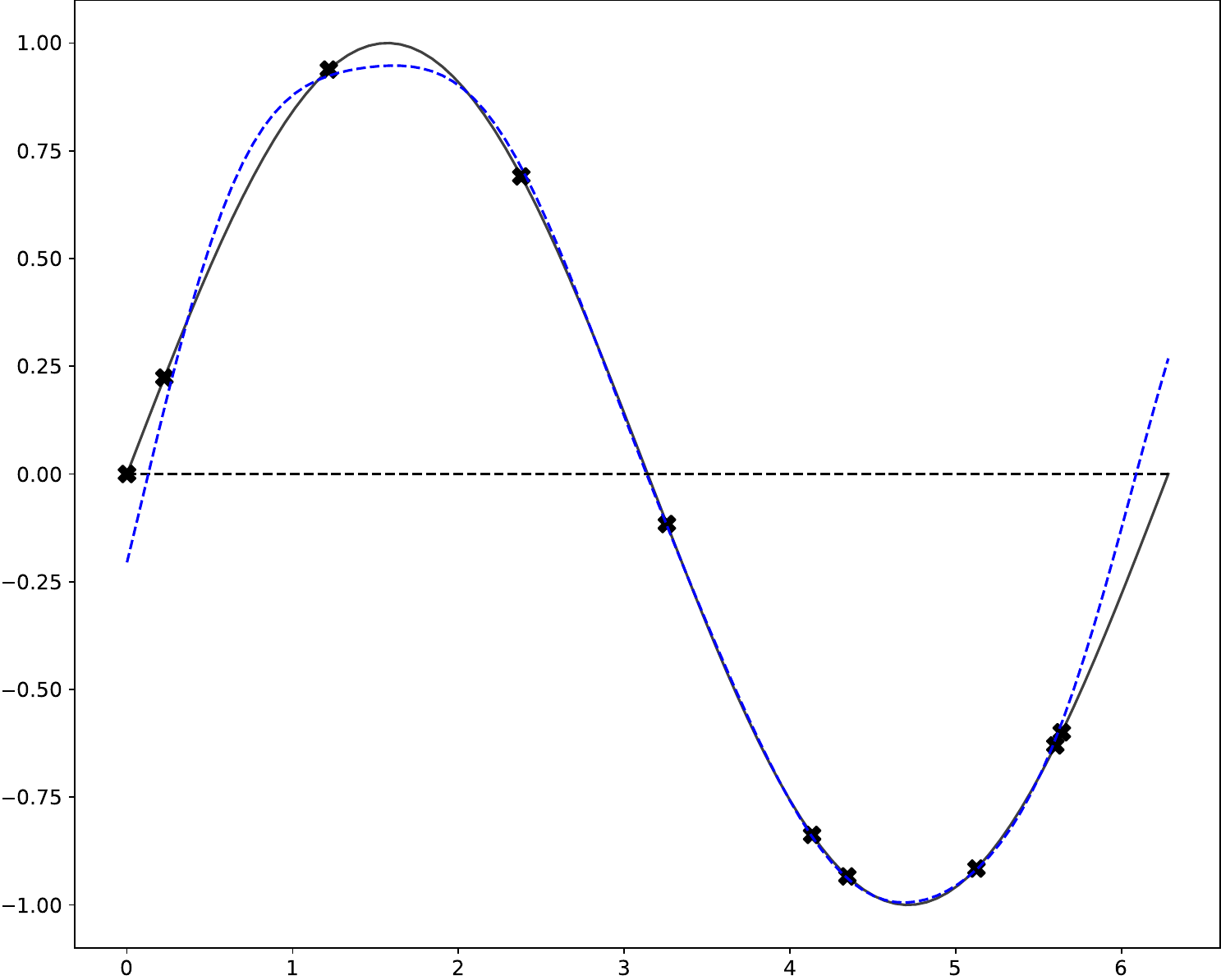} & \includegraphics[width=0.3\textwidth]{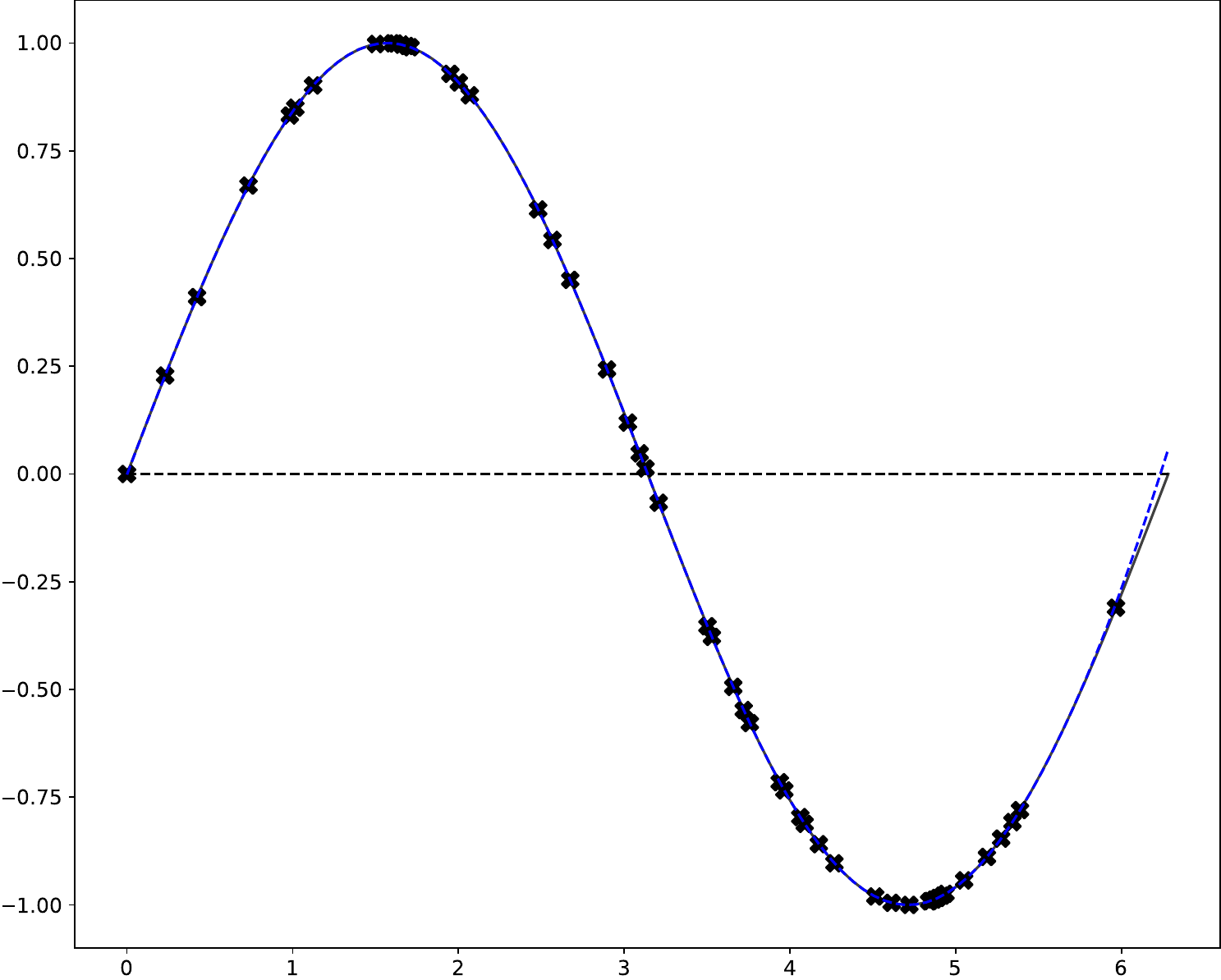} \\
   $n=2$ & $n=10$ & $n=50$
   \end{tabular}
   \caption{The effect of increasing observations on Lagrangian Kriging prediction (in blue) without uncertainty for $q=100$ prediction points. Fixed $\theta = 1$.}
   \label{fig: IncreasingObsLK}
   \end{figure}

Before we study exeriments in higher dimensions, we present some results specifically on Lagrangian Kriging since it loses the crucial, interpolation property. In Figure \ref{fig: ODEcomparison}, we chose $q = 10$ as the number of prediction locations. On augmenting $q$, Figure \ref{fig: CKn4p1000} illustrates that the prediction does not vary but it becomes smoother since we chose a finer set of prediction locations. Turning to the effect of the number of exact observations $n$, Figure \ref{fig: IncreasingObsLK} displays that Lagrangian Kriging converges to the exact function with increasing $n$. This is consistent with our expectations, since Kriging has traditionally been used as an interpolation model and given sufficient exact observations, it should converge to the Kriging prediction. Note that instead of performing LOOCV or minimal interpolation error hyperparameter tuning, we used a fixed lengthscale to better assess the effect of observations in isolation.

\subsection{2D scalar functions}
\label{sec: 2D scalar}
In 2D, we consider two real-valued functions, $f_1(x, y) = \cos(x)~\sin(y)$ and a harmonic function $f_2(x, y) = e^x~\sin(y) + 5$. We consider the observations of the functions alongside two sources of differential information, 
\begin{eqnarray}
\vecun \cdot \nabla f_i \coloneqq \left(\frac{\partial}{\partial x} + \frac{\partial}{\partial y}\right) f_i \quad  i=1,2~,\\
\nabla^2 f_i \coloneqq \left( \frac{\partial^2}{\partial x^2} + \frac{\partial^2}{\partial y^2} \right)f_i \quad i=1,2~,
\end{eqnarray} 
where the latter is the well-known Laplacian operator. We compare the model performances for a fixed number of observations $n=10$, $\vecun \cdot \nabla f$ at $p_1 = 50$ equispaced points and $\nabla^2 f$ at $p_2 = 100$ equispaced points. The number of prediction points for derivative observations based collocated simple Kriging, we predict on a grid of $q = 900$ points. Function $f_2$ is harmonic which means that it satisfies $\nabla^2 f_2 = 0$ everywhere in the domain. This is specifically chosen because we intend to apply these methods to fluid dynamics where the continuity equation translates to the Laplace's equation.

\begin{figure}[H]
\centering
\begin{tabular}{ccc}
\includegraphics[trim = {0 0 0 0.5cm}, clip, width=0.3\linewidth]{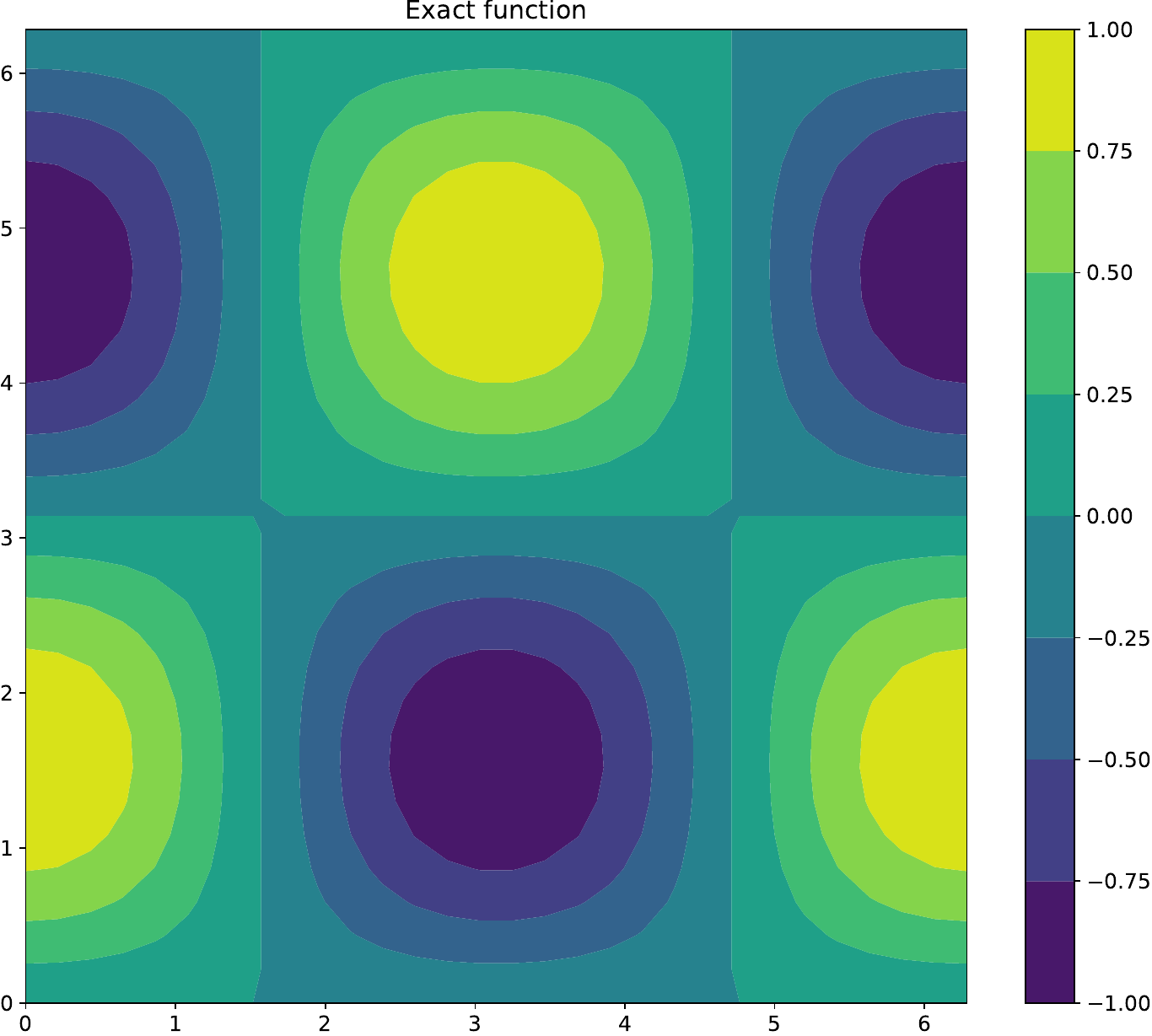} & \includegraphics[trim = {0 0 0 0.5cm}, clip,width=0.3\linewidth]{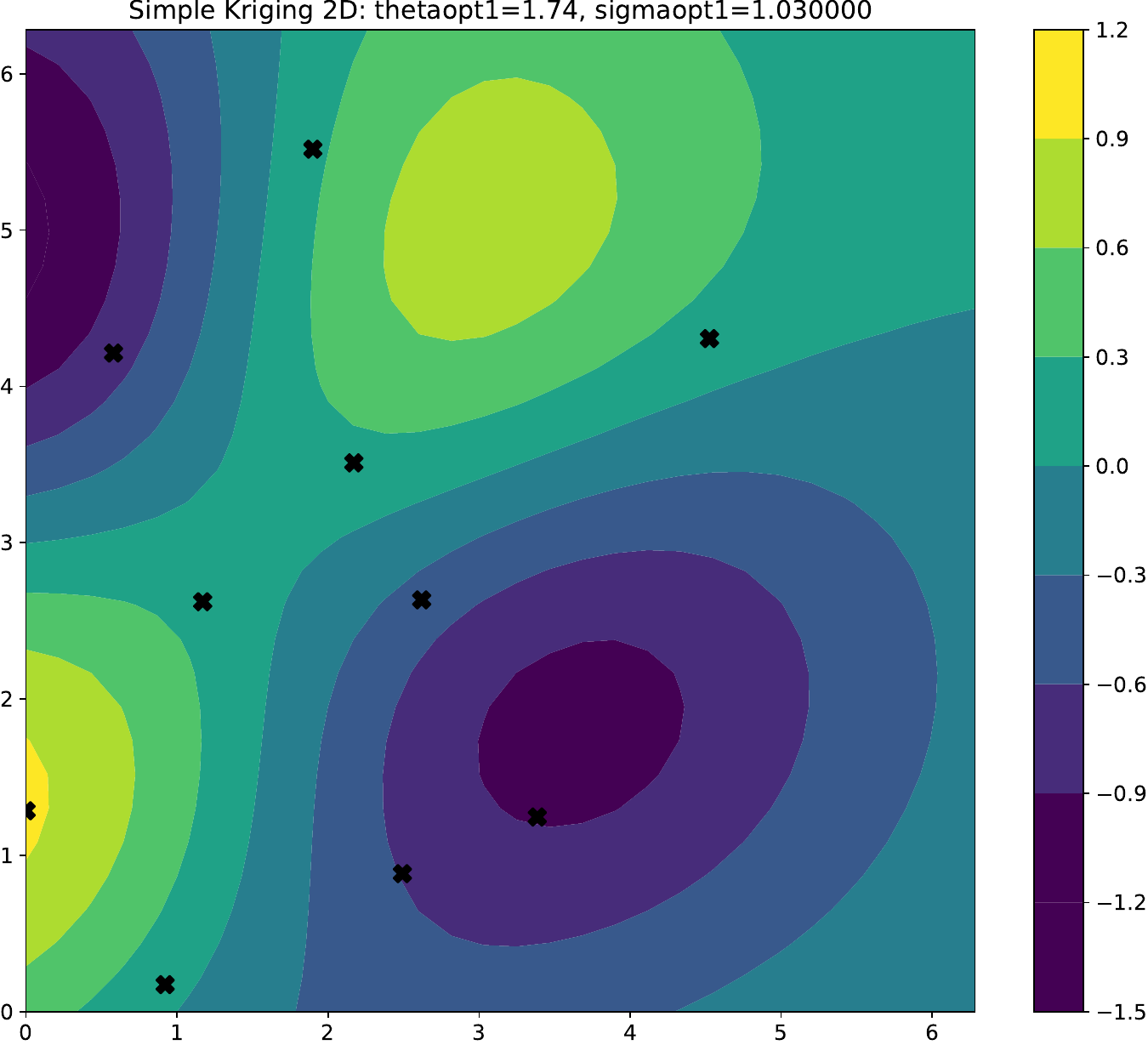} & \includegraphics[trim = {0 0 0 0.5cm}, clip,width=0.3\linewidth]{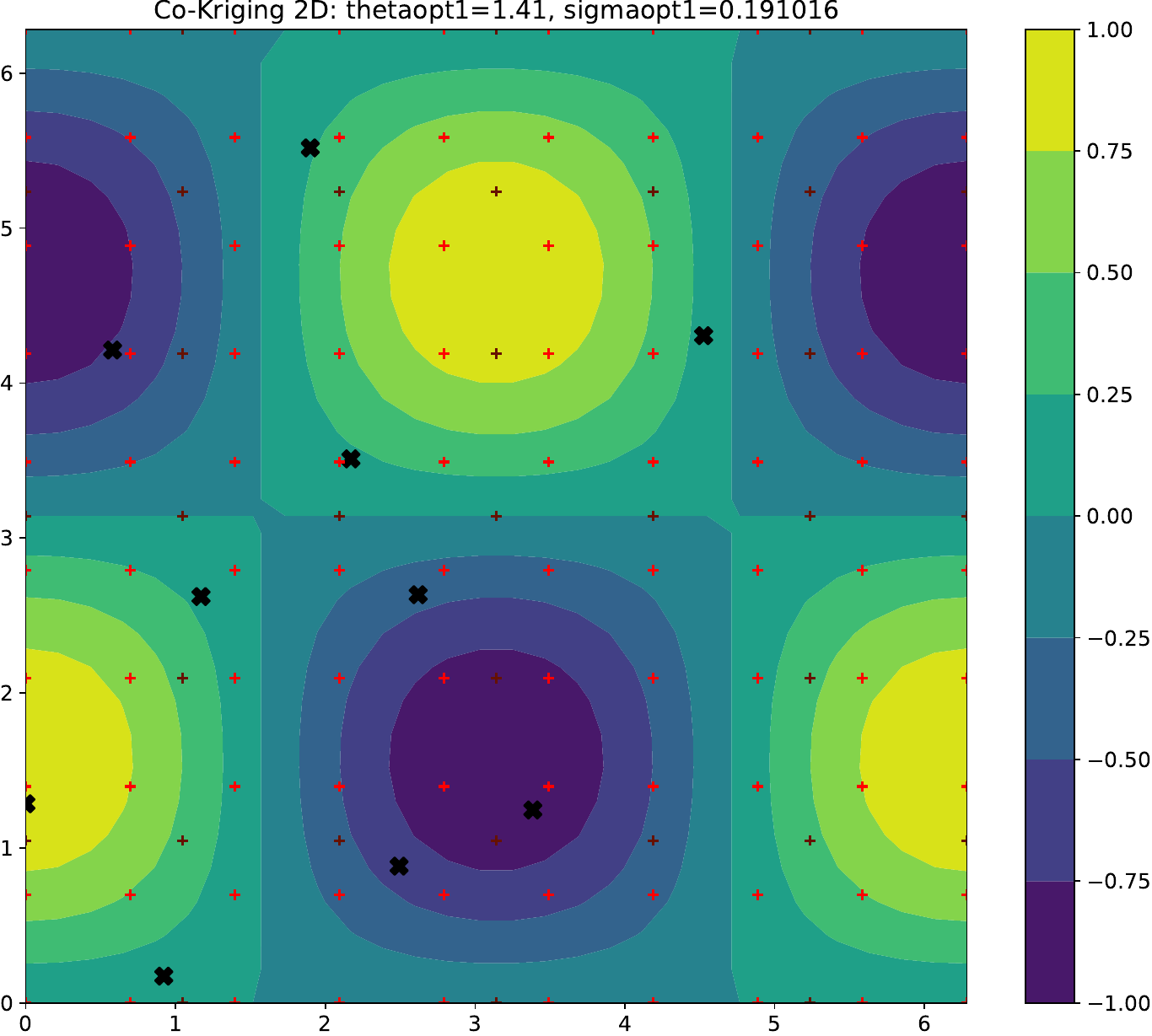} \\
(a) True $f_1(x, y)$ & (b) Simple Kriging, $\hat{\theta} =  1.74$ & (c) Collocated simple co-Kriging, $\hat{\theta} = 1.41$ \\[1em] %
\includegraphics[trim = {0 0 0 0.5cm}, clip,width=0.3\linewidth]{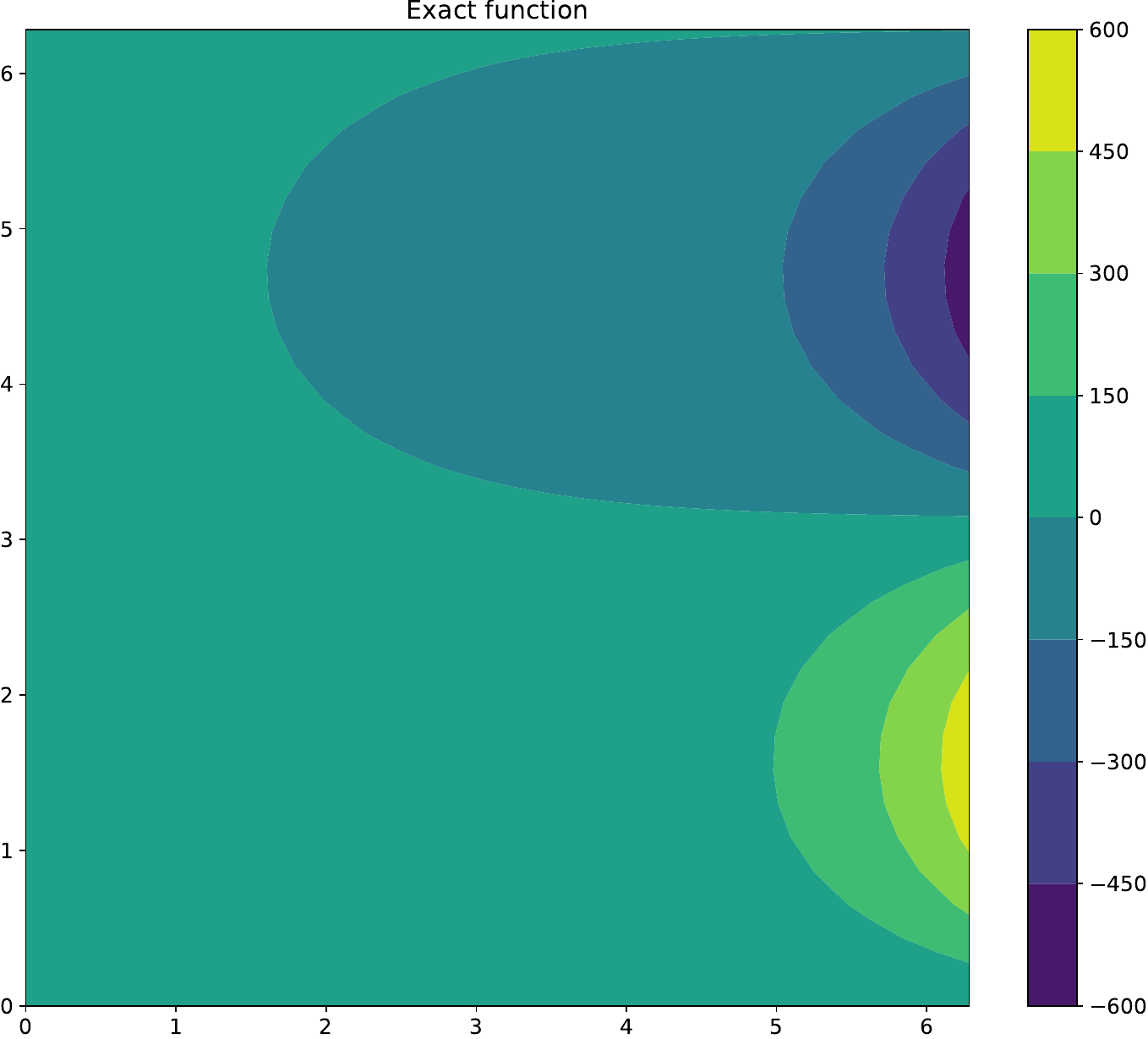} & \includegraphics[trim = {0 0 0 0.5cm}, clip,width=0.3\linewidth]{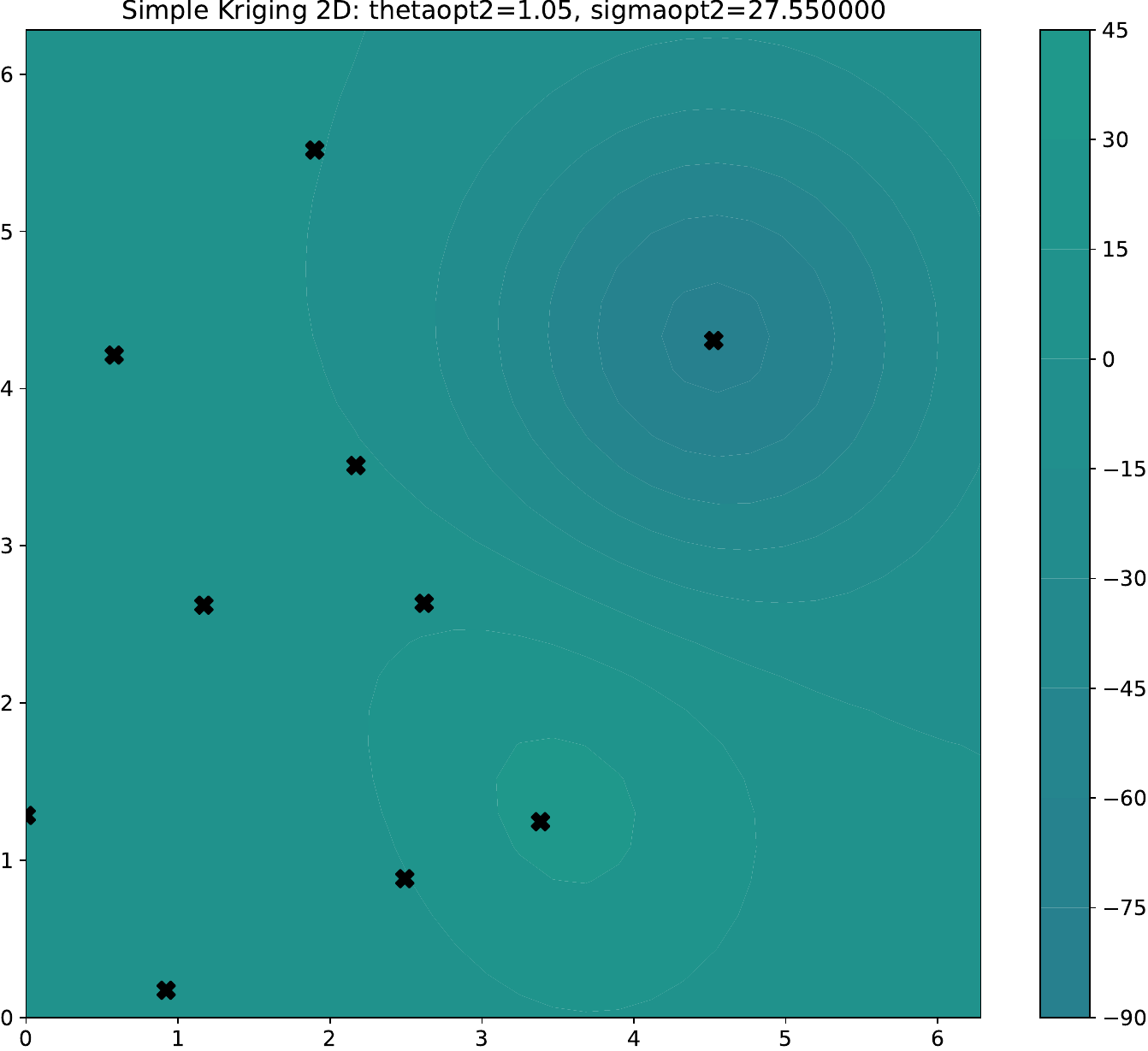} & \includegraphics[trim = {0 0 0 0.5cm}, clip,width=0.3\linewidth]{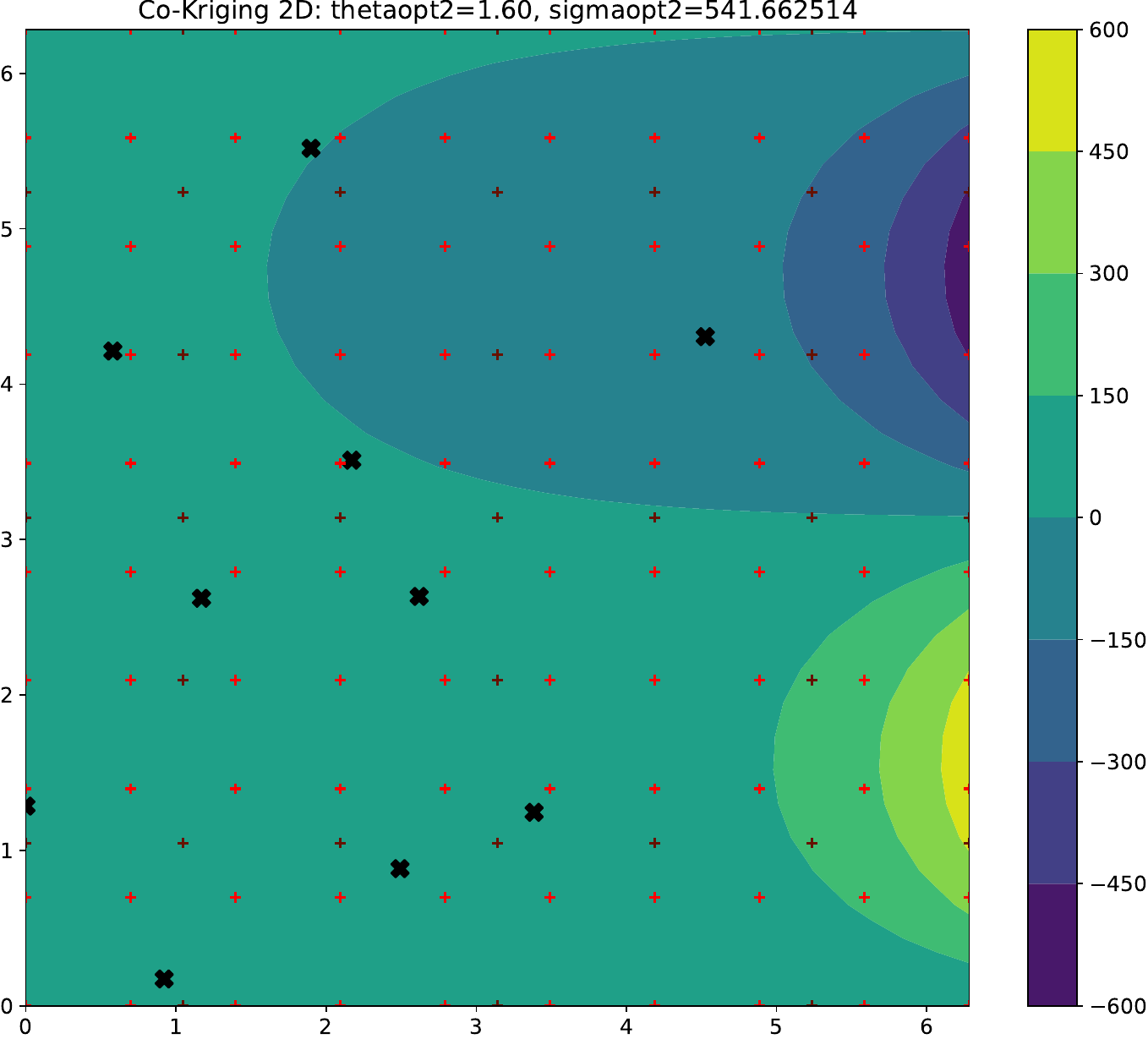} \\
(d) True $f_2(x, y)$ & (e) Simple Kriging $\hat{\theta} = 1.05$ & (f) Collocated simple co-Kriging, $\hat{\theta} = 1.60$. %
\end{tabular}
\caption{Comparing the simple Kriging vs. derivatives based collocated co-Kriging with $n=10$ observations (in black), gradient observations (in brown) and Laplacian observations (in red).}
\label{fig: func2DSKvsDerObs}
\end{figure}

\noindent In Figure \ref{fig: func2DSKvsDerObs}, for both $f_1$ and $f_2$, it is clear that co-Kriging meaningfully utilizes the linear differential information to gain a significant edge in performance over simple Kriging. The reconstruction from just a few exact observations is almost perfect. 

\noindent For Lagrangian Kriging, the constraints do not involve the functions $f_1$ and $f_2$ explicitly. This is the same situation as discussed in Remark \ref{remark: LK} and Remark \ref{remark: LK continued}. The constraints are composed of the first-order derivatives of $f_1, f_2$ or higher. Hence, the objective is independent of any $C2$ constraints in \eqref{eqn:LagrangianKrigOpt}. This implies that, under the absence of the unbiasedness constraint $C1$ as well, we are bound to get the same results as simple Kriging and hence we reuse the simple Kriging LOOCV optimal parameters. 
However, predicting derivatives might be of interest in some applications, for example, if you have limited observations of a function $f$ and you want to approximate the gradient or Hessian for iterative optimization algorithms, given that the function satisfies some differential constraints, it is possible to use the Lagrangian Kriging predictions for the higher-order derivatives involved in the constraints. 
\begin{figure}[H]
\centering
\begin{tabular}{cc}
\includegraphics[trim={0 0 0 1cm}, clip, width=0.45\linewidth]{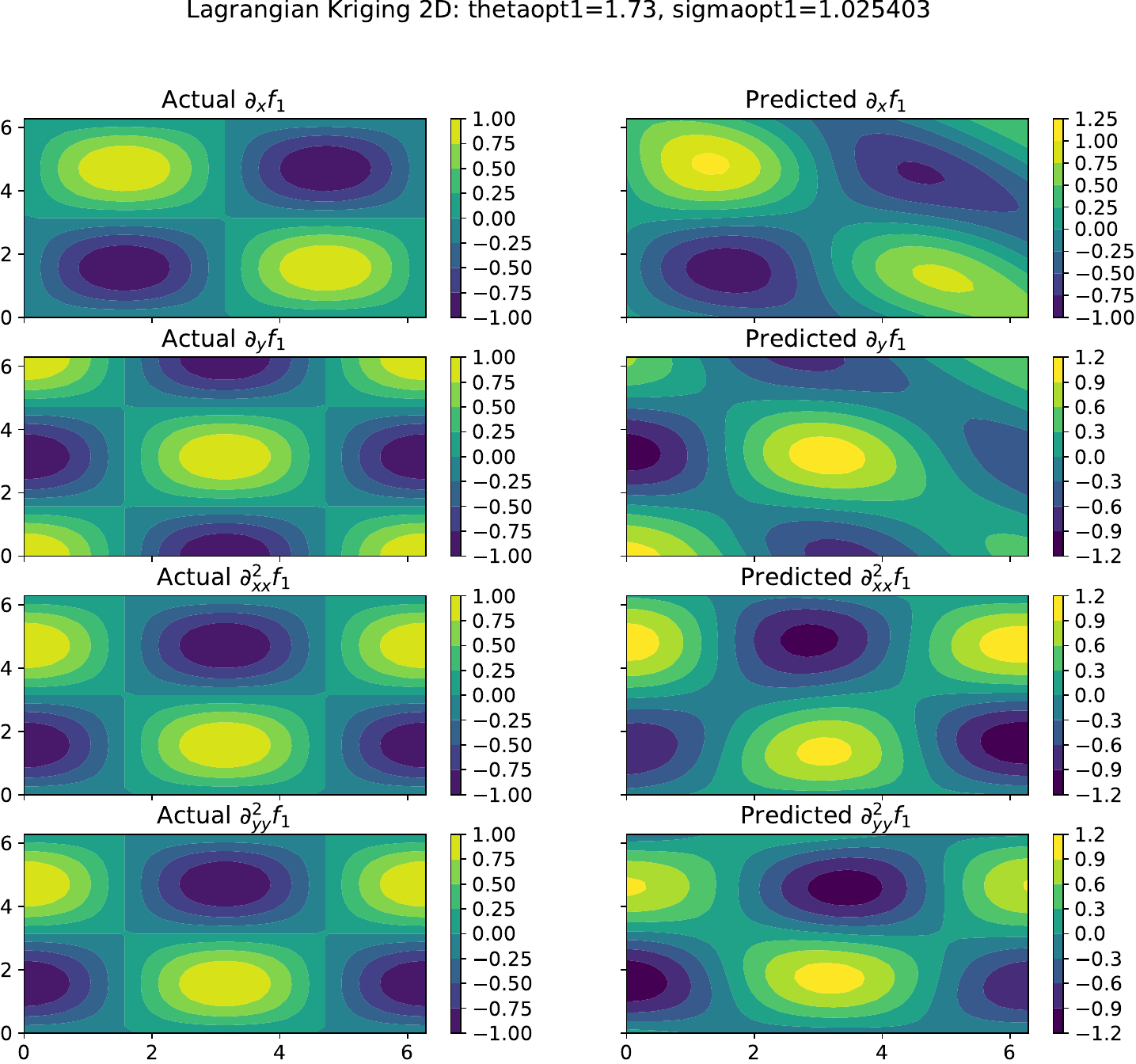} & \includegraphics[trim={0 0 0 1cm}, clip, width=0.45\linewidth]{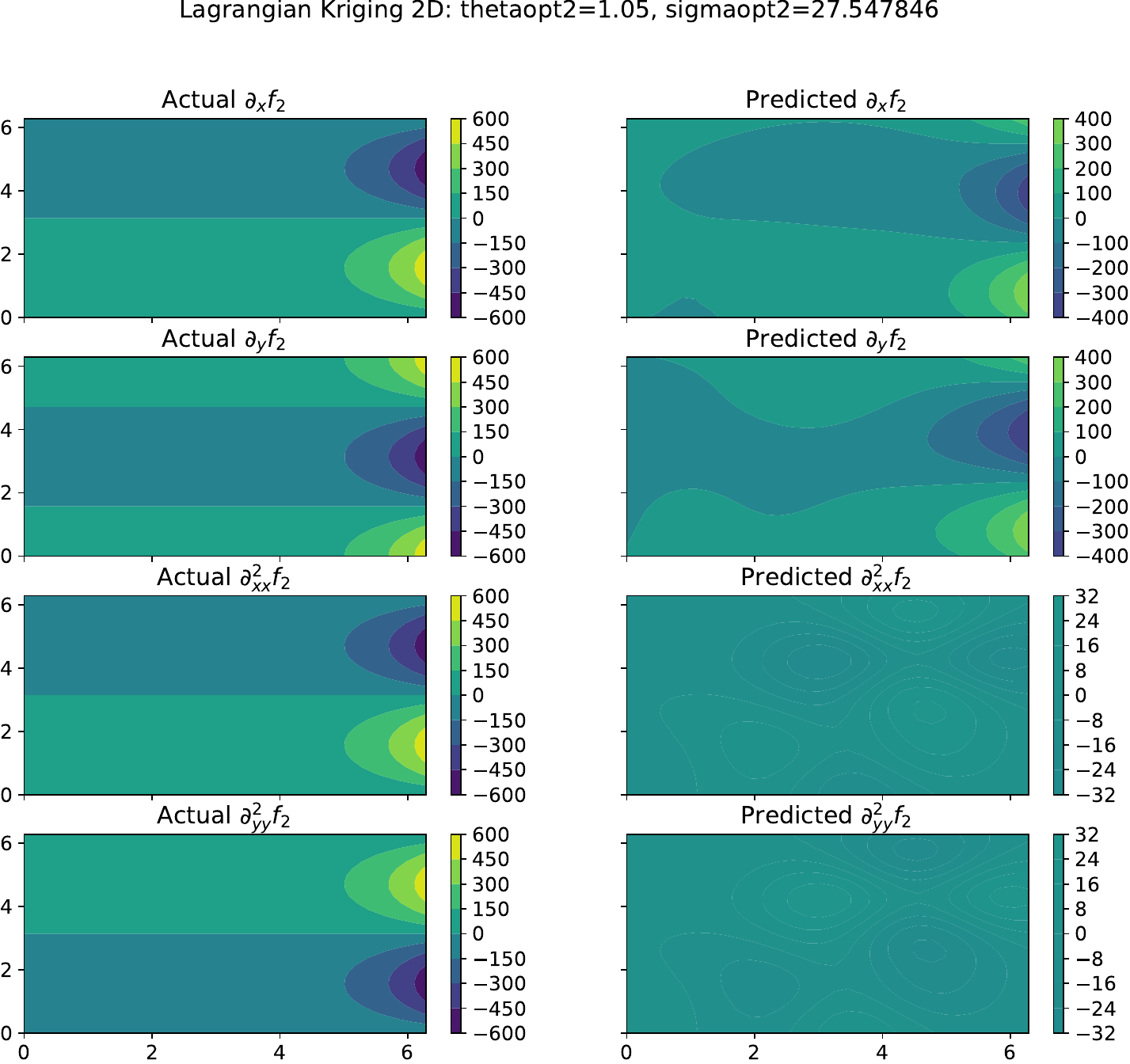} \\ 
(a) $f_1(x,y)$ partial derivatives & (b) $f_2(x,y)$ partial derivatives  
\end{tabular}
\caption{Lagrangian Kriging predictions for higher-order derivatives involved in the gradient-based and Laplacian operators.}
\label{fig: CKforhigherorderderivatives}
\end{figure}

Looking at Figure \ref{fig: CKforhigherorderderivatives}, Lagrangian Kriging produced approximations to the higher derivatives involved in the constraints. As seen before in the case of ODEs, the prediction quality may not be satisfactory especially for the second-order derivatives of $f_2(x,y)$ in Figure \ref{fig: CKforhigherorderderivatives}b. The prediction quality is likely to improve with more exact observations of the function. 
\subsection{Potential flow problem}
\label{sec: Flow prediction}
\begin{figure}[H]
\begin{tabular}{ccc}
\includegraphics[trim = {0 0 0 0.8cm}, clip, width=0.3\linewidth]{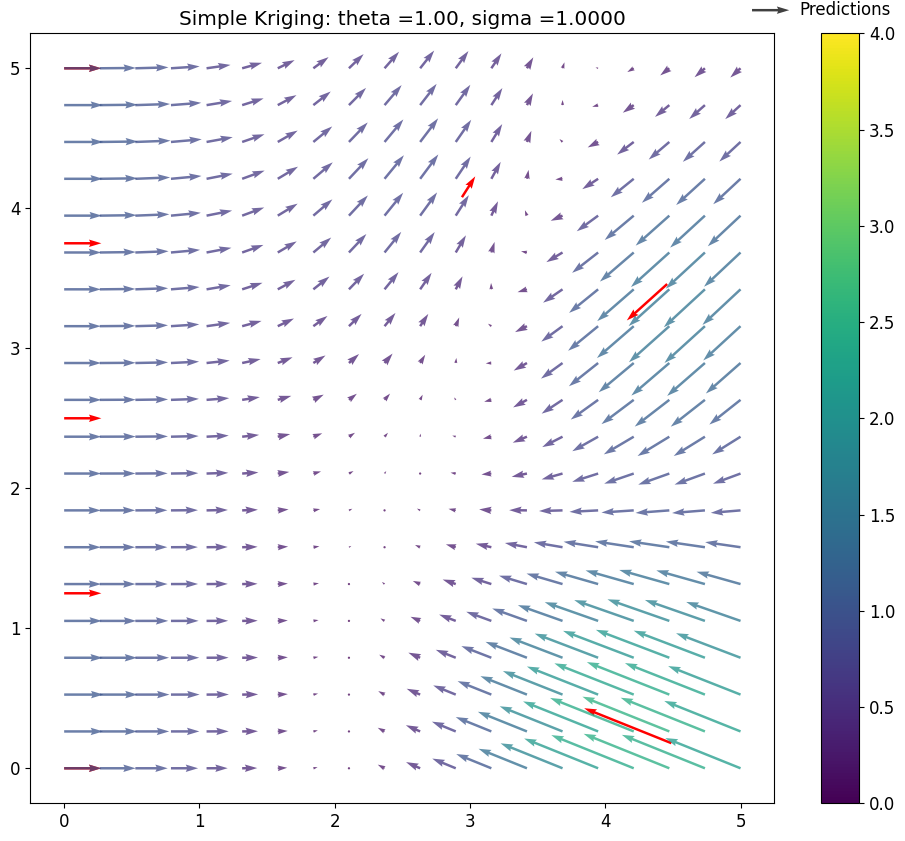} & \includegraphics[trim = {0 0 0 0.8cm}, clip, width=0.3\linewidth]{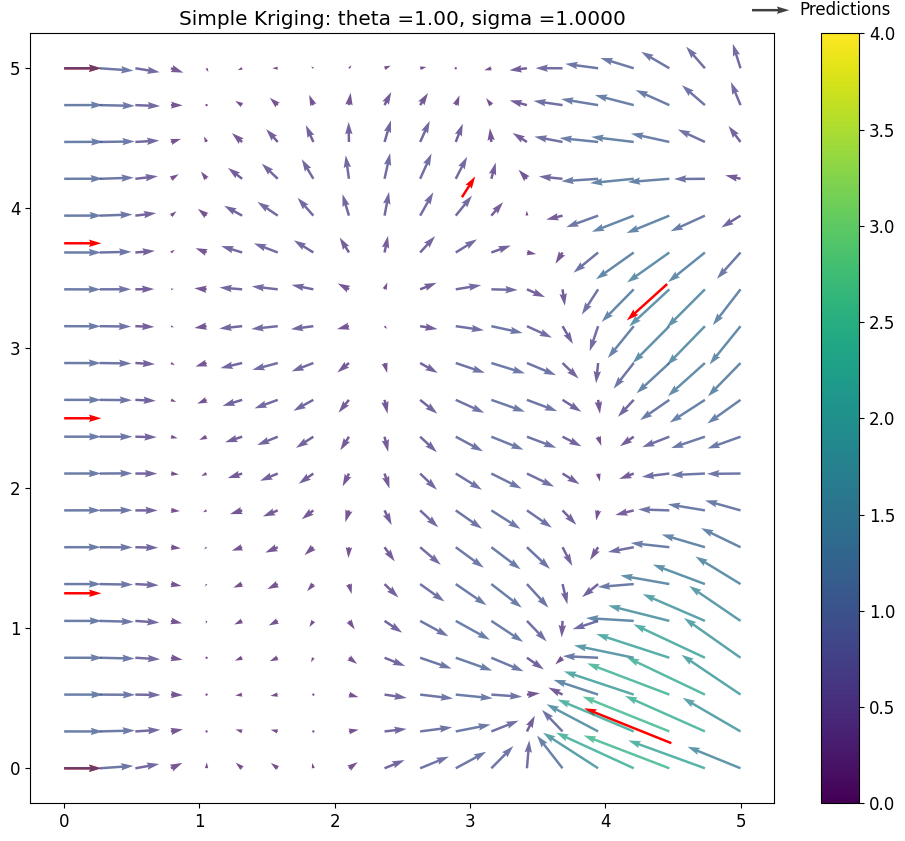} &
\includegraphics[trim = {0 0 0 0.8cm}, clip, width=0.3\linewidth]{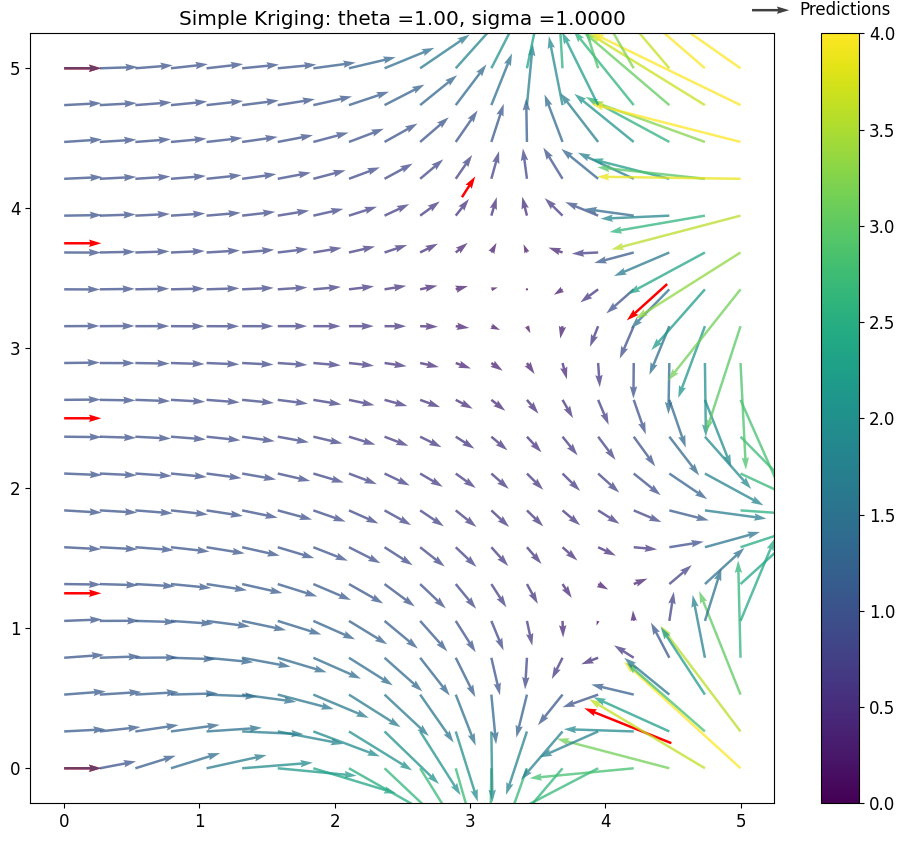} \\
(a) & (b) & (c)
\end{tabular}
\caption{Prediction based on random additional observations in the domain: (a) Simple Kriging without cross covariance between velocity components, (b) Simple Kriging (using $\phi$), (c) Co-Kriging with continuity. Both the color map and the arrow lengths represent the magnitude of velocity.}
\label{fig: flowwithrandomobs}
\end{figure}

\begin{figure}[hbtp]
\centering
\begin{tabular}{ccc}
\includegraphics[trim = {0 0 0 0.8cm}, clip, width=0.3\linewidth]{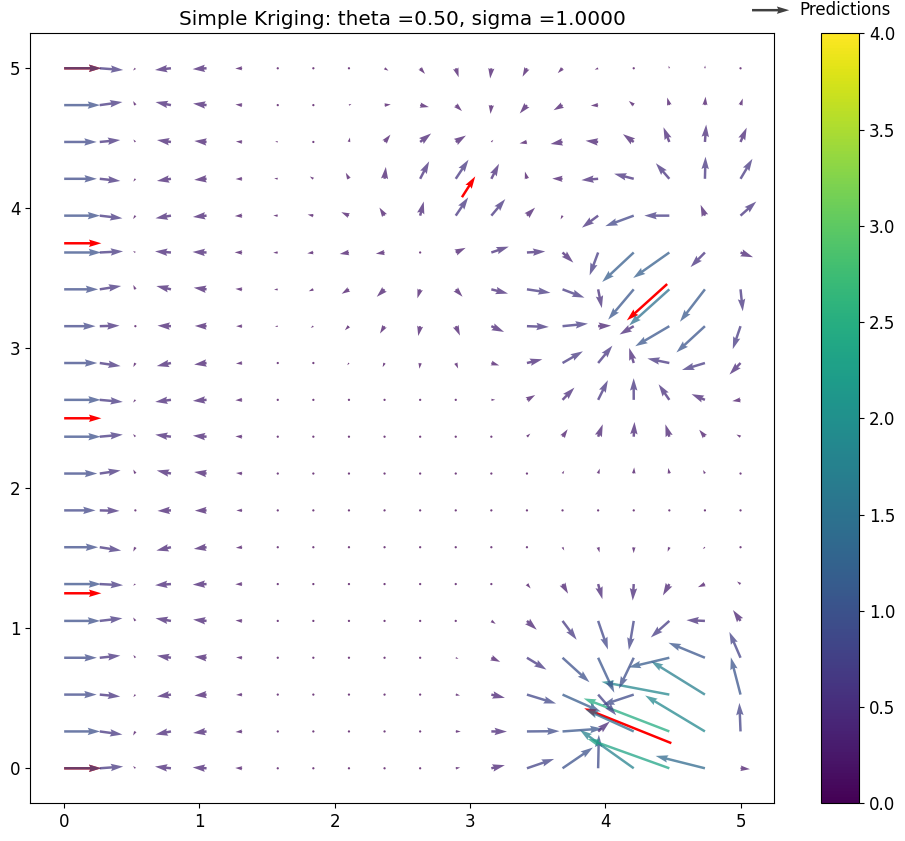} & \includegraphics[trim = {0 0 0 0.8cm}, clip, width=0.3\linewidth]{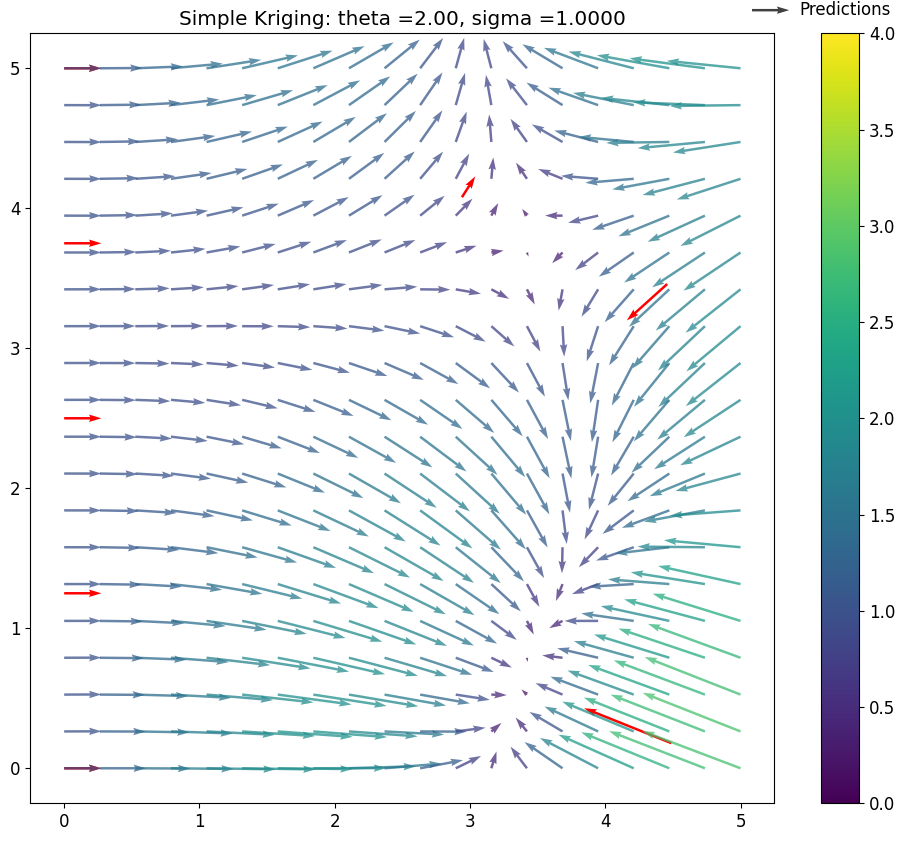} &
\includegraphics[trim = {0 0 0 0.8cm}, clip, width=0.3\linewidth]{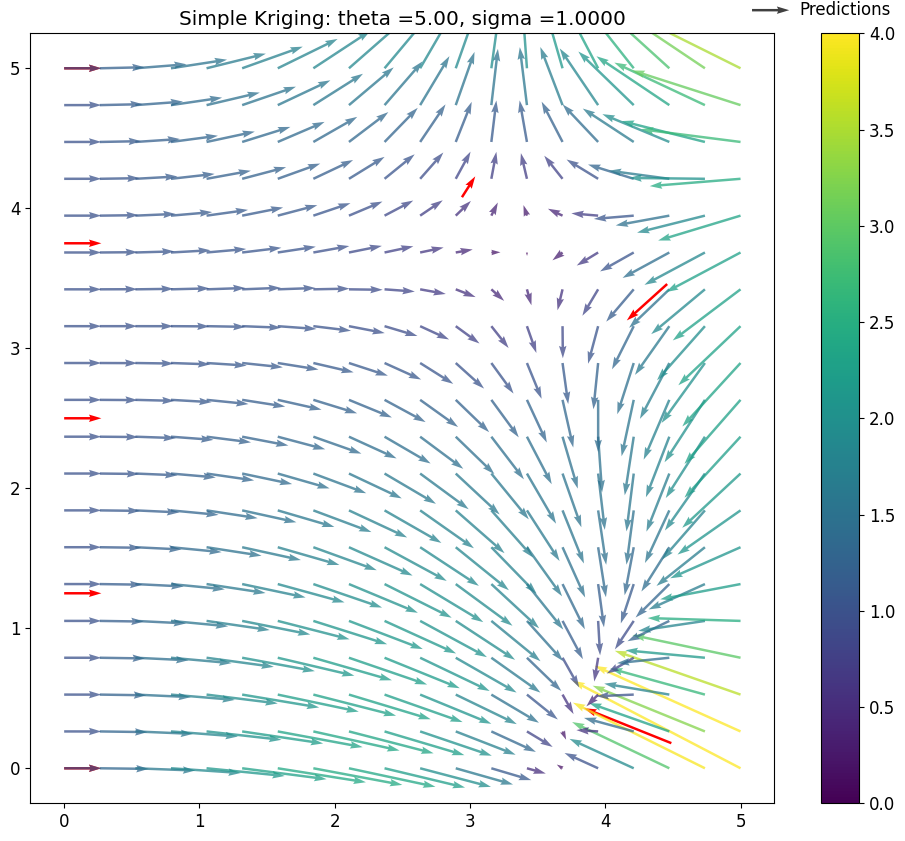} \\
(a) $\theta = 0.5$ & (b) $\theta = 2$ & (c) $\theta = 5$
\end{tabular}
\caption{Simple Kriging (using $\phi$) prediction based on random additional observations in the domain for varying lengthscale parameter (flow-reversal phenomena). Once again, the color and length of the arrows are proportional to the magnitude of velocity.}
\label{fig: thetaflowreversal}
\end{figure}

Our primary motivation is to address problems in computational fluid dynamics. As an introductory problem, we consider the steady-state flow around an obstacle that is incompressible, inviscid, with negligible external forces and we assume it satisfies a Neumann boundary condition on the obstacle boundary. These are simplifying assumptions that offers a linearized version of the Navier-Stokes equations. The potential flow formulation is frequently used to compute drag and lift coefficient of airfoils. It is described by the following equations,
\begin{eqnarray}
\label{eqn:energyconservation}
\rho (\bm{v} \cdot \nabla)\bm{v} &=& -\nabla P \\
\label{eqn:laplacianzero}
\nabla \cdot \bm{v} &=& 0,  \\
\label{eqn:Boundaries}
\bm{v} \cdot \bm{n} = 0 \quad \text{on $\Gamma$} &\quad\text{ and }\quad &\bm{v} = \bm{V}_{\infty} ,\; P = P_{\infty} \quad \text{far from the obstacle},
\end{eqnarray} 
where $\rho := \rho(\bm{x}, t)$ is the density of the fluid, $\bm{v} := \bm{v}(\bm{x}, t)$ is the velocity vector, $P := P(\bm{x}, t)$ is the pressure field and $\bm{n}$ denotes the unit vector normal to the obstacle. Further assuming that the flow is irrotational, i.e., $\nabla \times \bm{v} = 0$, allows to have a flow that derives from a potential, $\bm{v} = \nabla \phi$ where $\phi$ is a scalar field. Additionally, the pressure can be deduced from the velocity through Bernoulli's theorem, 
\begin{eqnarray*}
\frac{1}{2} \rho \Vert\bm{v}\Vert^2 + P &=& \frac{1}{2} \rho \Vert\bm{V}_{\infty} \Vert^2 + P_{\infty}
\end{eqnarray*} 
With the potential function, Equations~\eqref{eqn:laplacianzero}-\eqref{eqn:Boundaries} become
\begin{equation}
\label{eqn:PotentialPb}
\nabla^2 \phi = 0 \quad , \quad
\nabla \phi \cdot \bm{n} = 0 \quad \text{on $\Gamma$}\quad , \quad
\nabla \phi = \bm{V}_{\infty} \quad \text{at infinity}
\end{equation}
For future reference, we refer to all the simplifying fluid assumptions namely, incompressibile, invsicid and irrotational with negligible external forces as the hypothesis $\HF$. We start by modeling $\phi:\R^2\to \R$ as a random function with the squared exponential kernel as covariance function. We intend to use $n$ point observations of $\nabla \phi$ to predict the same in the rest of the domain while taking into consideration $\nabla^2 \phi = 0$ at specific collocation points in the domain and $\nabla \phi \cdot \bm{n} = 0$ at collocation points on the obstacle boundary. In the case of Lagrangian Kriging, there is no notion of collocation points and hence these differential equations will be treated as constraints at the prediction points. 
\newline\newline
\noindent 
The usual 2D velocity vector field $\bm{v}:\R^2 \to \R^2$ has two component functions, $v_x = \partial_x \phi$ and $v_y = \partial_y \phi$, the projections onto the Euclidean basis. Each velocity vector observation is a tuple $(v_x, v_y)$. The classical approach is to interpolate the flow using simple Kriging assuming there is no cross-covariance between $v_x$ and $v_y$ (independence). The potential function ($\phi$) formulation is of no use here since we assume independence. We use the same squared exponential kernel directly on the velocity components. The only exact observations are velocity vectors at infinity but this is not very useful since there is no way to account for the presence of an obstacle. Nonetheless, this is serves as a good baseline interpolation model to compare with. 
\newline\newline
\noindent As a first experiment (without an obstacle), we consider the task of predicting the flow based on random velocity vector observations. This is to demonstrate the difference between predictions of simple Kriging assuming independent velocity components (formerly described classical approach), simple Kriging using the cross-covariance relations arising from the potential function $\phi$ formulation (same as simple Kriging in the extended design space) and co-Kriging with $\nabla^2 \phi = 0$ at $p=100$ points. Figure \ref{fig: flowwithrandomobs} reveals some interesting features. We notice a peculiar flow reversal and other sophisticated flow patterns in Figure \ref{fig: flowwithrandomobs}b, all while retaining the interpolation property. On adding collocation points for co-Kriging, the continuity equation ensures that there is no back-flow but we observe that, the flow tends to enter or exit the domain erratically as shown in Figure \ref{fig: flowwithrandomobs}c. To investigate the flow patterns further, we assess the impact of changing the lengthscale parameter $\theta$ of the squared exponential kernel, see Figure \ref{fig: thetaflowreversal}. This clearly illustrates that for larger lengthscales, the peculiar flow patterns disappear and the observed velocity vectors have a stronger effect on the predictions. %

\begin{figure}[H]
\centering
\includegraphics[width=0.3\linewidth]{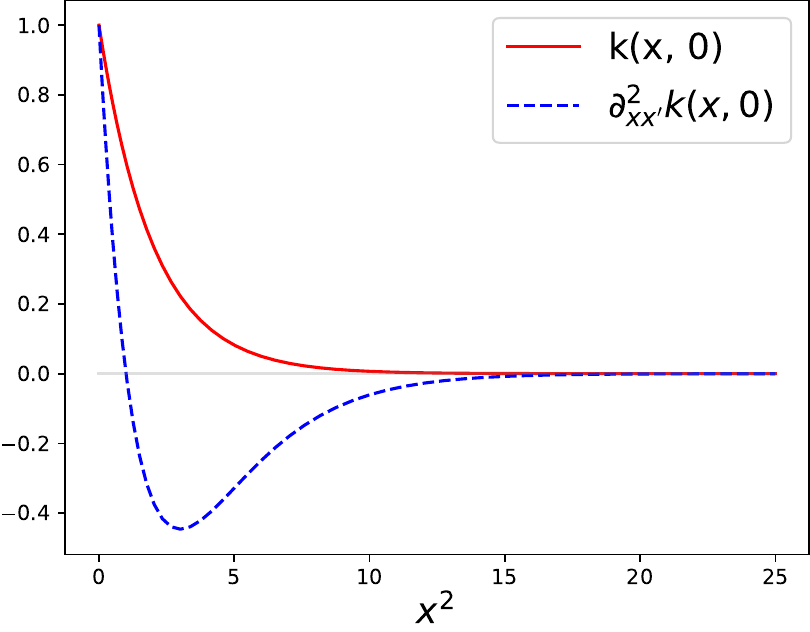}
\caption{The variation of covariance with respect to the squared distance for the squared exponential kernel ($\theta = 1, \sigma^2=1$).} 
\label{fig: covariancecurveSqExp}
\end{figure}

\noindent A possible contributor to the flow-reversal phenomena may be the negative correlations that appear when computing covariances using the differentiated kernel. This is illustrated for the squared exponential kernel in Figure \ref{fig: covariancecurveSqExp}. A significant portion of the double differentiated kernel is negative. Another possible interpretation is that the flow prediction is the consequence of appropriately placed fundamental potential flows (sources and sinks) such that their superposition is interpolating. The exact mechanism of how this can be facilitated by covariance relationships is unclear and beyond the scope of this paper and requires further research. For our experiments, we noticed that, spread-out observations in combination with the LOOCV hyperparameter tuning, counters this effect. 

\begin{figure}[H]
\centering
\begin{tabular}{ccc}
\includegraphics[trim = {0 0 0 0.5cm}, clip, width=0.3\linewidth]{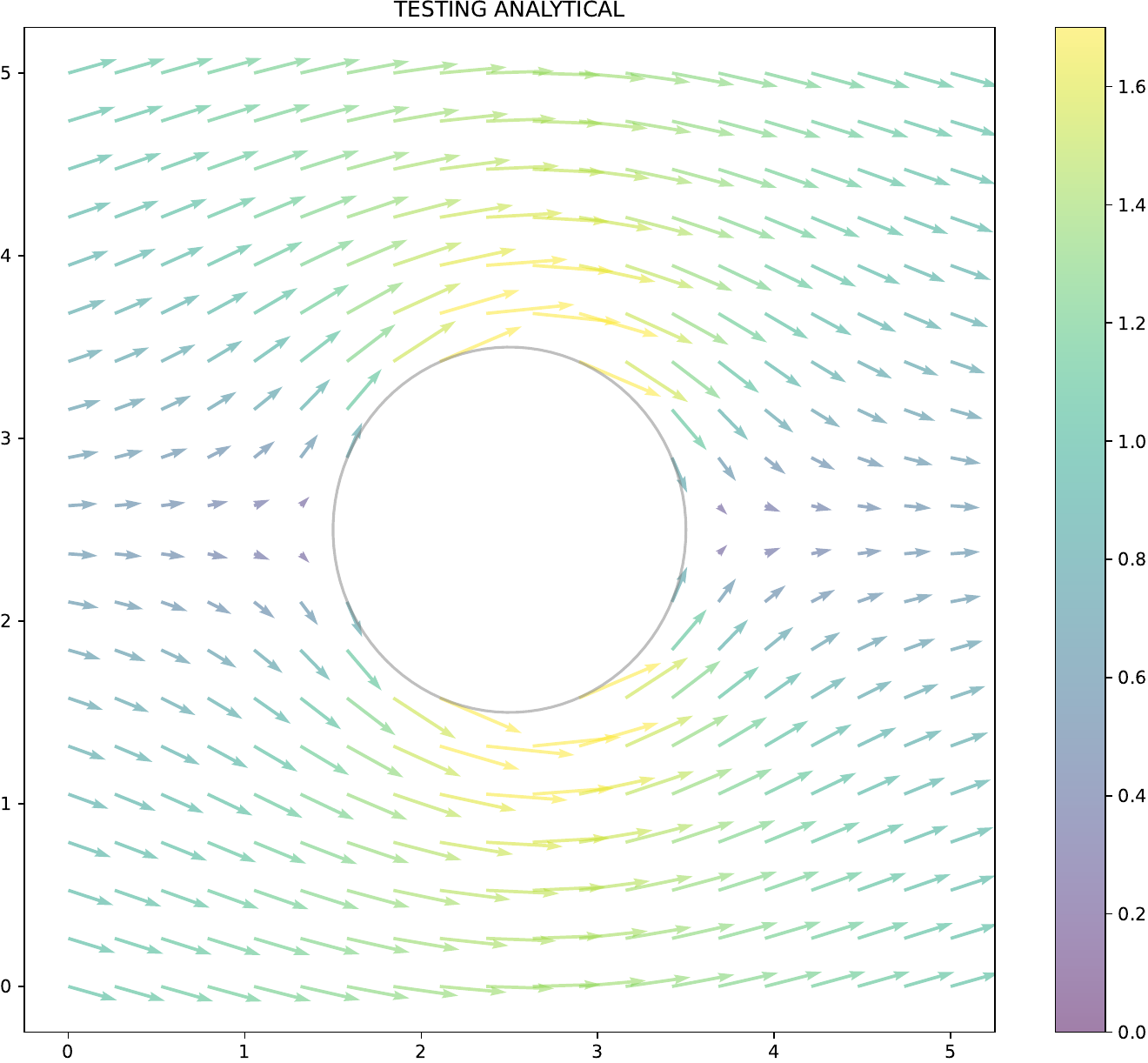} & 
\includegraphics[trim = {0 0 0 0.5cm}, clip, width=0.3\linewidth]{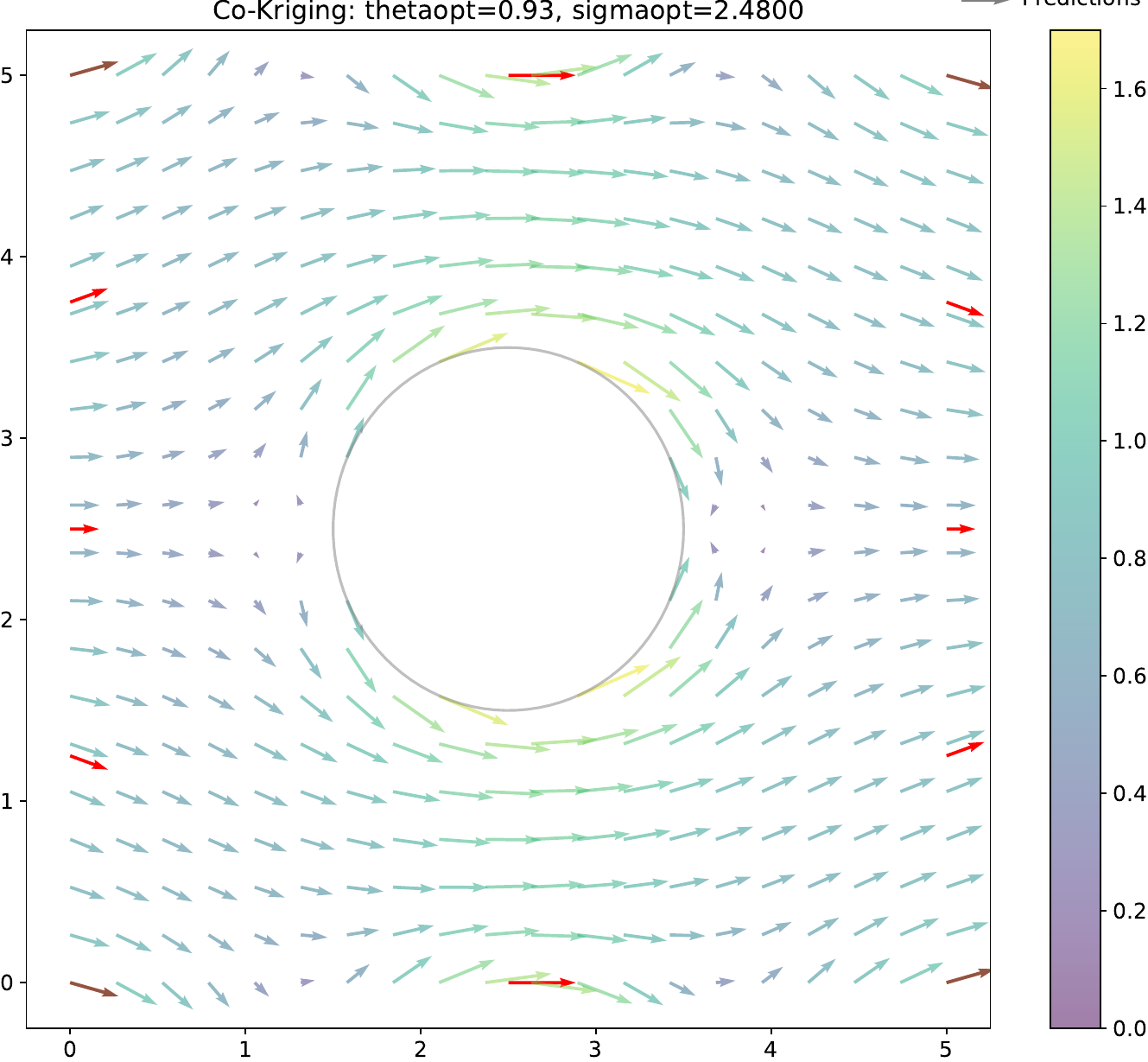} &
\includegraphics[trim = {0 0 0 0.5cm}, clip, width=0.3\linewidth]{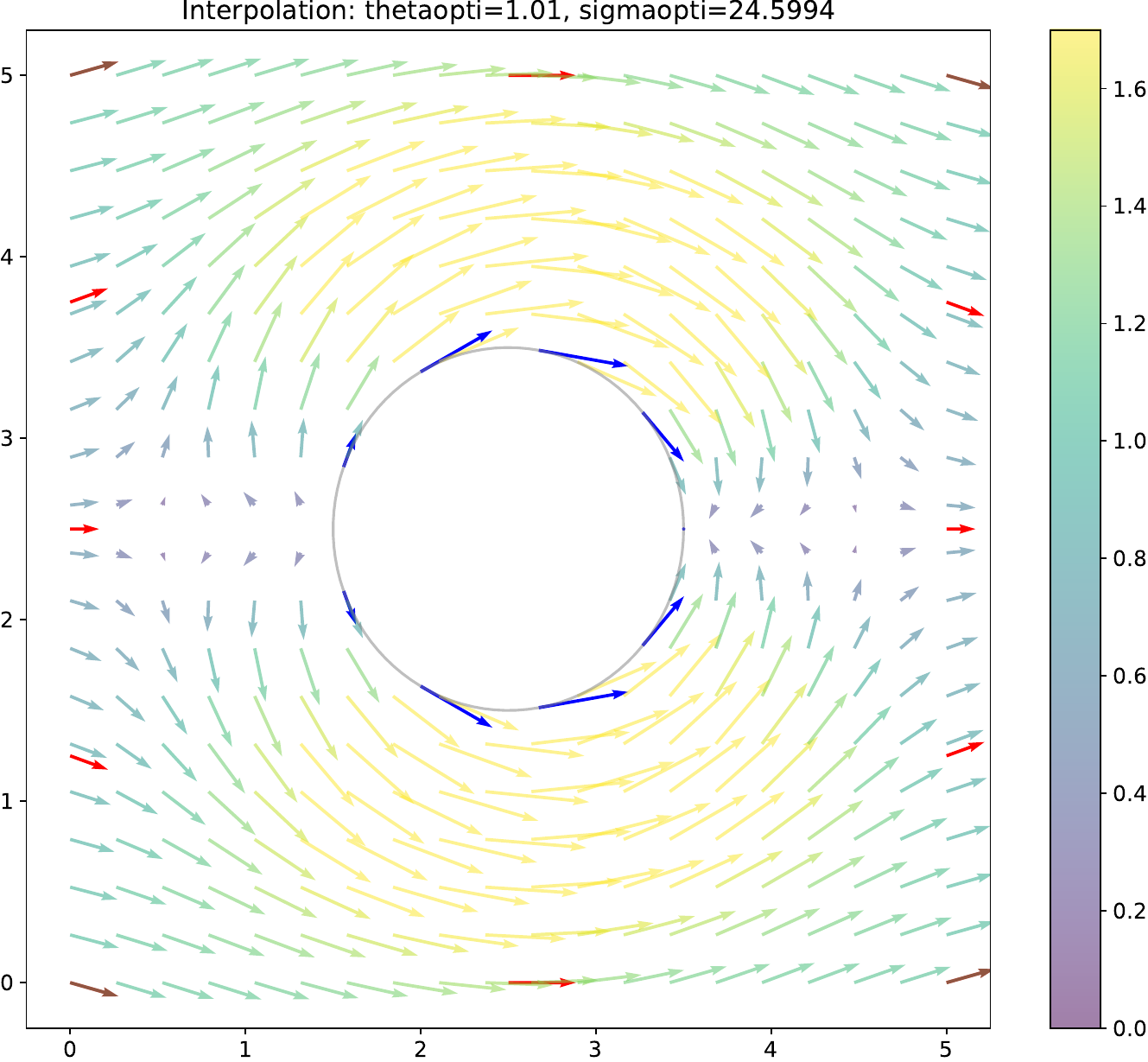} \\
(a) True flow under $\HF$ & (b) Co-Kriging & (c) Lagrangian Kriging 
\end{tabular}
\caption{Comparing the predicted flows. Observed velocity vectors in red. The LOOCV optimal co-Kriging parameters, $\hat{\theta} = 0.93$ and $\hat{\sigma} = 2.48$. The Lagrangian Kriging predictions (optimal lengthscale parameter, $\hat{\theta} = 2.01$) are shown in blue, using which, the rest is a simple Kriging based interpolation.}
\label{fig: finalfluidflows}
\end{figure}

\noindent For the obstacle geometry, we consider the cylindrical shape and since we have an analytical solution for a perfect flow around cylinder, we will use it to compare the two approaches. We use $n=12$ velocity vector observations, all observations are based on the analytical solution and kept relatively far from the obstacle. Using $q_1 = 10$ collocation points around the boundary of the obstacle to enforce the Neumann boundary condition and $q_2 = 88$ collocation points in the domain to ensure continuity. the collocation points translate to the prediction points in the case of Lagrangian Kriging, as seen before. 
\newline\newline
\noindent Specifically for Lagrangian Kriging, we cannot meaningfully consider the continuity equation does not directly involve the velocity vectors and we find ourselves in the situation described in Remark \ref{remark: LK} and Remark \ref{remark: LK continued}. which leaves us with only one constraint involving the velocity vectors, i.e., velocity vectors must be tangential around the obstacle boundary. This restricts the velocity predictions to the obstacle boundary only. For these predictions, we compute the optimal lengthscale $\hat{\theta}$ as the minimizer of the virtual LOOCV criterion \eqref{eqn: virtualLOOCVMSE} since the observation locations are unconstrained. For the Figure \ref{fig: finalfluidflows}c, we adopt a two-step strategy. Firstly, we make Lagrangian Kriging predictions along the $q_1$ boundary collocation/prediction points (shown in blue) and secondly, we interpolate the rest of the flow assuming, the original observations and the Lagrangian Kriging boundary predictions as exact observations, using simple Kriging with no cross-covariance among the two velocity components $v_x$ and $v_y$ as demonstrated in Figure \ref{fig: flowwithrandomobs}a. We tune the lengthscale parameter separately for the second interpolation step with $\hat{\theta} = 1.01$ minimizing the LOOCV MSE. Additionally, if useful, we can make direct Lagrangian Kriging predictions for $\partial^2_{xx} \phi$  and $\partial^2_{yy} \phi$ since they are part of the $q_2$ collocation/prediction points. This was of no immediate use to us and hence we omitted any figures about it. 

\begin{figure}[H]
\centering
\begin{tabular}{cc}
\includegraphics[trim = {0 0 0 0.45cm}, clip, width=0.3\linewidth]{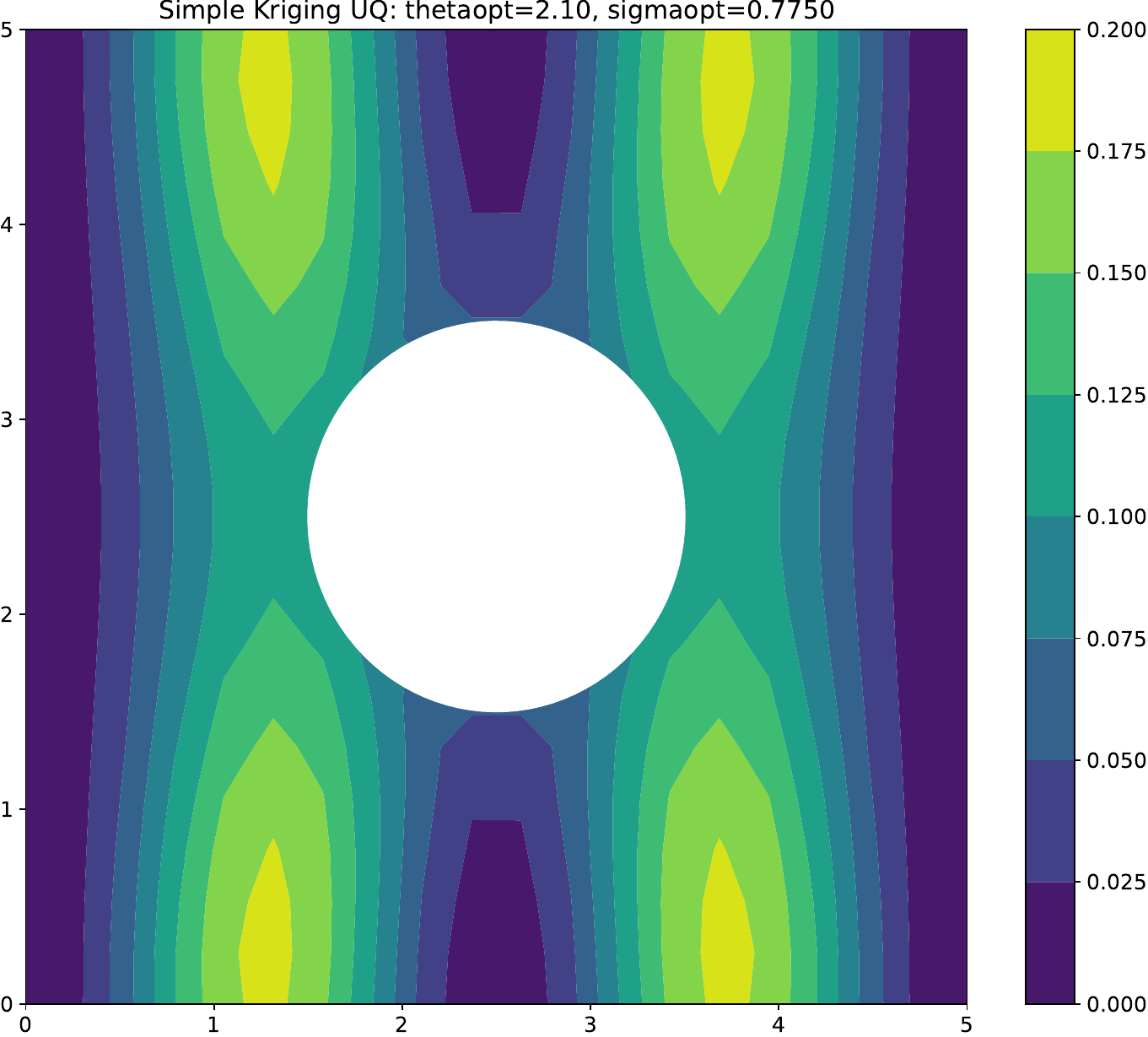} &
\includegraphics[trim = {0 0 0 0}, clip, width=0.3\linewidth]{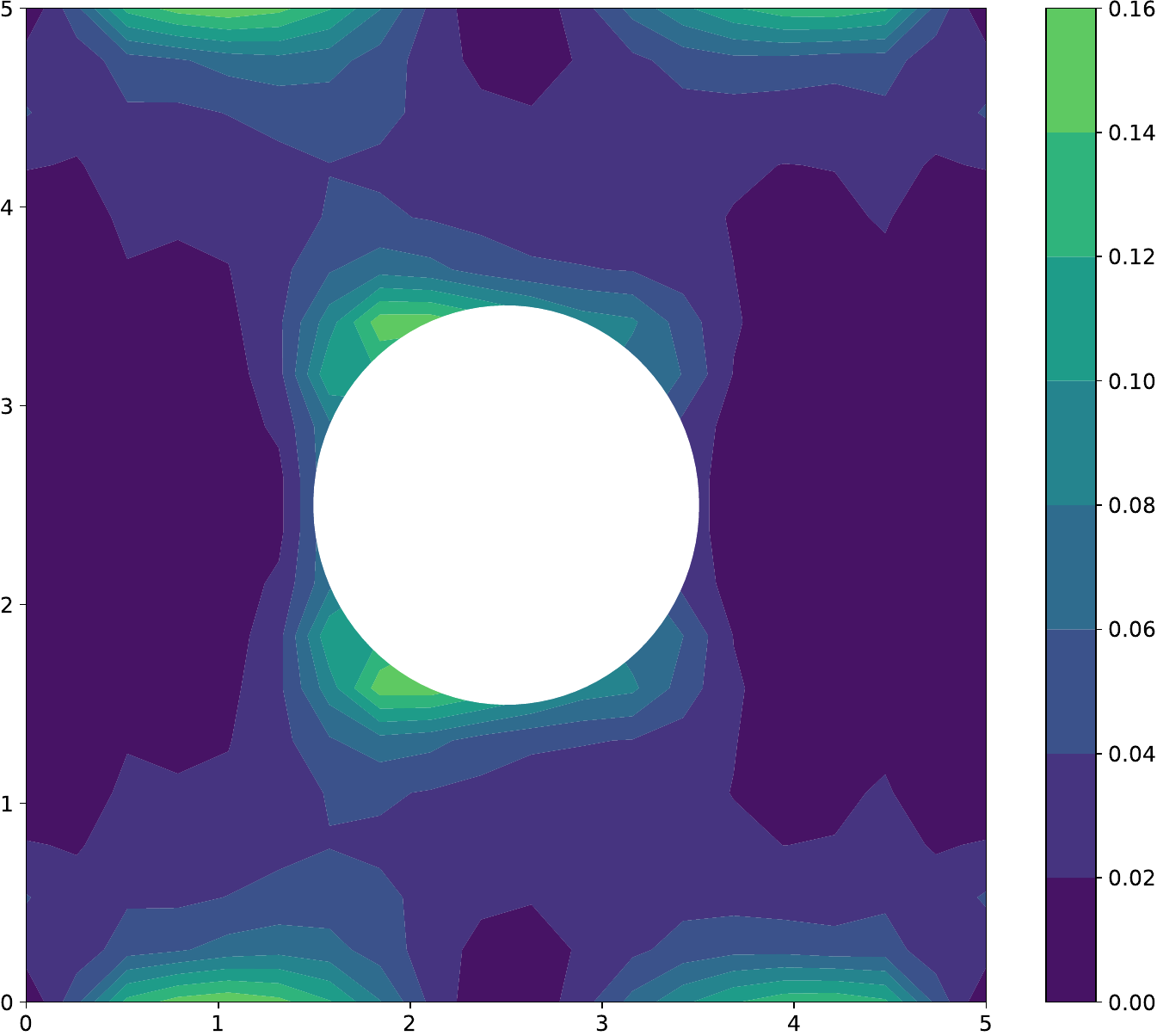} \\
(a) & (b) \\
\end{tabular}
\caption{Uncertainty in the squared magnitude of the predicted velocity vectors: (a) Simple Kriging with LOOCV optimal parameters $\hat{\theta} = 2.10$ and $\hat{\sigma}=0.77$, (b) Co-Kriging with LOOCV optimal parameters $\hat{\theta} = 0.93$ and $\hat{\sigma} = 2.48$.}
\label{fig: uncertaintyfluid}
\end{figure}

\begin{figure}[H]
\begin{tabular}{ccc}
\includegraphics[trim = {0 0 0 0.5cm}, clip, width=0.3\linewidth]{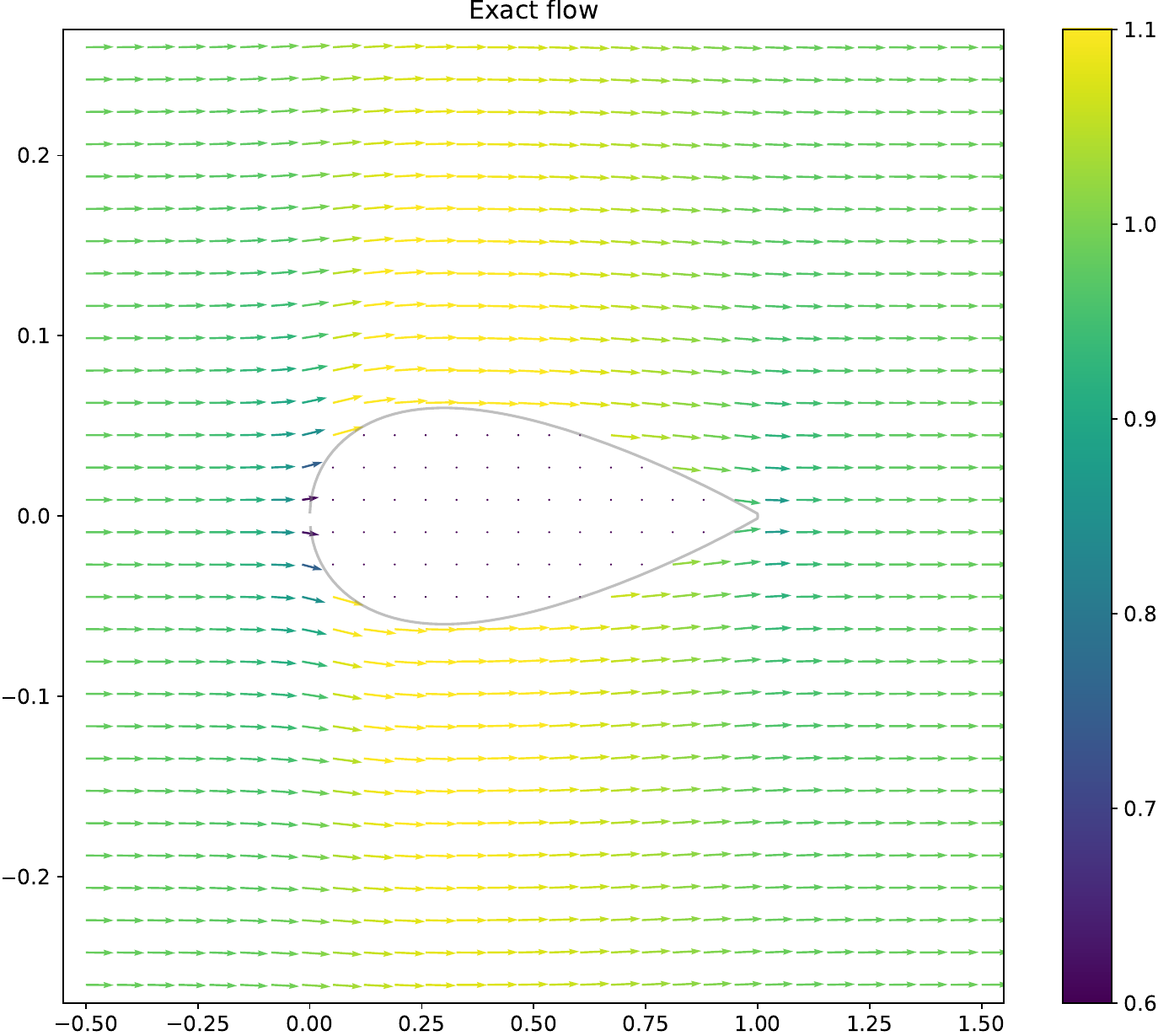} & 
\includegraphics[trim = {0 0 0 0.5cm}, clip,width=0.3\linewidth]{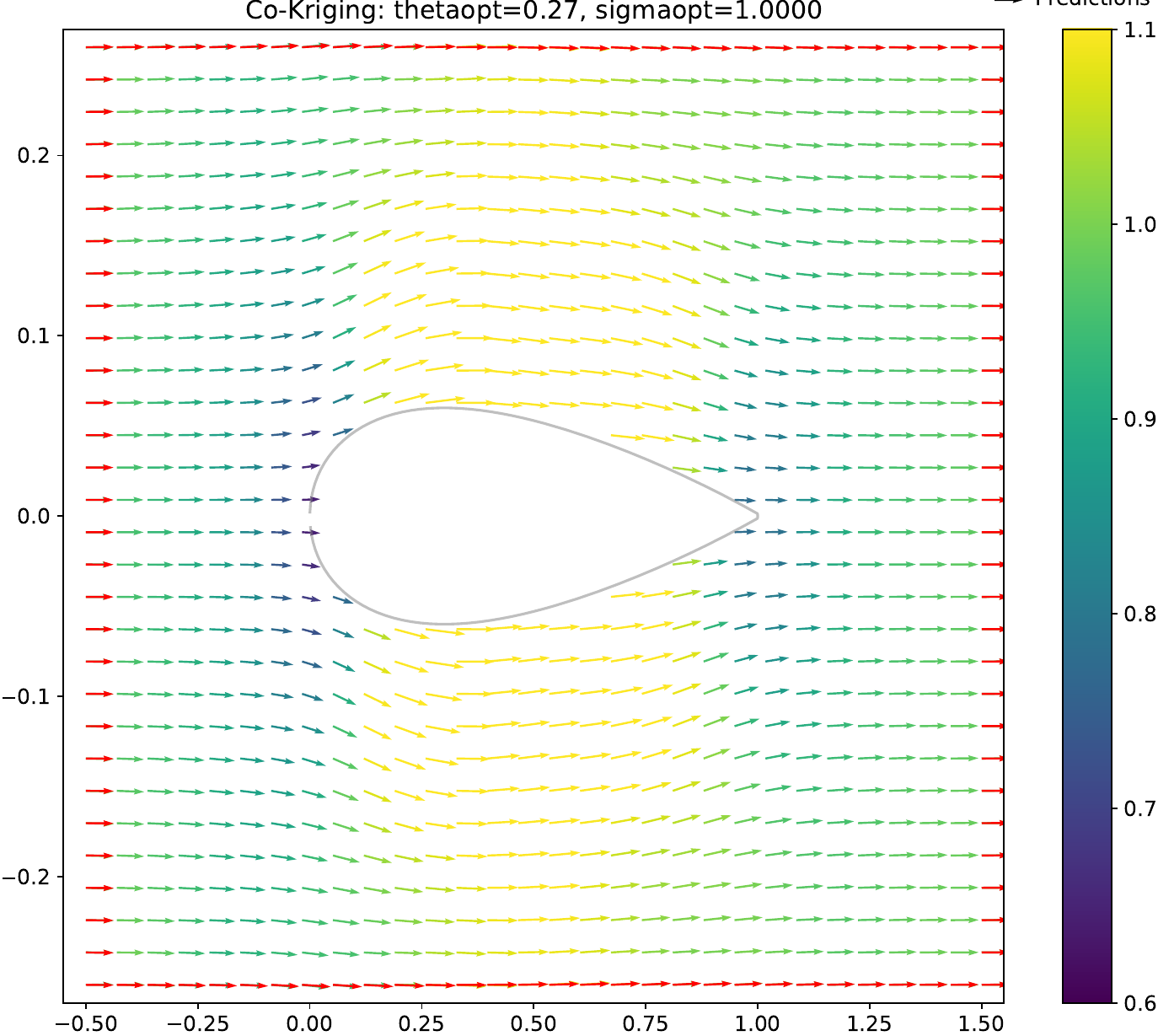} &
\includegraphics[trim = {0 0 0 0.5cm}, clip,width=0.3\linewidth]{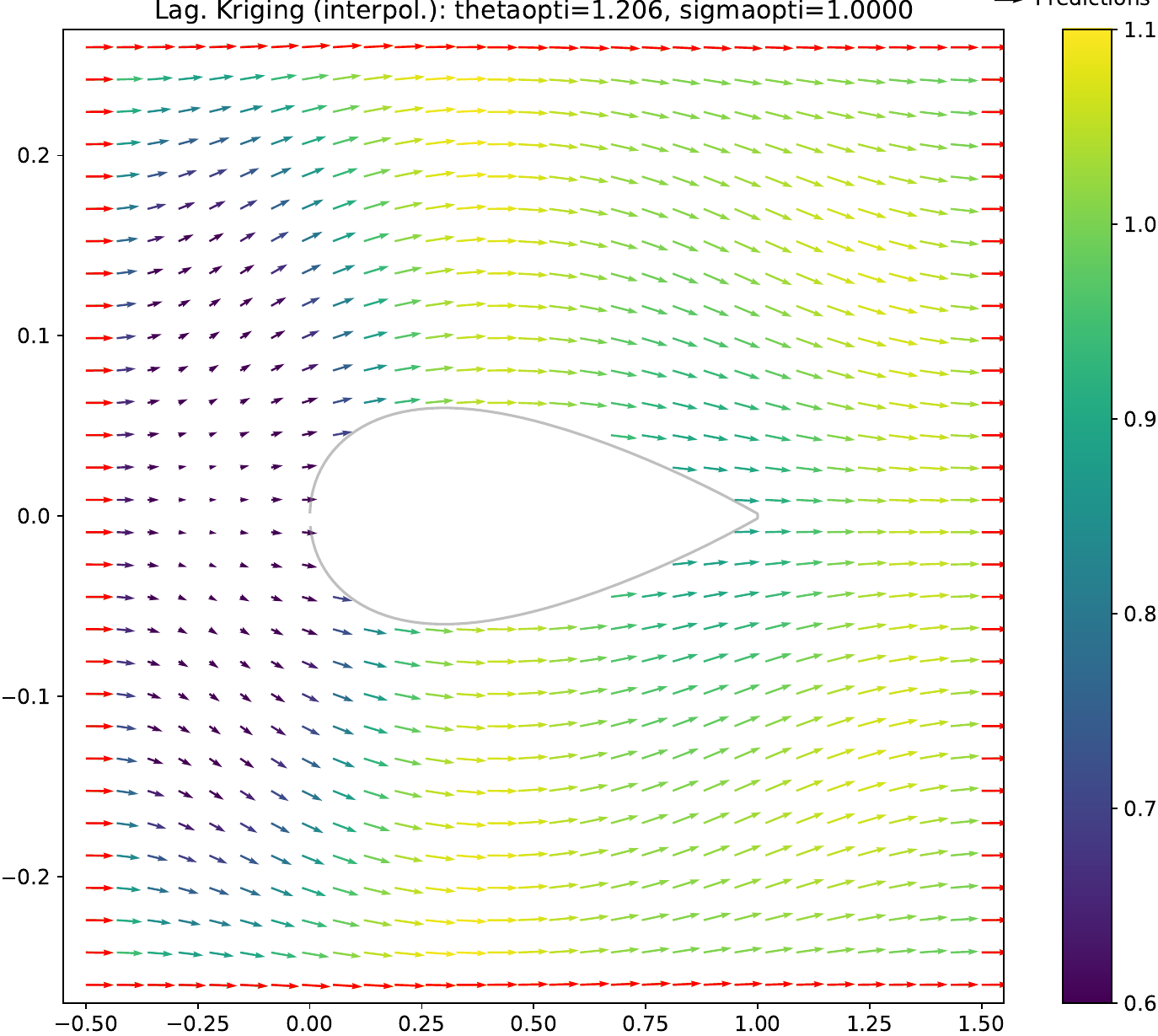} \\
& $\hat{\theta} = 0.27$ & $\hat{\theta} = 0.25$ \\
\includegraphics[trim = {0 0 0 0.5cm}, clip,width=0.3\linewidth]{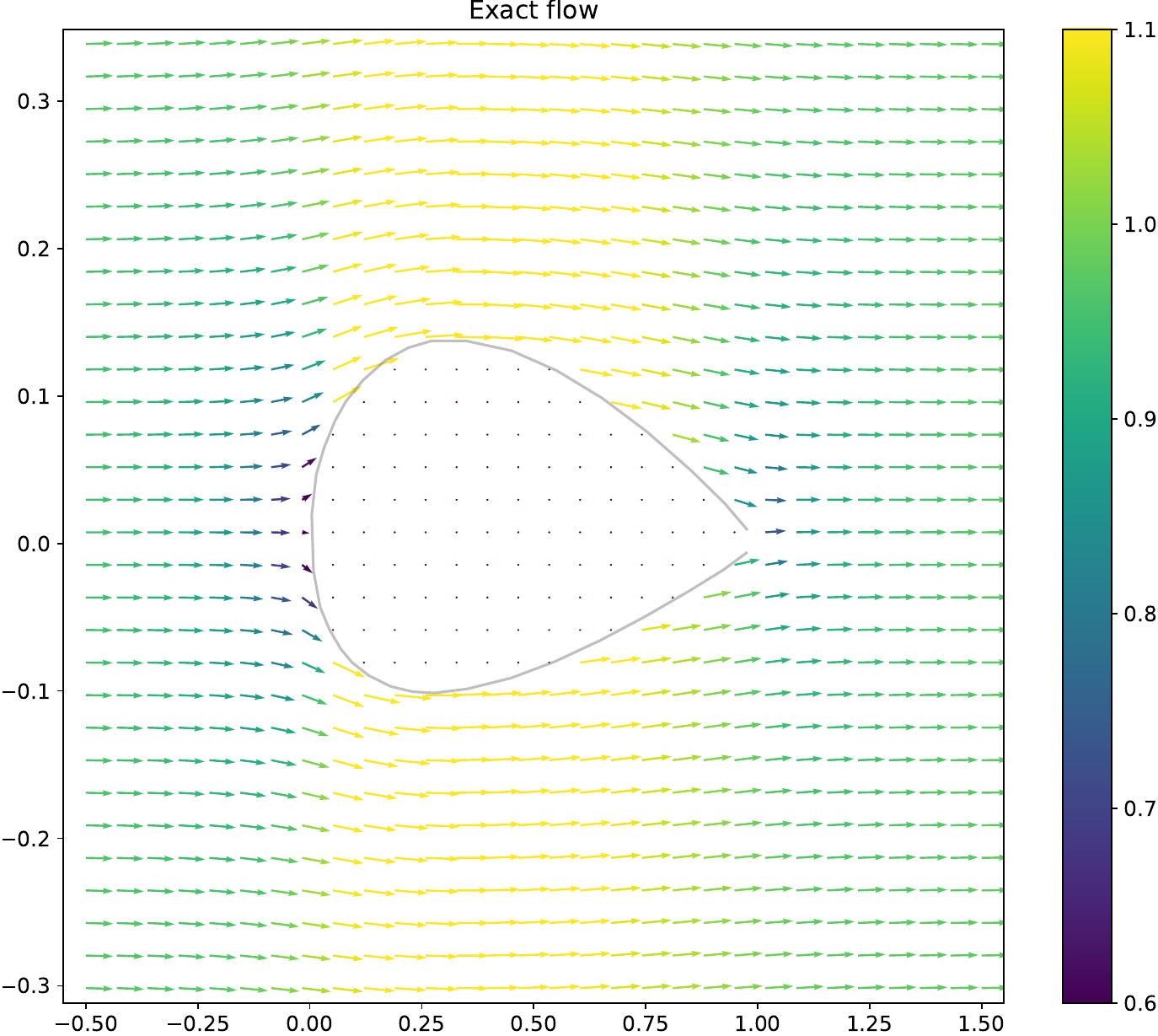} & 
\includegraphics[trim = {0 0 0 0.5cm}, clip,width=0.3\linewidth]{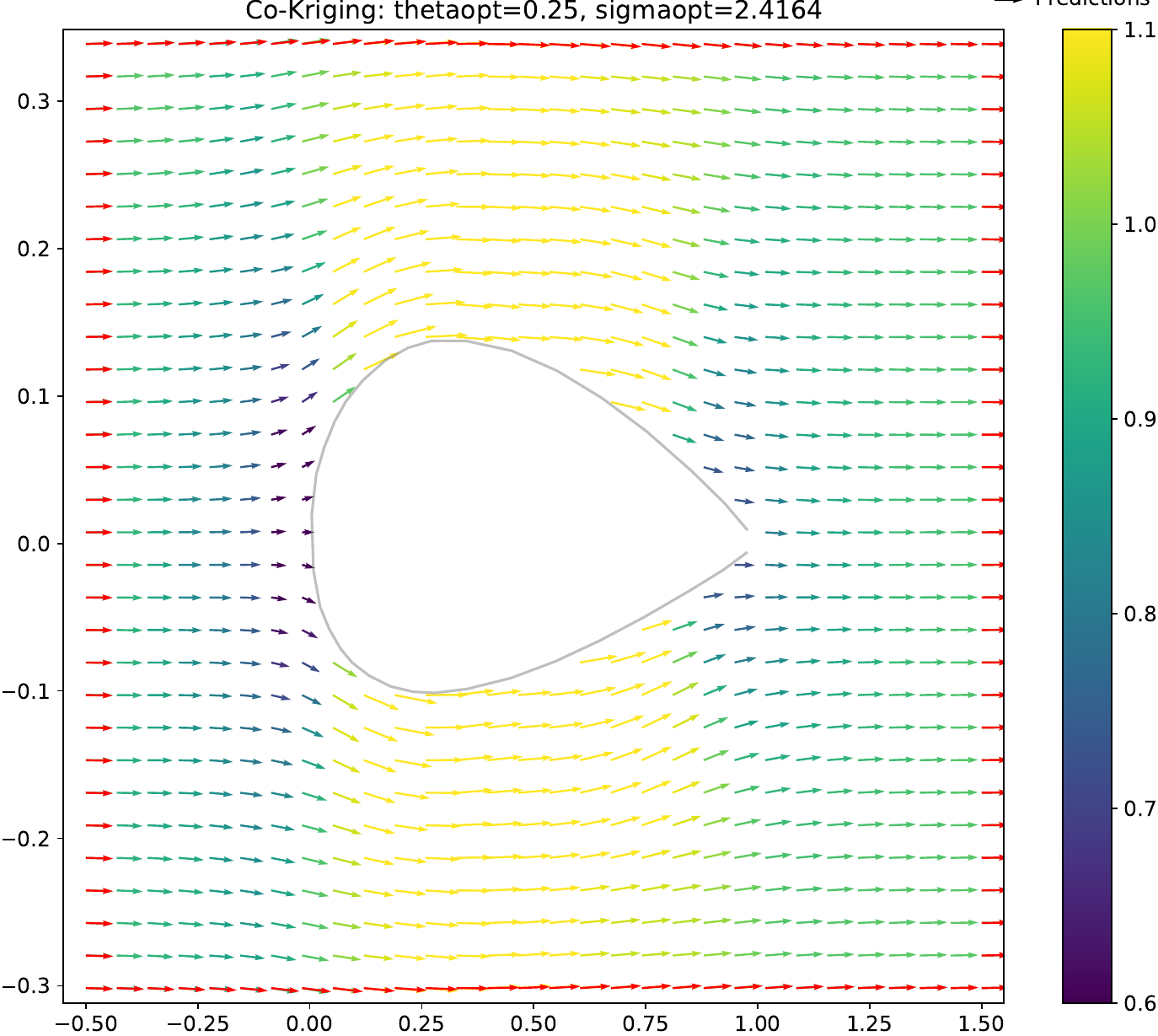} &
\includegraphics[trim = {0 0 0 0.5cm}, clip,width=0.3\linewidth]{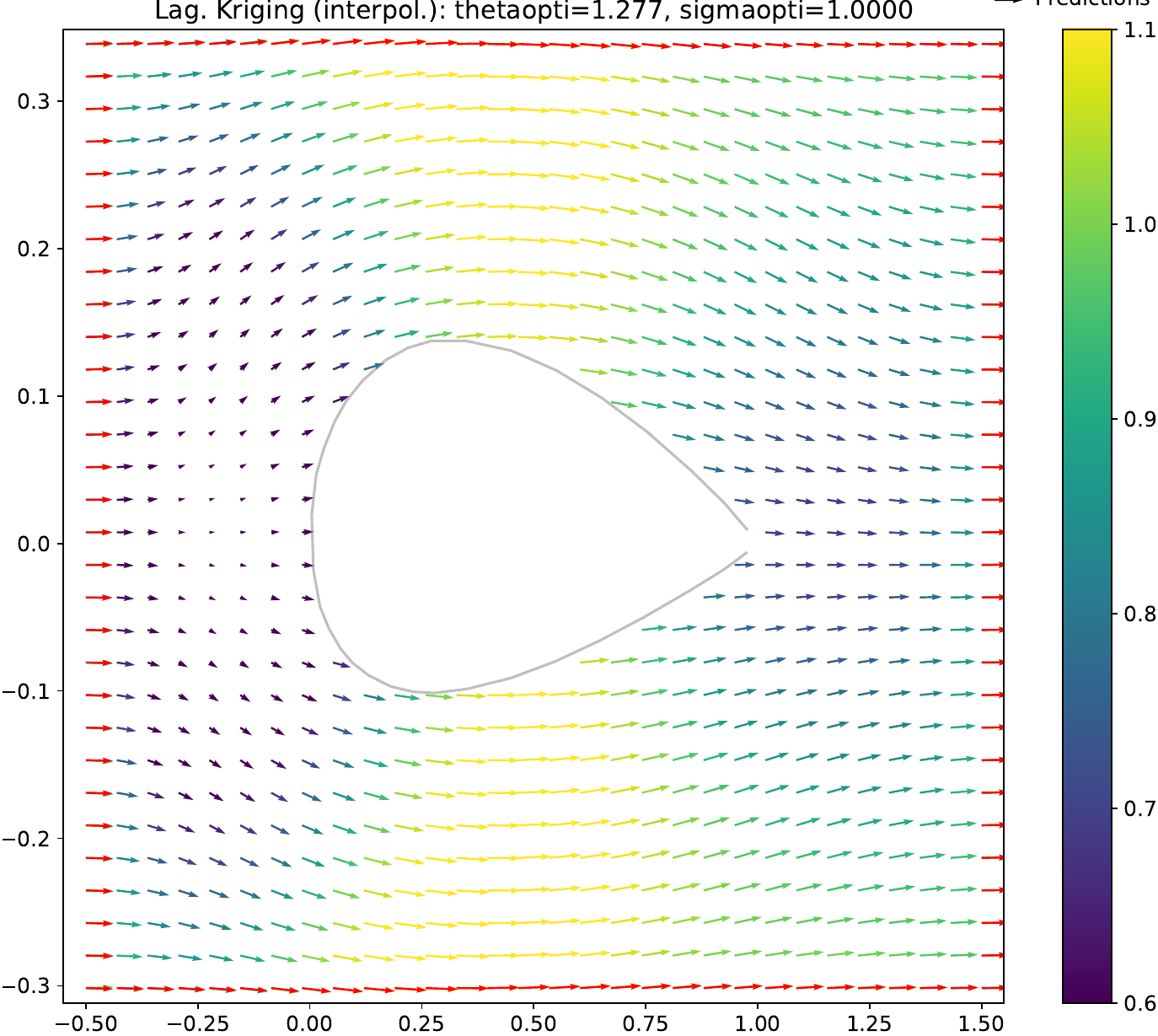} \\
& $\hat{\theta} = 0.25$ & $\hat{\theta} = 0.26$ \\ 
(a) Exact flow under $\HF$ & (b) Co-Kriging & (c) Lagrangian Kriging
\end{tabular}
\caption{The exact flow and co-Kriging and Lagrangian Kriging predictions for the airfoils NACA 0012 (top) and NACA 2424 (bottom). Observed velocity vectors marked in red.}
\label{fig: airfoil predictions}
\end{figure}

\noindent The flow prediction in Figure \ref{fig: finalfluidflows}b, tends to enter and exit the domain in the lack of velocity observations, which was also the case with Figure \ref{fig: flowwithrandomobs}c. On the bright side, the co-Kriging model uncertainty is consistent with this observation. Uncertainty spikes in the vicinity of the prediction locations where the flow enters and exits the domain as shown in Figure \ref{fig: uncertaintyfluid}b. %
To further explain the uncertainty Figure \ref{fig: uncertaintyfluid}, assume the model predictions of $\hat{v}_x$ and $\hat{v}_y$ as two dependent normal distributions. From Section \ref{sec: UQ}, recall that the squared sum of two dependent normal distributions is a generalized $\chi^2$ distribution. Finally, we computed the variance (and standard deviation) of this distribution to represent uncertainty in $\Vert \bm{v} \Vert^2_2 = v_x^2 + v_y^2$ as the $2\sigma$ value. We perform this procedure for both Simple Kriging in the extended design space and Co-kriging with differential information. From Figure \ref{fig: uncertaintyfluid}, it is clear that co-Kriging substantially reduces the model uncertainty across the domain compared to simple Kriging which solely accounts for the observations of velocity. Moreover, the co-Kriging uncertainty aligns well with the model performance in Figure \ref{fig: finalfluidflows}b, supporting its reliability in practice.\\[0.5em]%
\noindent We extend the same methodology to predict the potential flow around NACA airfoil geometries, particularly of interest to the aerospace engineering community and the aviation industry. %
We use the open source \texttt{UIUC Airfoil database} - \href{https://m-selig.ae.illinois.edu/ads/coord_database.html}{\texttt{\textcolor{blue}{https://m-selig.ae.illinois.edu/ads/coord\_database.html}}} to retrieve the airfoil coordinates and solve for the potential flow problem using the Boundary elements method (source panel's method, to be precise) \cite{Panelsmethod} which is frequently used for airfoils. We modified the python script for the source panel's method from \href{https://github.com/jte0419/Panel_Methods}{\texttt{\textcolor{blue}{https://github.com/jte0419/Panel\_Methods}}} to apply it to custom geometries. We test our methods on two NACA airfoils: NACA 0012 and NACA 2424. We added a dense set of observations at infinity (at both $+\infty$ and $-\infty$) alongside plenty of mid-domain observations, the velocity vectors marked in red in Figure \ref{fig: airfoil predictions}, to avoid the erratic domain exit-entry behaviour as discussed before. Total number of exact observations used $n = 80$, number of continuity collocation points (for co-Kriging) $p_1 \leq 100$, number of boundary collocation points (prediction points in the case of Lagrangian Kriging) $50 \leq p_2 \leq 30$ and the number of prediction points $q \approx 900$. The continuity collocation points for a uniform grid keeping in mind the scale difference between the $x-y$ axes. The discretization on the $x$ axis is four times the discretization on the $y$ axis to respect the scale difference. The boundary collocation points are equispaced points along the obstacle boundary. We also report the LOOCV optimal $\hat{\theta}$ parameters used in each case. In the case of Lagrangian Kriging, there is a second LOOCV optimization routine for the interpolation step with $\hat{\theta} = 1.20$ and $\hat{\theta} = 1.28$ for NACA 0012 and NACA 2424 respectively.       

\begin{figure}[H]
\begin{center}
\begin{tabular}{cc}
\includegraphics[trim = {0 0 0 0.5cm}, clip,width = 0.25\linewidth]{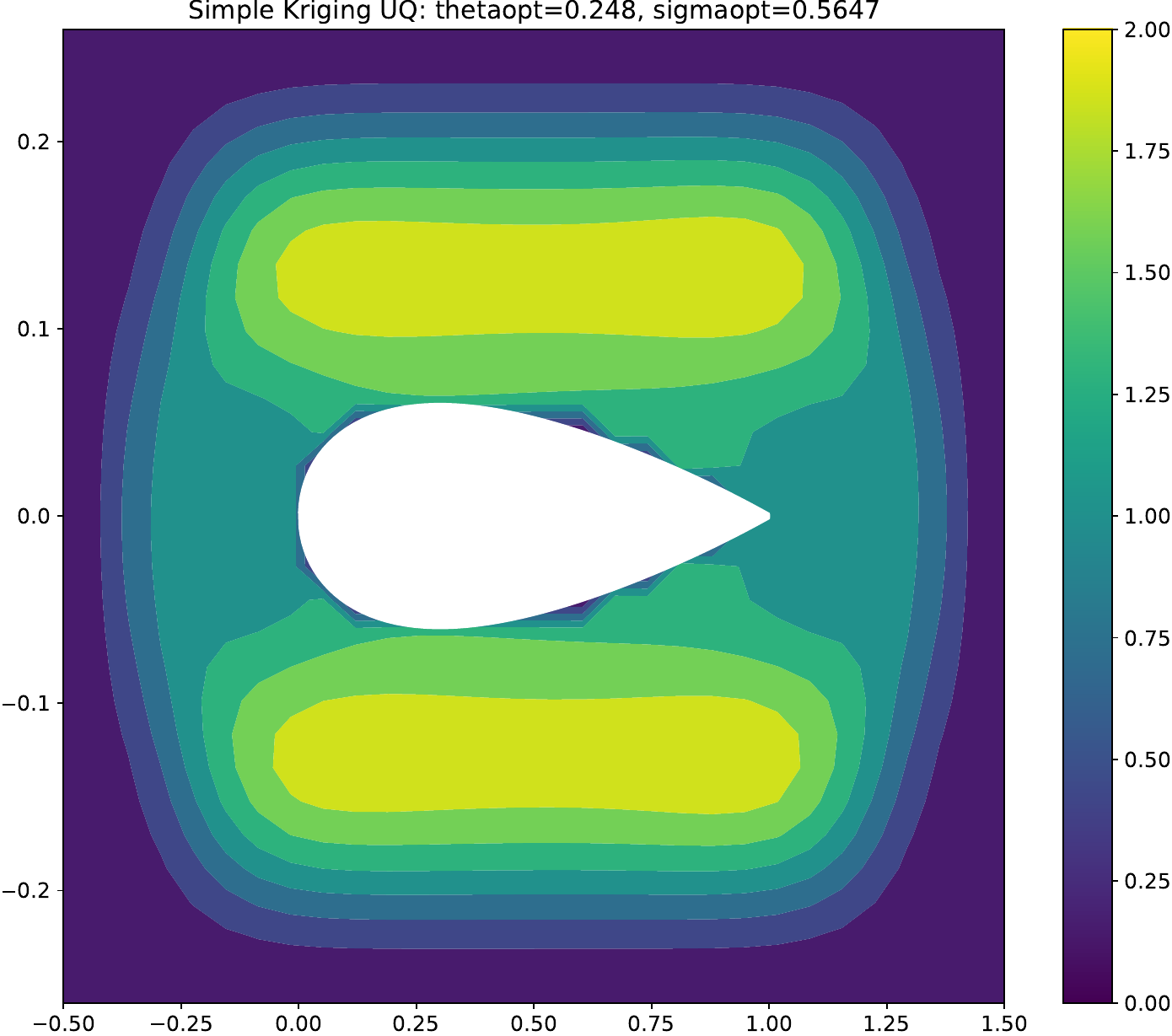} & 
\includegraphics[trim = {0 0 0 0.5cm}, clip,width = 0.25\linewidth]{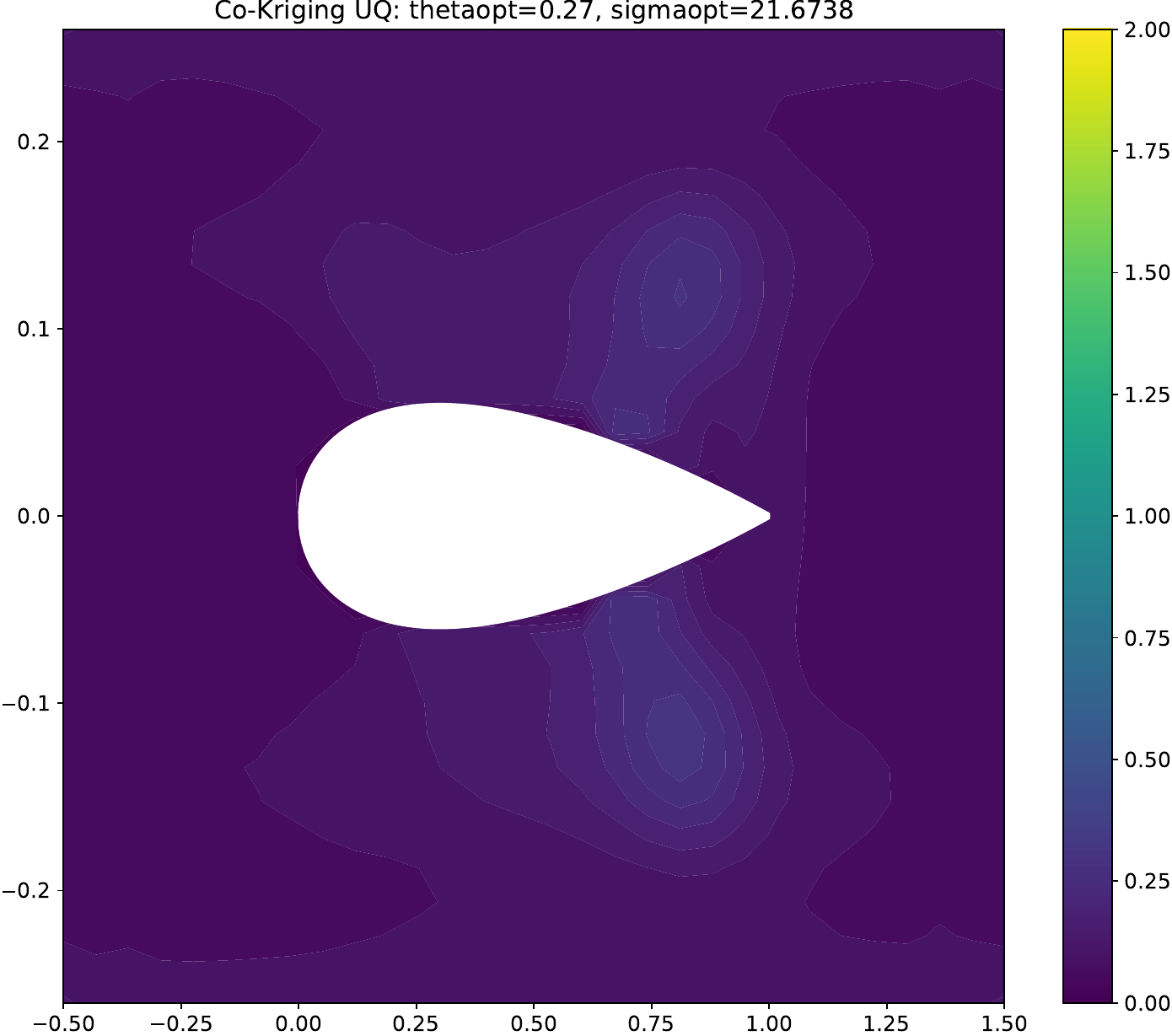} \\
$\hat{\sigma} = 0.56$ & $\hat{\sigma} = 21.67$ \\ 
\includegraphics[trim = {0 0 0 0.5cm}, clip,width = 0.25\linewidth]{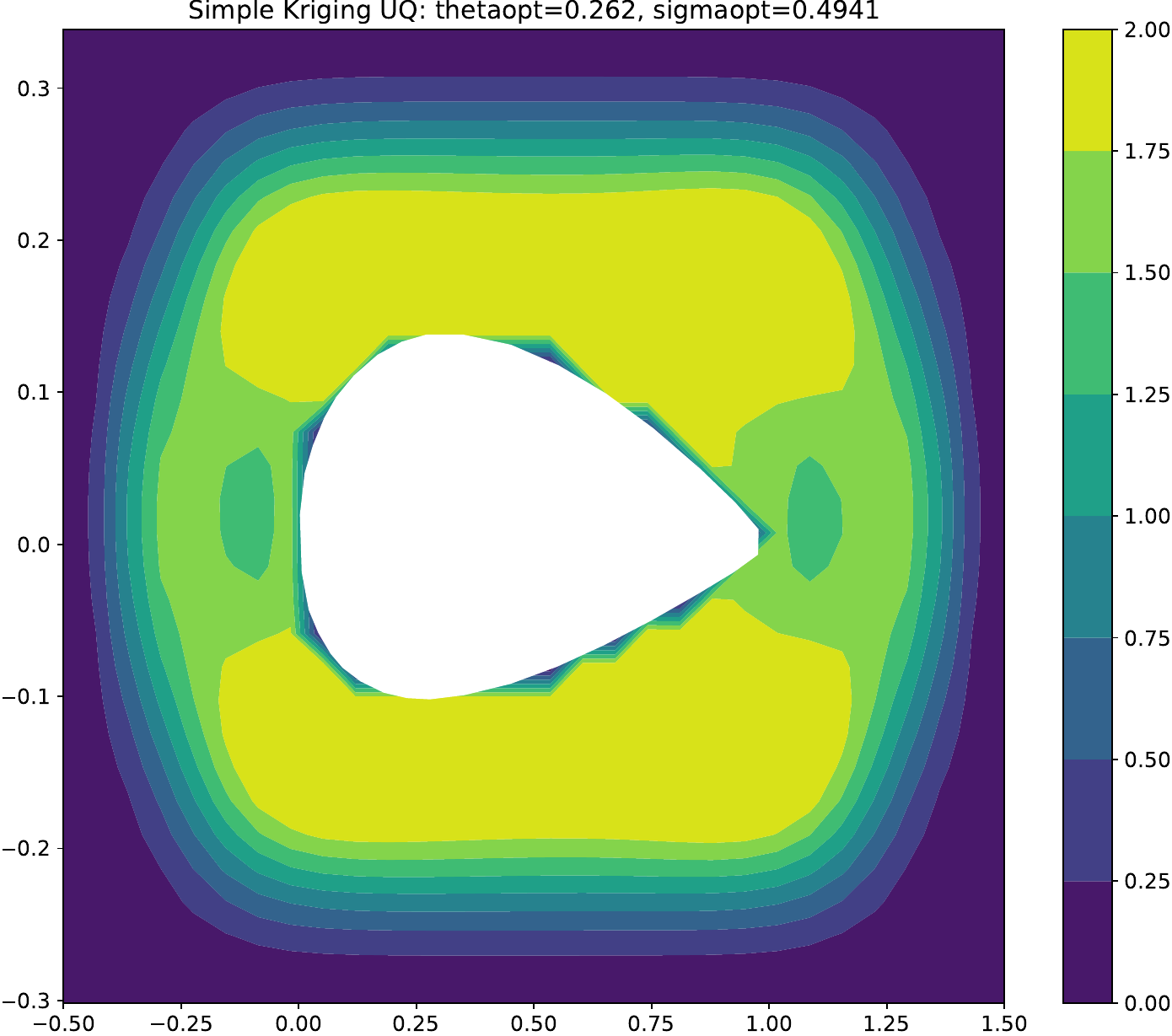} & 
\includegraphics[trim = {0 0 0 0.5cm}, clip,width = 0.25\linewidth]{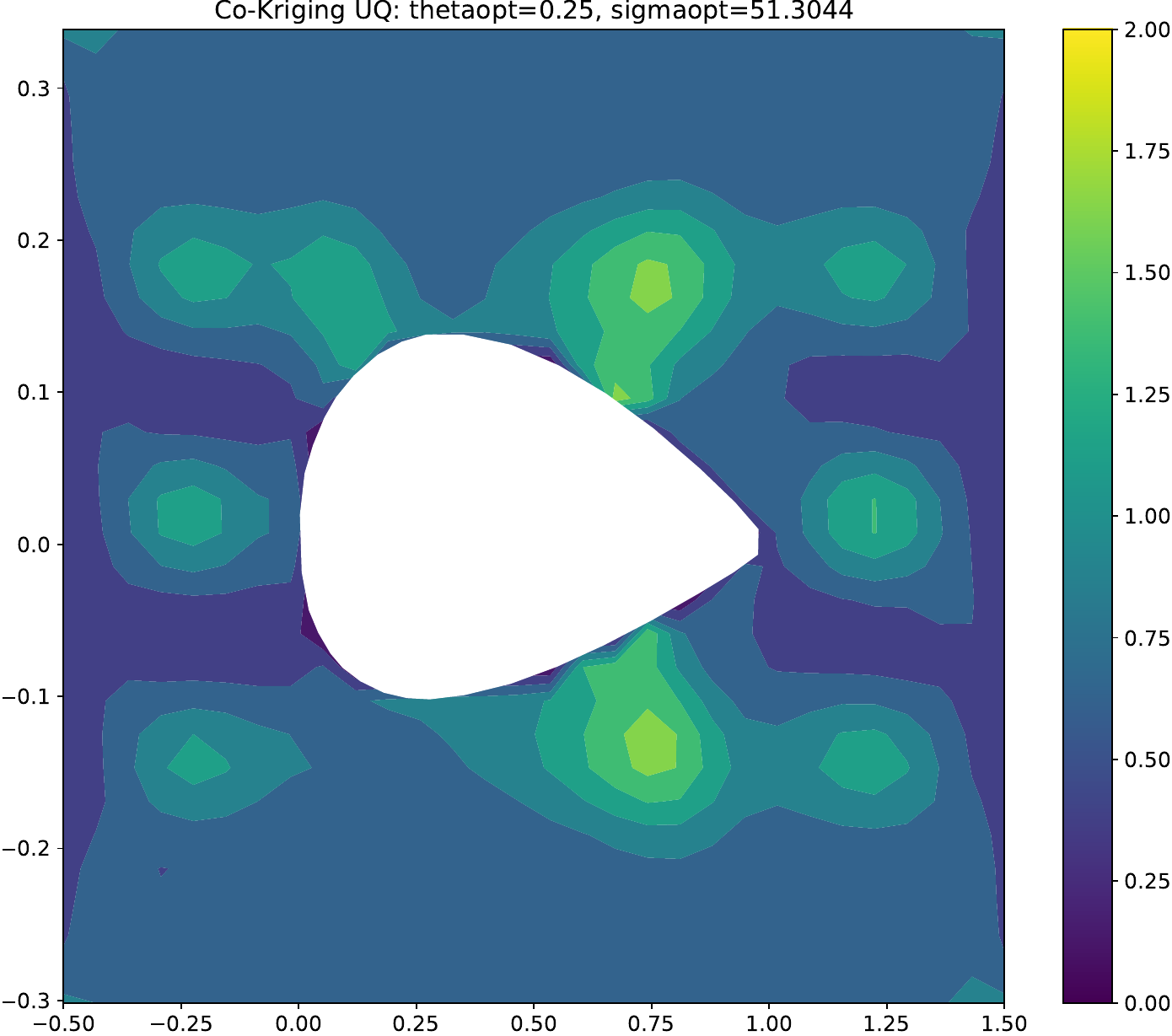} \\
$\hat{\sigma} = 0.49$ & $\hat{\sigma} = 51.30$ \\  
(a) Simple Kriging UQ & (b) Co-Kriging UQ
\end{tabular}
\end{center}
\caption{Comparing the uncertainty quantification ($\pm 2\sigma$) with simple Kriging vs. Co-Kriging for the airfoils NACA 0012 (top) and NACA 2424 (bottom) respectively.}
\label{fig: airfoil UQ}
\end{figure}

\noindent The potential flow reconstructions in Figure \ref{fig: airfoil predictions} are quite satisfactory considering that we have used an isotropic kernel. For both the airfoils, NACA 0012 and NACA 2424, the flows are tangential close to the boundary which implies that the models are well-informed of the boundary condition via collocation points. This demonstrates that our model can be generalized to non-trivial domains quite easily but it may need more careful, problem-specific kernel design, at times. Co-Kriging once again performs better than Lagrangian Kriging but this comes at a much higher computational cost. We also observed that the nugget values used
when inverting the $\matK^+$ matrix had a big influence in the final flow predictions. For co-Kriging, we used a high nugget value of $1e-4$ to avoid an ill-conditioned covariance matrix. This issue has been reported before by \cite{chen_solving_2021,Jidling2017}. For Lagrangian Kriging, we used nugget values varying from $1e-7$ to $1e-8$ for the interpolation step. We also computed the uncertainty in the squared magnitude of the velocity vectors as before and compare simple Kriging and co-Kriging in Figure \ref{fig: airfoil UQ}.     
 
\noindent From Figure \ref{fig: airfoil UQ}, it is clear that the uncertainty prediction with co-Kriging is more informative than simple Kriging for both airfoils. Once again, we had to use nugget values as high as $1e-6$ while inverting the $\matK^+$ matrix for reasonable $\hat{\sigma}$ values. The UQ, aligns well with the prediction quality. The co-Kriging model flags some locations above and below the airfoils as high uncertainty regions in all cases. The fluid flow in these regions dictates the lift force which is crucial in aerodynamic applications. Surprisingly however, the model is relatively more confident at the leading and the trailing edge of the airfoils where the prediction errors are higher as seen in Figure \ref{fig: airfoil predictions}. \\[0.5em]

\noindent To summarize, we have shown promising results with co-Kriging and Lagrangian Kriging for the potential flow prediction problem given minimal observations. However, we also faced some challenges namely, the flow-reversal and peculiar flow patterns which we linked with the lengthscale parameter. We also encountered the Lagrangian Kriging issue as described in Remark \ref{remark: LK} and Remark \ref{remark: LK continued} which we addressed by making boundary predictions at first, followed by simple Kriging (without cross-covariance). Overall, these models offer a flexible means of using boundary conditions (Dirichlet or Neumann) and any linear combination of mixed derivatives as additional information in the domain for physics-informed flow reconstruction in a linearized fluid dynamics settings with uncertainty quantification which, as demonstrated, corresponds to a measure of the prediction error.

\section{Conclusion}
\label{sec: conclusion}
The aim of the present research was to study extensions of the classical Kriging framework to a setting involving derivatives. In this context, we develop an extended design space in Section \ref{sec: extds} to handle the spatial coordinates and derivative order as a single mathematical object. Simple and ordinary Kriging in the extended design space are already flexible tools to take into account point-wise derivative observations or predict unobserved derivatives of the function. Furthermore, we observed that we can characterize linear differential equations (discretized) as matrix equations in the extended design space. We exploit this characterization to develop two Kriging methods capable of operating with general linear differential equations. The first method is a collocated co-Kriging methodology to account for linear differential information at collocation points as secondary variables in the co-Kriging predictor. %
In the second method, we exploit the optimization framework offered by Kriging, to inject differential information as constraints at points of prediction which are later resolved via Lagrange multipliers. The resultant predictor is referred to as Lagrangian Kriging. For hyperparameter tuning, we adopt the leave-one-out cross validation strategy which is ideal when there are no explicit likelihood expressions. This allows us to operate without making any prior assumptions about the distribution of the quantity of interest. We started by testing them on a synthetic ODE problem where we make several important observations about both models, see Section \ref{sec: ODE}. In summary, we observed that there is a trade-off between accuracy versus computational time. Co-Kriging offers accurate predictions, see Table \ref{tab: Performance1D}, with little information but suffers from significantly higher run-times, owing to numerous kernel derivative computations and large matrix inversion, as is evident from Table \ref{tab: Complexity} and Table \ref{tab: CPUtimes}. On the contrary, Lagrangian Kriging lacks in accuracy and loses the crucial interpolation property, a defining feature of Kriging methods but is much faster to implement. We briefly discuss a possible resolution to improve interpolation of Lagrangian Kriging via kernel hyperparameters that minimize the interpolation error, offering a noticeable improvement, as shown in Figure \ref{fig: ODEcomparison}d. %
\newline\newline
Beyond this, these results are further reaffirmed by experiments in two dimensions. Co-Kriging is successfully able to account for point-wise linear differential equations in the case of 2D scalar functions in Section \ref{sec: 2D scalar}. In the case of potential flow prediction, we observed spurious flow patterns when using simple Kriging, see Figure \ref{fig: flowwithrandomobs}, which we later attributed to the choice of lengthscale parameter in Figure \ref{fig: thetaflowreversal}. Fortunately, LOOCV based hyperparameter tuning served as quick resolution to this effect allowing us to conduct experiments on potential flow prediction in the presence of an obstacle. In this case, both co-Kriging and Lagrangian Kriging perform well, as seen in Section \ref{sec: Flow prediction} where we also provide co-Kriging uncertainty estimates (see Figure \ref{fig: uncertaintyfluid}). %
These experiments also reveal another surprising limitation of Lagrangian Kriging. In Section \ref{sec: 2D scalar}, we see that constraints that do not involve the function itself (the same order of derivative as the one to be predicted) are not accounted for by Lagrangian Kriging, a concise reiteration of Remark \ref{remark: LK}. We see this appear again in Section \ref{sec: Flow prediction} since Lagrangian Kriging cannot account for the continuity equation and requires a two-step diagnostic procedure, predicting at the obstacle boundary followed by interpolation everywhere else. %
\newline\newline
Moving on to some perspectives, it is desirable to introduce a more complete Lagrangian Kriging model that should be able to incorporate the dependence between prediction coefficients associated with different orders of derivatives of a function explicitly, by construction. To bring down the computational time of co-Kriging, \cite{High_Dimensiona_de_Roo_2021} presents a method to exploit the Kronecker product structure of the Gram matrices produced by stationary kernel derivatives, to accelerate the covariance matrix computation and its inversion. This is especially important in higher dimensions as the computational cost scales as a cubic with respect to dimension. Even though promising, the application is not entirely obvious for the co-Kriging model. %
\newline\newline  
An important aspect is that the scope of this work is not necessarily limited to the linear setting, unlike other related works using Gaussian processes such as \cite{Graepel2003234, Srkk2011, Raissi2017}. However, we should mention that the use of Gaussian processes for solving non-linear PDEs has been investigated in recent works, such as \cite{chen_solving_2021} where they enforce differential constraints at prediction. Although not pursued here, this represents a compelling alternative methodology. In this regard, extending Lagrangian Kriging to non-linear differential equations is relatively straightforward but an analytical resolution via Lagrange multipliers is not always guaranteed. In such a scenario, it is possible to utilize gradient-based numerical optimization schemes instead. Co-Kriging can also be extended to the non-linear setting provided that we can compute the covariance between non-linear transformations of random functions. This is a prominent advantage of distribution-free best linear prediction interpretation over the use of Gaussian processes which are usually constrained to linear problems. We postpone the aforementioned ideas of extension to non-linear problems as future work. 

\paragraph{Acknowledgement.}
\label{sec: acknowledgements}
This research was conducted with the support of the consortium in Applied \\ Mathematics CIROQUO, gathering partners in technological and academia in the development of advanced methods for Computer Experiments. This research was partly funded by a CIFRE grant (convention n°2024/0423) established between the ANRT and Stellantis. We would also like to extend our gratitude to the UMR 6158, CNRS lab Laboratory of Informatics, Modelling and Optimization of the Systems (LIMOS) and Institute Fayol - École des Mines de Saint-Étienne for providing the necessary computational equipment.

\bibliographystyle{plain}

\appendix
\section{Appendix}
\label{sec: appendix}
\subsection{Proofs}
\printProofs
\newcommand{\kernel}{k(\vecx,\vecx')}
\newcommand{\delxi}{\partial_{x_i}} 
\newcommand{\delxti}{\partial_{x'_i}}
\newcommand{\delmi}{\partial_{m_i}}
\newcommand{\negpow}[1]{(-1)^{\delta'_{#1 x}}}
\newcommand{\constant}{\frac{1}{\theta^2}}
\newcommand{\Neg}[1]{\mathrm{neg}~(#1)}
\subsection{Kernel derivative expressions}
\label{sec: kernelderiavtives}
Since we need to quantify the covariance between derivatives of random field, we provide some analytical expressions we used in our implementations. For the squared exponential kernel given in \eqref{eqn: sqexp}, we derive the general expression for the first four orders of derivatives. For ease of notation, we denote, $k(\vecx,\vecx'\vert \sigma^2, \theta) = k(\vecx,\vecx')$. We begin by stating these formulas in their full form and provide the mathematical working later. 

\paragraph{First-order derivative}
\begin{equation}
\delmi \kernel = (-1)^{\delta'_{mx}} ~ \frac{(x_i - x_i')}{\theta^2} ~ \kernel 
\end{equation}

\paragraph{Second-order derivative}
\begin{equation}
\partial^2_{m^1_i m^2_j} \kernel = (-1)^{\delta'_{m^1x} + \delta'_{m^2x}} \left( h_i h_j - \constant  \delta_{ij} \right)  \kernel
\end{equation}

\paragraph{Third-order derivative}
\begin{equation}
\partial^3_{m^1_i m^2_j m^3_l}k = (-1)^{(\sum_{z=1}^3 \delta'_{m^zx})}  \big( h_ih_jh_l  
- \constant \left(h_i \delta_{jl} +    
h_j \delta_{il} + h_l \delta_{ij} \right) \big) \kernel
\end{equation}

\paragraph{Fourth-order derivative}
\begin{equation*} 					
\end{equation*}\\[-4em]					%
\begin{eqnarray}
\partial^4_{m^1_i m^2_j m^3_l m^4_o} k &=& (-1)^{(\sum_{z=1}^4 \delta'_{m^zx})}  \Big[ h_ih_jh_lh_o \nonumber\\
&-& \frac{1}{\theta^2} \big( h_i h_j \delta_{lo} + h_i h_l \delta_{jo}
+ h_i h_o \delta_{jl} 
+ h_j h_l \delta_{io} + h_j h_o \delta_{il} + h_l h_o \delta_{ij} \big) \nonumber\\
&+& \frac{1}{\theta^4} \big( \delta_{ij} \delta_{lo} + \delta_{il} \delta_{jo} + \delta_{io} \delta_{jl} \big)
\Big] \kernel
\end{eqnarray}
where $m^{z} \in \{x, x'\}$, $i, j, l, o \in \{1,\ldots,d\}$, $h_i \coloneqq \frac{(x_i - x_i')}{\theta^2}$, $\delta'_{m^zx} = 1$ if $m^z = x$ and $0$ otherwise and $\delta_{ij}$ is the usual Kronecker-delta. We will re-explain the notation and why we need it in the mathematical working that follows. 
\newline\newline
\noindent
An important observation is that for the first order derivatives, there exists an eigenfunction relationship that we can exploit,
$$
\delxi \kernel = - \frac{(x_i - x_i')}{\theta^2} ~\kernel
$$
for $i \in \{1,\ldots,d\}$  where $\vecx, \vecx' \in \R^d$. If instead this was with respect to the second variable, 
$$
\delxti \kernel =  \frac{(x_i - x_i')}{\theta^2} ~\kernel
$$
So, the general \textbf{first-order derivative} is basically, 
\begin{equation}
\label{eqn: firstorder}
\delmi \kernel = (-1)^{\delta'_{mx}} ~ \frac{(x_i - x_i')}{\theta^2} ~ \kernel
\end{equation}
where $\delta'_{mx} = 1$ if $m$ is the first variable $x$ and $\delta'_{mx} = 0$ otherwise, i.e., $m$ is the second variable $x'$. Since the expressions are about get a lot more complicated, we introduce a shorthand $h_i \coloneqq \frac{(x_i - x_i')}{\theta^2}$. For the \textbf{second-order derivative} consider, 
\begin{equation*}
\partial_{m^1_i} \left( \frac{\partial_{m^2_j} k}{k} \right) = \partial_{m^1_i} \left( \negpow{m^2} h_j \right)
\end{equation*}
This can be further simplified by using the division rule, 
\begin{eqnarray*}
\frac{k ~\partial^2_{m^1_i m^2_j}k -  \partial_{m^2_j} k ~\partial_{m^1_i} k }{k^2} =  - \delta_{ij} \constant  \negpow{m^2} \negpow{m^1} \\
\partial^2_{m^1_i m^2_j} \kernel = (-1)^{\delta'_{m^1x} + \delta'_{m^2x}} \left( h_i h_j - \constant  \delta_{ij} \right)  \kernel
\end{eqnarray*} 
where $m^1,m^2 \in \{x, x'\}$ and $i, j \in \{1,\ldots,d\}$ and we used $-\partial_{m^1_i} h_j = (-1)^{\delta'_{m^1x}}~\frac{\delta_{ij}}{\theta^2}$ which is straightforward to check. In the case of the \textbf{third-order derivative}, the expressions become more complicated. To be able to represent them meaningfully, we represent the recurring quantity $(-1)^{\delta'_{m^1x} + \ldots + \delta'_{m^ix}}$ as $\mathrm{neg}~(m^1:m^i)$. We proceed in a similar fashion, 
\begin{eqnarray*}
\partial_{m^1_i} \left( \frac{\partial^2_{m^2_j m^3_l} k}{k} \right) &=& \Neg{m^2,m^3} ~ \partial_{m^1_i} \left( h_j h_l \right) \\
\frac{k ~\partial^3_{m^1_i m^2_j m^3_l}k -  \left(\partial^2_{m^2_j m^3_l} k \right) ~(\partial_{m^1_i} k) }{k^2} &=& \Neg{m^2,m^3} ~ \left(h_l \partial_{m^1_i} h_j + h_j \partial_{m^1_i} h_l\right) \\
\frac{k ~\partial^3_{m^1_i m^2_j m^3_l}k -  \Neg{m^1:m^3} \left( h_ih_jh_l - \frac{h_i \delta_{jl}}{\theta^2} \right)k^2}{k^2} &=&  -\Neg{m^1:m^3} \left(\frac{h_l \delta_{ij}}{\theta^2} + \frac{h_j \delta_{il}}{\theta^2}\right) \\ 
\partial^3_{m^1_i m^2_j m^3_l}k &=&  \Neg{m^1:m^3} \big( h_ih_jh_l  \\ 
&-& \constant \left(h_i \delta_{jl} +    
h_j \delta_{il} + h_l \delta_{ij} \right) \big) \kernel
\end{eqnarray*}
We can continue on like this for the \textbf{fourth-order derivative} as well. The expression for the fourth-order is given by, 
\begin{eqnarray*}
\partial^4_{m^1_i m^2_j m^3_l m^4_o} k &=& \Neg{m^1:m^4} \left[ h_ih_jh_lh_o - \frac{1}{\theta^2} \sum_{\bm{c} \in \mathcal{C}}  h_{c_1} h_{c_2} \delta_{c_3c_4} + \frac{1}{\theta^4} \sum_{\bm{p} \in \mathcal{P}} \delta_{p_1p_2} \delta_{p_3p_4}
\right] \kernel
\end{eqnarray*}
where the set $\mathcal{C}$ represents the set of distinct pairs $\bm{c} = \{c_1, c_2\}$ (and $\{c_3, c_4\} \in \mathcal{C}$ as $\bm{c}^c$) chosen sequentially, $c_1 \in \{i, j, l, o\} = \mathcal{I}$ and $c_2 \in \mathcal{I} - \{c_1\}$ and $\mathcal{P}$ is the set of partitions of $\mathcal{I}$ into disjoint pairs, $\{\{\bm{p}, \bm{p}^c\} ~\big\vert~ \bm{p} \in \mathcal{I}, |\bm{p}| = 2\}$ with $\{p_1, p_2\} = \bm{p}$ and $\{p_3, p_4\} = \bm{p}^c$. The sets are of sizes, $\vert \mathcal{C} \vert = \binom{4}{2} = 6$ elements and $\vert \mathcal{P} \vert = 3$. 
\newline\newline
\noindent
There are clear patterns that emerge as we evaluate higher derivatives and might allow us to evaluate the general $n^{th}$ order derivative. The prescribed algorithm of differentiation can also be applied to other stationary kernels such as the Mat\`ern kernel and the main difference would show in the formula of $h_i$ and its derivatives.       

\subsection{Lagrangian Kriging working example.}
\label{sec: WorkingExampleLKmath}
To solve the Lagrangian Kriging problem with one observation and at one prediction point as presented in Examples \ref{example: LKWorkingExample}, we introduce a minor change in notation, $\alpha_1 = \alpha(\vecs\etoile), \alpha_2 = \alpha(\vecs^c)$. With this, at first, we expand the Kriging predictive error in \eqref{eqn: ConstrainedKrigOpt}, $\Delta(\vecs\etoile) + \Delta(\vecs^c)$, which can be expanded using \eqref{eqn: krigingopttraceform}, 
$$
\alpha_1^2 \Var{Z(\vecs^1)} - 2 \Cov{Z(\vecs^1), Z(\vecs\etoile)} \alpha_1 + \alpha_2^2 \Var{Z(\vecs^1} - 2 \Cov{Z(\vecs^1), Z(\vecs^c)} \alpha_2
$$
ignoring the terms not involving the Kriging coefficients, i.e., $\sigma^2(\vecs\etoile)$ and $\sigma^2(\vecs^c)$. We can solve the constrained minimization problem by constructing the Lagrangian with multiplier $2\lambda \in \R$,  
$$ 
\underset{\alpha_1, \alpha_2 \in \R}{\min} \cbind{\Delta(\vecs\etoile) + \Delta(\vecs^c) - 2\lambda ~(\alpha_1 + \alpha_2)} 
$$
and setting the gradient equal to zero with necessary optimality conditions gives, 
$$
\accolade{
\alpha_1 \Var{Z(\vecs^1)} - \Cov{Z(\vecs^1), Z(\vecs\etoile)} - \lambda = 0 \\
 \alpha_2 \Var{Z(\vecs^1)} - \Cov{Z(\vecs^1), Z(\vecs^c)} - \lambda = 0 \\ 
 \alpha_1 + \alpha_2 = 0 
}
$$
solving for the coefficients $\alpha_1, \alpha_2$ and the Lagrange multiplier $\lambda$ is a fairly simple exercise, and on doing so we get,
$$
\lambda = -\left(\frac{\Cov{Z(\vecs^1), Z(\vecs\etoile)} + \Cov{Z(\vecs^1), Z(\vecs^c)}}{2} \right)
$$ 
and consequently the prediction $\Zlk(\vecs\etoile) = \alpha_1 Z(\vecs^1)$ is, 
$$
\Zlk(\vecs\etoile) = \left(\frac{\Cov{Z(\vecs^1), Z(\vecs\etoile)} - \Cov{Z(\vecs^1), Z(\vecs^c)}}{2 \Var{Z(\vecs^1)}} \right) Z(\vecs^1)
$$

\end{document}